\documentclass[nonblindrev,dvipsnames]{informs4}
\OneAndAHalfSpacedXI%
\renewcommand{\arraystretch}{0.8}
\usepackage{mathtools}
\usepackage{bm}
\usepackage{scalerel}
\usepackage{setspace}
\usepackage{subfigure}
\usepackage{capt-of}
\usepackage{graphicx}
\usepackage{makecell, multirow}
\usepackage{tablefootnote}
\usepackage{fix-cm}
\usepackage{etoolbox}

\usepackage[shortlabels]{enumitem}
\usepackage{bbm}

\usepackage{booktabs}
\usepackage{float}
\usepackage{hyperref}
\hypersetup{
  hidelinks,
  pdftitle={A Markovian Approach for Cross-Category Complementarity in Choice Modeling},
  pdfauthor={Omar El Housni, Shuo Sun, Rajan Udwani}
}
\usepackage{cleveref}
\DeclareFontShape{OT1}{ptm}{m}{scit}{<->ssub * ptm/m/sc}{}
\makeatother

\usepackage{natbib}
 \bibpunct[, ]{(}{)}{,}{a}{}{,}%
 \def\bibfont{\small}%

\TheoremsNumberedThrough     %
\ECRepeatTheorems

\EquationsNumberedThrough    %

\usepackage[ruled,vlined, noend]{algorithm2e}
\usepackage{tikz}
\usetikzlibrary{positioning, fit, shapes, arrows.meta, calc}
\usepackage{xcolor}

\newcommand{\NN}{\mathbf{N}}
\newcommand{\EE}{\mathbb{E}}
\newcommand{\MM}{\mathbf{M}}
\newcommand{\bfE}{\mathbf{E}}

\newcommand{\ALG}{\operatorname{ALG}}
\newcommand{\bfv}{\mathbf{v}}
\newcommand{\bfpsi}{\boldsymbol{\psi}}
\newcommand{\bfr}{\mathbf{r}}

\newcommand{\bfl}{\boldsymbol{\lambda}}
\newcommand{\bfrho}{\boldsymbol{\rho}}
\newcommand{\1}{\mathbbm{1}}

\makeatletter
\renewcommand{\footnoterule}{%
  \kern 10pt
  \hrule \@width \columnwidth
  \kern 0.6pt
}
\makeatother



\makeatletter
\def\Dskips#1#2#3#4{\abovedisplayskip 0.5em \@plus#2em \@minus#3em
   \abovedisplayshortskip \z@ \@plus#2em
   \belowdisplayshortskip 0.4em \@plus#2em \@minus#2em
   \belowdisplayskip \abovedisplayskip}%
\makeatother
\AtBeginDocument{\normalsize}
\setlist{itemsep=1pt, topsep=3pt, parsep=1pt}
\SetAlFnt{\smallXnoCE}

\begin{document}

\TITLE{A Markovian Approach for Cross-Category Complementarity in Choice Modeling}

\ARTICLEAUTHORS{%
\AUTHOR{Omar El Housni}
\AFF{School of Operations Research and Information Engineering, Cornell Tech, Cornell University}
\AUTHOR{Shuo Sun}
\AFF{Department of Industrial and Systems Engineering, University of Wisconsin--Madison}
\AUTHOR{Rajan Udwani}
\AFF{Department of Industrial Engineering and Operations Research, University of California, Berkeley}
}

\ABSTRACT{%
In modern retail and e-commerce, customers often purchase multiple items across distinct product categories, exhibiting both substitution and complementarity. We consider the cross-category assortment optimization problem where retailers jointly determine assortments across categories to maximize expected revenue. Most prior work on the topic either overlooks complementarity or proposes models that lead to intractable optimization problems, despite being based on the multinomial logit (MNL) choice model. We propose a sequential multi-purchase choice model for cross-category choice that incorporates complementarity through a Markovian transition structure across categories, while allowing general Random Utility Maximization (RUM)-based choice models to capture the within-category substitution. We develop an Expectation-Maximization algorithm for estimation, and a polynomial-time algorithm for unconstrained assortment optimization that yields the optimal solution when the within-category substitution follows a Markov chain choice model.
Furthermore, we introduce an empirical metric to quantify the strength of complementarity across product categories and conduct extensive numerical experiments using both synthetic data and a large-scale transaction-level dataset from a major U.S. grocery store. Our model improves predictive accuracy and model fit in both synthetic and real-data experiments and improves expected revenue in synthetic settings with complementarity. It also reveals intuitive market structures such as brand-loyal cross-category purchasing.
Overall, we believe that our model provides a theoretically grounded and practical framework for modeling complementarity and making better cross-category assortment decisions.
}%

\maketitle
\section{Introduction}\label{sec:Intro}

Customers frequently purchase products from multiple categories within a single transaction. In retail settings, for instance, a shopper may select products from several different categories. A large-scale study of U.S. and U.K. supermarkets found that $83\%$ of transactions involved products from more than one category \citep{tanusondjaja2016understanding}. This phenomenon extends beyond retail to service industries such as airlines, where a traveler first books a flight and then makes a sequence of choices among ancillary services. When purchases span multiple categories, a choice in one category often changes a customer’s preferences in another. The marketing and economics literature provides strong empirical evidence of these cross-category dependencies, commonly referred to as \emph{complementarity} effects \citep{mulhern1991implicit,Aurier2014,KWAK201519}.

In many settings, customers’ choices unfold sequentially across related product categories. Shoppers typically begin by selecting a focal or primary product, such as a cell phone, television, cake mix, or flight, which then shapes their choices among complementary products, such as a phone case, television bracket, cake frosting, or ancillary service. This directional dependence is often referred to in the marketing literature as asymmetric complementarity \citep{manchanda1999shopping,lee2013direct}. For instance, the price of cake mix has been shown to substantially affect frosting sales, whereas the reverse effect is considerably weaker \citep{mulhern1991implicit,manchanda1999shopping}. Similar one-way relationships have been documented between spaghetti and spaghetti sauce \citep{walters1991assessing}, as well as between laundry detergent and fabric softener \citep{manchanda1999shopping}. Asymmetric complementarity is also prominent in digital ecosystems, where the purchase of an e-reader drives subsequent e-book sales \citep{li2019intertemporal}. Such patterns may arise from the logic of product usage, in-store merchandising strategies, or consumers’ cognitive decision-making processes. Moreover, customers often have a higher willingness to pay for primary products than for complementary products, causing the influence of the primary product to dominate the reverse effect \citep{ke2022cross}.

Optimizing assortments across interdependent categories poses a significant managerial challenge in many industries. In grocery retail, for instance, jointly managing assortments of primary products such as cake mixes and complementary products such as frosting can increase both sales and profitability. In e-commerce, coordinating assortments of core products and compatible accessories at the fulfillment-center level can improve customer experience and facilitate cross-selling. For example, offering phones together with compatible chargers or protective cases can increase the likelihood of multi-product purchases. Similarly, auto dealerships must decide not only which vehicle models and trims to stock, but also how to design a profitable portfolio of accessories and extended service plans. In the travel industry, assortment decisions for airline routes and fare classes are closely linked to the design of ancillary offerings, such as seat upgrades and partner hotel bookings, which can contribute substantially to profitability.

Despite its practical importance and the long tradition of studying complementarity in economics and marketing, cross-category complementarity remains underexplored in the assortment optimization literature. Most existing work focuses on substitution within a single category, typically using variants of the multinomial logit (MNL) model \citep{talluri2004revenue}, nested logit model \citep{Davis2014}, or Markov chain (MC) choice model \citep{blanchet2016markov}. However, ignoring cross-category dependencies can lead to suboptimal assortments and substantial profit losses \citep{ghoniem2016optimizing}. To date, no standard framework exists for modeling and optimizing assortments in general cross-category settings. Existing models that incorporate complementarity often rely on heuristics without theoretical guarantees or become computationally intractable as the number of products and categories grows \citep{ghoniem2016optimizing,ke2022cross,chen2022assortmentoptimizationmultivariatemnl,tulabandula2023}. Moreover, most prior studies use the MNL model to represent substitution within each category, even though the broader discrete choice literature offers richer alternatives that may be more suitable in practice. Finally, prior work has largely emphasized algorithm design under a given model, with less attention to behavioral realism and ease of estimation.

In this paper, we address these gaps by proposing a practical yet theoretically grounded model that captures cross-category complementarity, remains tractable for large-scale assortment optimization, and supports efficient parameter estimation. Its modular structure also allows us to incorporate richer choice models for within-category substitution.

\subsection{Summary of Our Results}

Our contributions are threefold. First, we introduce a general and estimable multi-purchase choice model that captures cross-category complementarity while remaining tractable for large-scale assortment optimization. Second, we establish the tractability of joint assortment optimization across multiple categories when within-category substitution follows an MC choice model. {Third, we validate the model using both synthetic experiments and large-scale retail data, demonstrating improved predictive performance and actionable managerial insights into the strength and heterogeneity of cross-category complementarity. The synthetic experiments also demonstrate meaningful revenue gains from modeling complementarity.}
\smallskip

\noindent \textbf{Model.} At the core of our work is a novel Markovian multi-purchase choice model in which customers make sequential choices across multiple product categories. To build intuition, we first describe the model for two categories: a primary category and a secondary category. Within each category, choices are governed by a general random-utility-based model, such as an MC model or the simpler MNL model. After the customer selects a product in the primary category, this choice induces a product-specific distribution over latent initial preference states in the secondary category.  Each latent state corresponds to the product that anchors the customer’s initial preference in the secondary category, rather than a product that the customer necessarily observes or encounters. If this product is offered, the customer selects it directly. Otherwise, the customer substitutes among the available alternatives according to the within-category choice model, conditional on the latent initial preference. Thus, the upstream choice determines the entry point of the downstream substitution process, while the subsequent substitution mechanism is shared across upstream products. This shared substitution mechanism does not eliminate the influence of the upstream choice. Different upstream products induce different distributions over latent initial states, and, when the downstream substitution model is sufficiently expressive, these differences propagate into  heterogeneous downstream choice probabilities.

This approach differs from prior work, which typically extends single-purchase MNL models by assuming that the primary choice redefines the entire set of preference weights in the secondary category. Instead, we retain a fixed within-category choice model for the secondary category, with the primary choice determining the latent entry point of the downstream substitution process rather than reshaping the entire preference structure. Our framework also extends naturally to settings in which categories are connected through a directed acyclic graph (DAG); see Section~\ref{sec:dag}.  This separation allows the downstream substitution model to be substantially richer than MNL and avoids the restrictive independence of irrelevant alternatives (IIA) property of the MNL model.
\smallskip

\noindent \textbf{Assortment optimization.}
A key advantage of our approach is that it enables tractable unconstrained assortment optimization under the MC choice model, even with an arbitrary number of categories. At first glance, the problem appears challenging because assortments must be selected jointly across categories and the resulting decisions are inherently interdependent. Tractability arises because, after appropriate adjustments to product revenues, the optimization problem decouples across categories. In particular, the choice in one category affects only the initial arrival probabilities in the next category’s MC model, and an optimal assortment for an MC model is invariant to these initial probabilities. Consequently, the overall problem can be solved sequentially, category by category, following an appropriate ordering of the categories.

This scalability represents a substantial improvement over existing multi-purchase choice models, which are often restricted to two categories. In contrast, our model accommodates an arbitrary number of categories without sacrificing tractability. Single-category assortment optimization under a general random utility maximization (RUM) model does not admit a polynomial-time constant-factor approximation algorithm \citep{aouad2018approximability}. We therefore focus on the MC model, which is both expressive and tractable for single-category assortment optimization. Other within-category models with tractable single-category assortment problems, such as the nested logit model, could also be considered. However, whether joint optimization remains tractable under such non-Markovian choice models \citep{selena,ma2023assortment} is an open and challenging direction for future research.
\smallskip

\noindent \textbf{Empirical validation and managerial insights.} %
To assess the practical value of our framework, we examine whether it can recover meaningful cross-category complementarity from transaction data and whether modeling this dependence improves operational decisions. Treating each category independently would simplify both estimation and optimization, but it would overlook cross-category dependencies that may be important in practice. We therefore develop an efficient expectation-maximization (EM) algorithm and evaluate the model using both real and synthetic data.

Our empirical analysis is based on a dataset of nearly six million item-level purchase records from a major U.S. grocery retailer. We first introduce a new metric for quantifying cross-category complementarity and use it to identify product categories with strong and heterogeneous complementarity dependencies.
For these categories, our model improves model fit and predictive accuracy relative to both benchmarks that ignore cross-category dependence and previous MNL-based models of complementarity. The estimated parameters also reveal interpretable product-level relationships. For example, customers who purchase a particular brand of cake mix exhibit a substantially stronger preference for the corresponding brand of frosting. Such estimates can inform operational decisions beyond assortment optimization, including coordinated pricing, promotions, and store layout.

We complement this empirical analysis with synthetic experiments in which we systematically vary the strength of cross-category complementarity. These experiments allow us to evaluate both predictive performance and the revenue impact of the resulting assortment decisions. Results show that
our model consistently outperforms benchmarks in both predictive performance and expected revenue. Importantly, the revenue gains increase with the strength of complementarity: in settings with
strong cross-category dependence, our model delivers $4\%$--$7\%$ higher expected revenue relative to independent choice models that ignore complementarity. In particular, the comparison with \textsc{IndMC}, which retains MC within-category substitution but omits cross-category dependence, shows the value of combining Markovian cross-category transitions with sufficiently expressive within-category substitution, as discussed in our model introduction.

\emph{Outline. }
The remainder of the paper is organized as follows. Section~\ref{sec:model} introduces the choice model for two categories and extends it to a general DAG. Section~\ref{sec:3-opt} studies assortment optimization, establishing tractability for the unconstrained problem and characterizing the effect of cardinality constraints.  Section~\ref{sec:estimation} develops the EM estimation procedure and discusses convergence and identifiability. Sections~\ref{sec:synthetic_experiments} and~\ref{sec:real_data} evaluate the model using synthetic experiments and real-world retail data, respectively. Sections~\ref{sec:cardinality_experiments} and~\ref{sec:bidirectional} present computational results for cardinality-constrained assortment optimization and a bidirectional extension for parallel categories. Section~\ref{sec:conclusion} concludes with a summary of findings and directions for
future research.

\subsection{Related Work}\label{sec:related}

The literature on single-purchase assortment optimization is extensive. Efficient algorithms have been developed under various choice models, such as MNL \citep{talluri2004revenue, gallego2004managing,rusmevichientong2010dynamic}, mixed MNL \citep{Bront2009}, nested logit \citep{Davis2014}, and Markov chain choice models \citep{blanchet2016markov,Feldman2017, desir2020constrained, udwani2025submodular}. For comprehensive surveys, we refer the reader to \cite{kok2015assortment}, \cite{strauss2018}, and \cite{heger2024assortment}.
In this section, we focus on reviewing related work on complementarity modeling and assortment optimization with complementarity. We also review multi-purchase choice models and sequential choice models that only focus on substitution effects.

\noindent \textbf{Complementarity in marketing and economics. }
Modeling product complementarity has a long history in marketing and economics, where the literature focuses on estimation and interpretation of model parameters rather than optimization.
 This stream of work typically classifies complementarity as either \textit{symmetric} or \textit{asymmetric}, depending on whether purchase order influences the complementarity effect \citep{seetharaman2005models,lee2013direct,Aurier2014, KWAK201519}. For instance, \cite{lee2013direct} propose a direct utility model with a latent decision sequence across categories to measure the asymmetric complementarity effects. \cite{ruiz2019shopper} develop a sequential probabilistic model of customers' shopping baskets, where customers select products sequentially and each choice is conditional on products already in the basket. However, these models are primarily designed for prediction and estimation, and do not readily extend to assortment optimization.

\noindent \textbf{Choice models and assortment optimization with product complementarity. } Our work is most closely related to recent research on assortment optimization under asymmetric complementarity, which typically assumes a sequential purchase structure where an initial choice influences subsequent ones. The closest paper is \cite{ke2022cross}, which models a two-stage process: customers first select a primary product using an MNL model, then make a secondary choice from another category using an MNL model whose parameters depend on the first selection. The authors consider two settings: a personalized one, where the secondary assortment can be tailored to the initial purchase, and a non-personalized one, where the same secondary assortment is offered to all customers. The personalized setting is straightforward to solve in both \cite{ke2022cross} and our framework (see Section~\ref{sec:optimization} for details). We focus instead on the non-personalized setting, which is more relevant for retailers that cannot adjust assortments in real time, such as brick-and-mortar stores, and is also more technically challenging. In this setting, \cite{ke2022cross} show that the unconstrained assortment problem is NP-hard even for two categories. In contrast, our model admits a tractable algorithm for the unconstrained case, accommodates more general within-category choice models (e.g., the MC model), and extends to an arbitrary number of categories. Other works have explored the joint assortment and pricing problem for complementary categories, proposing mixed-integer programming formulations \citep{ghoniem2016optimizing} or stylized models \citep{rodriguez2011assortment}.

A separate line models asymmetric complementarity between products in a single-purchase setting \citep{feng2018substitutability,lo2019}, whereas in our model complementarity arises from actual purchases across categories.
A related stream considers personalized recommendations conditioned on a given primary choice \citep{chen2024assortment,ban2024selling}. A key distinction is that in our model, the primary product is not exogenously given but is part of the endogenous customer choice process. Moreover, we focus on determining fixed assortments offered to all customers, rather than an assortment that adapts dynamically to their prior selections.

For symmetric complementarity, recent research considers assortment optimization under variants of multivariate MNL (MVMNL) models. In these models, substitution and complementarity are captured through a `bundle’ utility that combines individual product utilities with pairwise interaction terms \citep{chen2022assortmentoptimizationmultivariatemnl,Jasin2024,tulabandula2023}.
\cite{chen2022assortmentoptimizationmultivariatemnl} study a two-category problem and propose a $0.74$-approximation algorithm for the unconstrained case, but show that the capacitated version admits no constant-factor approximation (under the Exponential Time Hypothesis). \cite{Jasin2024} develop a Fully Polynomial Time Approximation Scheme (FPTAS) for a variant with group-based interactions for both uncapacitated and capacitated assortment problems, allowing an arbitrary number of groups. \cite{tulabandula2023} examine another MVMNL variant where customers can buy at most $K$ products and show the unconstrained problem is NP-hard even when $K=2$.

Currently, all models incorporating complementarity rely on MNL-based structures. In contrast, our model accommodates more general within-category choice models, such as the MC model. Moreover, whereas prior work establishes NP-hardness even in the unconstrained case, we show that our framework admits tractable unconstrained optimization when within-category choices follow an MC model.

\noindent \textbf{Multi-purchase and sequential choice models with substitution effects only. }
A separate stream of literature studies multi-purchase choice behavior driven by only substitution. These models typically govern basket composition by assuming customers select a predetermined or random number of top-utility products \citep{benson2018, bai2024multi, abdallah2024multi}, purchase any item exceeding a utility threshold \citep{gallego2019threshold}, or maximize a value function over bundles that encodes substitutability \citep{immorlica2021combinatorial}. Because these frameworks exclude complementarity effects, their assortment optimization approaches fundamentally differ from ours. Finally, there is a stream of work that considers choice models where customers view the assortments in multiple stages sequentially \citep{flores2019assortment, liu2020assortment, feldman2021assortment, gao2021assortment}. In these studies, customer choices in each stage follow an MNL model. This setting is clearly distinct from ours, as it assumes a single purchase, captures only substitution effects, and treats product preferences in each stage as independent.
\section{Model}\label{sec:model}
In this section, we formally introduce our multi-category choice model. We begin by defining the model for two categories and then extend it to more general, multi-category structures.

\subsection{Choice Model with Two Categories}\label{sec:model-two} At a high level, our model uses a general ranking-based choice model (RCM) to capture substitution effects within each category, and we capture the complementarity effect between categories through a Markovian transition structure.
Recall that in RCM, customer preferences are represented by a probability distribution $\mathbf{p}$ over a set $\boldsymbol{\Sigma}$ of all rankings over the universal set $\NN^+ = \NN \cup \{0\}$. Each ranking $\sigma \in \boldsymbol{\Sigma}$ is a permutation of $\NN^+$, where a lower value corresponds to a higher preference. That is, $\sigma(i) < \sigma(j)$ implies product $i$ is preferred over product $j$. Given an offered assortment $S \subseteq \NN$, a customer draws a ranking $\sigma$ according to $(p_{\sigma})_{\sigma \in \boldsymbol{\Sigma}}$ and selects the most preferred product. The choice probability of product $i \in S$ is:
 \(\phi(i,S)=\sum_{\sigma \in \boldsymbol{\Sigma}} \1\left\{i= \arg \min_{j\in S\cup\{0\}}\sigma(j)\right\} p_{\sigma}\).
It is known that RCM is equivalent to RUM-based choice models \citep{block1959random,farias2013nonparametric}.

Let $A$ and $B$ denote two product categories, with ground sets $\NN_A = \{1,\dots,n_A\}$ and $\NN_B = \{1,\dots,n_B\}$, respectively. Each category contains a no-purchase option, denoted $0_A$ for $A$ and $0_B$ for $B$. We refer to the full set of options as $\NN_A^+=\NN_A\cup\{0_A\}$  and $\NN_B^+=\NN_B\cup\{0_B\}$. For assortments $S_A \subseteq \NN_A$ and $S_B \subseteq \NN_B$, we define the extended assortments $S_A^+ = S_A \cup \{0_A\}$ and $S_B^+ = S_B \cup \{0_B\}$. We use $\phi_A$ and $\phi_B$ to denote the RCM models for $A$ and $B$, respectively.
\smallskip

\begin{figure}[htbp]
\vspace{-5pt}
\centering
\begin{tikzpicture}[scale=0.58, transform shape,
    node distance=1cm and 1cm,
    category/.style={
        rectangle, draw, thick, rounded corners,
        minimum width=4.5cm, minimum height=6.2cm,
        label={[font=\bfseries]above:#1}
    },
    assortment_a/.style={
        rectangle, draw, thick, fill=blue!10,
        minimum width=3.5cm, minimum height=3.5cm,
        label={[anchor=north west]north west:$S_{A}$}
    },
    assortment_b/.style={
        rectangle, draw, thick, fill=blue!10,
        minimum width=2.5cm, minimum height=2.5cm,
        label={[anchor=north west]north west:$S_{B}$}
    },
    product/.style={
        ellipse, draw, thick, fill=white, minimum size=0.9cm
    },
    no_purchase/.style={
        ellipse, draw, thick, fill=gray!20, minimum size=0.9cm
    },
    customer/.style={
        circle, draw, thick, minimum size=0.5cm,
        append after command={
            \pgfextra{
                \draw[thick] (\tikzlastnode.south) -- ++(0,-0.5);
                \draw[thick] (\tikzlastnode.south) ++(-0.25,-0.25) -- ++(0.5,0);
                \draw[thick] (\tikzlastnode.south) ++(0,-0.5) -- ++(-0.25,-0.25);
                \draw[thick] (\tikzlastnode.south) ++(0,-0.5) -- ++(0.25,-0.25);
            }
        }
    },
    main_arrow/.style={
        -Stealth, thick, red
    },
    sub_arrow/.style={
        -Stealth, thick, red, densely dotted
    }
]
\usetikzlibrary{calc}

\node[category=$\NN_A^+$] (cat_a) {};
\node[assortment_a] (S_a) at (cat_a.center) {};

\node[product] (i) at ($(S_a.center) + (0, -0.2cm)$) {$i$};
\node[product] (1) [above=0.4cm of i] {$k$};

\node[no_purchase] (0a) at ($(S_a.south)!0.5!(cat_a.south)$) {$0_A$};

\node[customer,
      label={[yshift=-1.1cm]below:Customer}]
      (start) at ($ (cat_a.north west) + (-2.2cm, -2.8cm) $) {};

\coordinate (leftHand)  at ($(start.south) + (-0.25cm,-0.25cm)$);
\coordinate (rightHand) at ($(start.south) + ( 0.25cm,-0.25cm)$);
\foreach \x in {1,i,0a} {
    \draw[-Stealth, thick] (rightHand) -- (\x.west);
}

\node[draw=blue, very thick, circle, fit=(i), inner sep=-3pt] (choice_i) {};

\node[category=$\NN_B^+$] (cat_b) [right=of cat_a, node distance=6.5cm] {};
\node[assortment_b] (S_b) at (cat_b.center) {};

\node[product] (m) at (S_b.center) {$j$};

\node[product] (j) at ($(S_b.north) + (0cm, 0.6cm)$) {$\ell$};

\node[no_purchase] (0b) at ($(S_b.south)!0.5!(cat_b.south)$) {$0_B$};
\node[draw=blue, very thick, circle, fit=(j), inner sep=-3pt] (choice_i) {};

\draw[main_arrow, bend left=20] (i.east) to node[above, xshift=6pt,  yshift=5pt] {$\lambda_{i,\ell}$} (j.west);
\draw[main_arrow, bend right=25] (i.east) to node[above, xshift=-1pt,  yshift=0pt] {$\lambda_{i,j}$} (m.west);
\draw[main_arrow, bend right=40] (i.east) to node[above, xshift=17pt,  yshift=-7pt] {$\lambda_{i,0_B}$} (0b.west);

\draw[sub_arrow, bend left=40] (j.east) to node[right, yshift=-2pt] {} (m.east);
\draw[sub_arrow, bend left=50] (j.east) to node[right, yshift=-6pt, xshift=-2pt] {} (0b.east);

\end{tikzpicture}
\caption{Diagram of the choice model with two categories.}
\label{fig:choice_model_diagram}
\vspace{-10pt}
\end{figure}
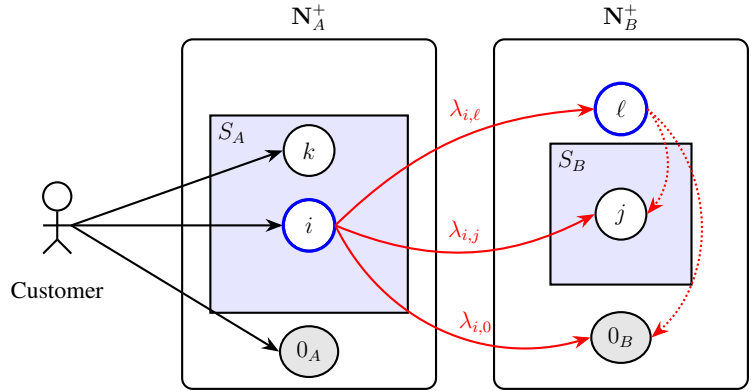

\noindent \textbf{Our Model.} Now, we introduce our two-category model, which is visualized in Figure~\ref{fig:choice_model_diagram}. The choice process proceeds as follows: the customer first selects $i \in S_A^+$ according to $\phi_A$. Unlike prior work that uses this choice to induce a new set of parameters for category~$B$, our model uses $i$ to determine an \emph{initial attraction} toward some $\ell \in \NN_B^+$. This initial attraction is governed by a vector of transition probabilities $\bfl^i = (\lambda_{i,\ell})_{\ell \in \NN_B^+}$, which satisfies $\sum_{\ell \in \NN_B^+} \lambda_{i,\ell} = 1$. The dependence of $\bfl^i$ on $i$ is the key feature that captures heterogeneous complementarity.
{
We interpret $\ell$ as a \emph{latent initial state} in category~$B$, rather than as an item that is necessarily displayed or observed by the customer. Note that we allow $\ell \notin S_B^+$, which captures the idea that a purchase in category~$A$ may induce a preference for a specific complementary product (e.g., a particular brand or color), even if that item is not available in the offered assortment. In this sense, $\ell$ determines the entry point of the downstream substitution process rather than an observed browsing action. As with single-category Markov chain choice models, this transition is not intended to represent the customer’s mental deliberation process; instead, it provides a parsimonious way to encode how upstream choices shape initial downstream preferences for all products \citep{blanchet2016markov}. The transition vector $\bfl^i$ is a preference distribution over the full product universe and is therefore assortment-independent; the offered assortment affects the realized downstream choice through the substitution model $\phi_B(\cdot\mid \ell)$.
}
Once the initial attraction is determined, the subsequent choice behavior in category~$B$ depends on the attracted product $\ell$. Thus, the upstream choice continues to influence downstream substitution through the entry point it induces, while the residual substitution process is shared across upstream choices. If the attracted product $\ell$ is available ($\ell \in S_B^+$), the customer selects it and exits. If $\ell \notin S_B^+$, the customer substitutes among the products in $S_B^+$ according to $\phi_B$, conditioned on $\ell$ being the top-ranked product in $\NN_B^+$. This Markovian structure ensures that, conditional on $\ell$, the final choice in category~$B$ is independent of $i$.

Formally, let $\boldsymbol{\Sigma}_\ell$ be the set of rankings of
$\NN_B^+$ that place $\ell$ first, and write
$z_\ell=\sum_{\sigma\in\boldsymbol{\Sigma}_\ell}p_\sigma$.
We specify a ranking distribution $q_\ell$ on
$\boldsymbol{\Sigma}_\ell$ for every initial state $\ell$. If
$z_\ell>0$, we use the usual conditional distribution,
$q_\ell(\sigma)=p_\sigma/z_\ell$. If $z_\ell=0$, the ranking
distribution $\mathbf{p}$ does not determine substitution from $\ell$.
For the unrestricted RCM, we complete the model by drawing a
ranking from $\mathbf{p}$ and moving $\ell$ to first place while preserving
the relative order of the other options; $q_\ell$ is the
distribution of the resulting ranking. Given an initial attraction to an unavailable product
$\ell\notin S_B^+$, the probability of choosing $j\in S_B^+$ is
\[
\phi_B(j,S_B\mid\ell)
=\sum_{\sigma\in\boldsymbol{\Sigma}_\ell}
q_\ell(\sigma)
\1\!\left\{j=\arg\min_{k\in S_B^+}\sigma(k)\right\}.
\]
Combining the direct purchase and substitution cases, the overall probability of choosing $i \in S_A^+$ and then $j \in S_B^+$ is
\begin{equation}\label{eq:2-prob}
P(i,j,S_A,S_B) = \phi_A(i,S_A)\bigl(\lambda_{i,j} + \sum_{\ell \in \NN_B \setminus S_B} \lambda_{i,\ell}\, \phi_B(j,S_B \mid \ell)\bigr).
\end{equation}
The first term corresponds to the customer being directly attracted to $j$, while the second term accounts for attraction to an unavailable product $\ell$, followed by substitution to $j$.

{
We emphasize that adopting a Markovian transition structure between categories~$A$ and~$B$ does not imply that the upstream choice in category~$A$ loses its influence on category~$B$. Rather, this influence is captured through the initial arrival product $\ell$, which determines the entry point of the downstream substitution process. Since $\phi_B(\cdot\mid \ell)$ is conditioned on $\ell$, substitution behavior can vary with the upstream choice through this entry point. To illustrate, consider category~$A$ as mobile phones and category~$B$ as phone cases. If a customer purchases a black iPhone in category~$A$, it is natural that they initially prefer a black iPhone case. If that specific case is unavailable, the Markovian structure implies that they substitute in a manner consistent with other customers who also begin with a preference for black iPhone cases. When the within-category choice model $\phi_B$ is sufficiently expressive, it can capture such structured substitution behavior. For instance, this customer may switch to an iPhone case of a different color but is unlikely to switch to a Samsung case. If $\phi_B$ is a Markov chain choice model, this structure can be captured by setting the transition probability from a black iPhone case to a Samsung case to zero. }

We also acknowledge that, relative to models that allow the downstream choice model to vary freely with the upstream choice (e.g., a separate $\phi_B^i$ for each upstream choice $i$), our framework imposes more structure: we assume $\phi_B$ is a fixed ranking-based model and the upstream choice affects category~$B$ only through the initial arrival product. However, as we demonstrate in later sections, this design captures realistic complementarity patterns while enabling tractable optimization. Furthermore, our comparison with existing models confirms that our formulation remains structurally distinct and is not a simpler special case of prior work (see Section \ref{sec:special-mnl} for details).

Observe that if $\bfl^i$ is the same for all $i$, then cross-category complementarity vanishes, and the two-category model reduces to two independent choice models. Note that ranking-based choice models include the MNL and MC models as special cases. While our framework allows for general RCMs to capture substitution, even the single-category unconstrained assortment problem under a general RCM is NP-hard to approximate \citep{aouad2018approximability}. Therefore, we focus on settings where the within-category choice models $\phi_A$ and $\phi_B$ are either MNL or MC, for which single-category assortment optimization is tractable. {We refer to the specifications in which both within-category models are MNL or both are MC as \textsc{MarkovMNL} and \textsc{MarkovMC}, respectively.} Our choice of these two particular models is further motivated by two factors. First, using the MNL model enables a direct comparison with the existing literature, which is predominantly MNL-based. Second, as we show, the structure of the MC model is critical for achieving tractable assortment optimization within our framework.

\subsubsection{Special Case: Multinomial Logit (MNL)}\label{sec:special-mnl}

We first analyze the special case where both $\phi_A$ and $\phi_B$ are MNL models, to establish a clear connection with the existing literature. In a standard single-purchase MNL formulation, each product $k$ has a preference weight $v_k$. The probability of choosing product $k$ from an assortment $S$ is
$
\phi(k,S) = \frac{v_k}{\sum_{i\in S} v_i + v_0},
$
where $v_k$ is the preference weight of product $k$ and $v_0$ is the weight of the no-purchase option (typically normalized to $1$). We use superscripts $A$ and $B$ to denote the weights for categories $A$ and $B$, respectively.

In our two-category model with $\phi_A$ and $\phi_B$ as MNL, the customer first selects product $i$ from $S_A^+$ with probability $\phi_A(i,S_A)=\frac{v_i^A}{\sum_{k\in S_A}v_k^A+1}$, following the standard MNL model. We then consider the conditional choice probability for product $j$ in $S_B^+$, given that $i$ was chosen in category $A$. With probability $\lambda_{i,j}$, the customer is directly attracted to $j$ and purchases it. Alternatively, with probability $\lambda_{i,k}$, the customer is initially attracted to a product $k \notin S_B$ and then substitutes according to $\phi_{B}(\cdot |k)$. In the case of the MNL
model, the Independence of Irrelevant Alternatives (IIA) property implies that this distribution is independent
of $k$: $\phi_B(j,S_B\mid k)=\phi_B(j,S_B)$ for every $k\in\NN_B\setminus S_B$ and $j\in S_B^+$. We show this formally in Lemma~\ref{prop:rcm-mnl} in Appendix~\ref{app:rcm-spec}. Thus, the customer substitutes according to the same MNL distribution $\phi_B$. Combining with \eqref{eq:2-prob}, we obtain the following lemma.
    \begin{lemma}\label{lem:choice_prob}
In the MNL special case, given assortments $S_A \subseteq \NN_A$ and $S_B \subseteq \NN_B$, the joint probability of selecting $i \in S_A^+$ in category $A$ and then $j \in S_B^+$ in category $B$ is \(
P(i,j,S_A,S_B) =
\frac{v_i^A}{V^A(S_A)+1}
\bigl(
\lambda_{i,j} +
\frac{v_j^B}{V^B(S_B)+1} \sum_{m \in \NN_B \setminus S_B} \lambda_{i,m}
\bigr)\),
where $V^A(S_A) = \sum_{k\in S_A} v_k^A$ and $V^B(S_B) = \sum_{k\in S_B} v_k^B$.
\end{lemma}
    {\bf Comparison to other MNL-based models.} Now we show our MNL-based model is structurally different from previous two-category models, which are all MNL-based. In \citet{ke2022cross}, the choice from category $A$ induces a separate conditional MNL on $B$, while in our model it affects only the initial attraction probabilities in $B$.
The multivariate MNL model of \citet{chen2022assortmentoptimizationmultivariatemnl} captures symmetric complementarity by assigning a preference weight to each bundle, with the choice probability proportional to this weight. In Appendix~\ref{app:connect} we show that none of these model classes contains another. Moreover, our category model is also distinct from the single-category model proposed in  \citet{cao2023revenue}, where the choice probability is a convex combination of an MNL and an independent demand component. While Lemma~\ref{lem:choice_prob} shows that the choice probability in $B$ also has a mixture form, our mixture proportion is endogenous to the assortment $S_B$ and the first-stage choice $i$, rather than a fixed parameter.

\subsubsection{Special Case: Markov Chain (MC)}\label{sec:special-mc} We first recall the definition of an MC model. A standard Markov chain choice model is defined by an initial arrival probability vector $\boldsymbol{\psi} = (\psi_i)_{i \in \NN^+}$ and a transition matrix $\bfrho=(\bfrho_{i,k})_{i,k\in \NN^+}$ \citep{blanchet2016markov}. The choice process begins with an initial attraction to product $i\in \NN^+$ with probability $\psi_i$. If $i\in S^+$, or if it is the no-purchase option, the process terminates and $i$ is chosen. Otherwise, if $i\notin S^+$, the customer transitions to a new product $k\in \NN^+$ with probability $\bfrho_{i,k}$, and this process repeats until an available product or the outside option is reached. For our two-category setting, $\phi_A$ and $\phi_B$ are parameterized by $(\bfpsi^A,\bfrho^A)$ and $(\bfpsi^B,\bfrho^B)$, respectively. In our two-category model, after a customer selects $i$ from $S_A^+$, they form an initial attraction to an item $j$ from category $B$. If $j\notin S_B^+$, then the customer transitions from $j$ to $\ell\in \NN_B^+$ with probability $\bfrho^B_{j,\ell}$.
Overall, the substitution process is equivalent to an MC model that retains the transition matrix $\boldsymbol{\rho}^B$ but sets the initial arrival to product $j$. Thus, conditioned on $i$ being chosen from $A$, the choice process in category $B$ is an MC model with initial arrival distribution $\bfl^i$ and transition matrix $\bfrho^B$.
\begin{lemma}\label{lem:choice-mc}
    When the within-category model $\phi_B$ is a Markov chain choice model with parameters $(\bfpsi^B,\bfrho^B)$, conditioning on $i$ being chosen from category $A$, the induced choice process in category $B$ is equivalent to a Markov chain choice model with parameters $(\bfl^i, \bfrho^B)$.
\end{lemma}

Our model for category $B$ differs from a standard MC model in a critical way: instead of relying on a single fixed initial distribution $\bfpsi^B$, it uses $\bfl^i$ that depends on the first-stage choice $i$ from category $A$. This dependency allows the initial choice to directly shape preferences in category $B$, thereby capturing complementarity. The substitution dynamics within $B$, governed by $\bfrho^B$, remain the same as in the standard MC model. If all $\bfl^i$ were identical across $i\in \NN_A^+$, the complementarity would disappear, and the model in category $B$ would reduce to a standard MC choice model.

\subsection{Beyond Two Categories}\label{sec:dag}
\
A key feature of our model is its extensibility to settings with more than two categories. We represent cross-category dependencies using a directed acyclic graph (DAG) $G=(\mathcal{V},\mathcal{E})$, where each node $U\in\mathcal{V}$ represents a product category and each directed edge $(U,W)\in\mathcal{E}$ indicates that the choice in category $U$ may influence the subsequent choice in category $W$.  We call a category with no incoming edges a \emph{root category}. To introduce the main ideas, we first focus on DAGs in which every non-root category has a single parent. The extension to categories with multiple parents is described in Appendices~\ref{app:multi-parent} and~\ref{app:alg-dag}; all of our optimization results continue to hold under this more general construction.

In the DAG, each edge $(U,W)\in \mathcal{E}$ is associated with its own cross-category transition matrix $\boldsymbol{\lambda}^{U,W}=(\lambda^{U,W}_{i_U,j_W})_{i_U\in \NN_U^+,\, j_W\in \NN_W^+}$. After the customer chooses an outcome $i_U\in S_U^+$ in category $U$, including the no-purchase outcome $0_U$, the customer forms an initial attraction to an option $j_W\in\NN_W^+$ in category $W$ with probability $\lambda^{U,W}_{i_U,j_W}$. This initial attraction plays exactly the same role as in our two-category model. If $j_W\in S_W^+$, the customer selects it directly; otherwise, the customer substitutes among the available options in $S_W^+$ according to the within-category choice model $\phi_W$, conditional on $j_W$ being the initial preferred option. Thus, once $j_W$ is determined, the subsequent choice in $W$ is conditionally independent of the upstream outcome $i_U$.

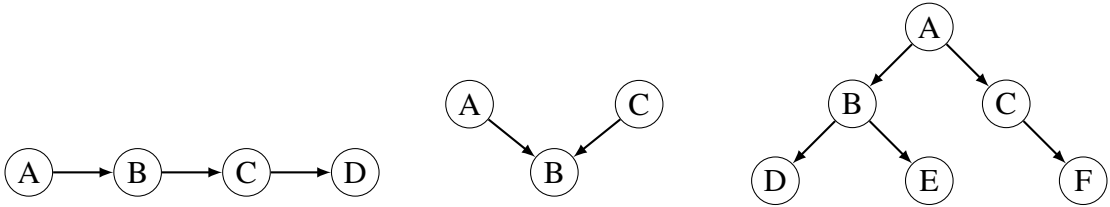
\begin{figure}[htbp]
\vspace{-0.8em}
\centering
\tikzstyle{vertex}=[circle, draw, minimum size=1.5em, inner sep=0pt]
\tikzstyle{arrow}=[-latex, thick]

\begin{tikzpicture}[node distance=0.8cm]
    \node[vertex] (A1) {A};
    \node[vertex] (B1) [right=of A1] {B};
    \node[vertex] (C1) [right=of B1] {C};
    \node[vertex] (D1) [right=of C1] {D};

    \draw[arrow] (A1) -- (B1);
    \draw[arrow] (B1) -- (C1);
    \draw[arrow] (C1) -- (D1);

    \node[vertex] (B2) [right=2cm of D1] {B};
    \node[vertex] (A2) [above left=0.45cm and 0.67cm of B2] {A};
    \node[vertex] (C2) [above right=0.45cm and 0.67cm of B2] {C};

    \draw[arrow] (A2) -- (B2);
    \draw[arrow] (C2) -- (B2);

    \node[vertex] (A3) [above right=0.57cm and 3.4cm of C2] {A};
    \node[vertex] (B3) [below left=of A3] {B};
    \node[vertex] (C3) [below right=of A3] {C};
    \node[vertex] (D3) [below left=of B3] {D};
    \node[vertex] (E3) [below right=of B3] {E};
    \node[vertex] (F3) [below right=of C3] {F};

    \draw[arrow] (A3) -- (B3);
    \draw[arrow] (A3) -- (C3);
    \draw[arrow] (B3) -- (D3);
    \draw[arrow] (B3) -- (E3);
    \draw[arrow] (C3) -- (F3);
\end{tikzpicture}
    \caption{Examples of cross-category dependency structures.}
    \label{fig:graph}
    \vspace{-10pt}
\end{figure}

This representation accommodates several common purchasing patterns, as illustrated in Figure~\ref{fig:graph}. A chain structure (left) represents a sequential process, such as building a PC, where the choice in each category influences the next one. A convergent structure (middle) represents several upstream categories that may influence a common downstream category, such as purchases of a television or gaming console affecting the choice of a soundbar. A tree structure (right) captures a primary choice that branches into several downstream decisions, such as a flight purchase influencing subsequent hotel and rental-car choices.

When a category $W$ has multiple parents, each realized parent outcome induces an initial-state distribution in $W$. We aggregate these parent-specific distributions into a single pooled distribution; Appendix~\ref{app:multi-parent} gives the formal construction. The customer then enters $W$ according to this pooled distribution and follows the same within-category choice model $\phi_W$, thereby generating a single choice in $W$.

\section{Assortment Optimization}\label{sec:3-opt}
In this section, we formulate the multi-category assortment optimization problem, where the goal is to jointly select assortments across categories to maximize total expected revenue. We begin with the two-category case, showing that the unconstrained problem can be solved efficiently when within-category choices follow MC models. We then extend the result to a general DAG. Finally, we establish the computational hardness of cardinality-constrained assortment optimization.
\subsection{Assortment Optimization for Two Categories}\label{sec:optimization}
Consider first the two-category case. Let $r_i^A$ and $r_j^B$ denote the unit prices of products $i\in \NN_A$ and $j\in \NN_B$, with $r_{0_A}^A=r_{0_B}^B=0$. The problem is to select $S_A\subseteq \NN_A$ and $S_B\subseteq \NN_B$ to maximize expected revenue:
\begin{equation}\label{eq:opt}
\max_{S_A\subseteq \NN_A, S_B\subseteq \NN_B}\sum_{i\in S_A^+}\sum_{j\in S_B^+}(r_i^A+r_j^B)P(i,j,S_A,S_B),
\end{equation}
where $P(i,j,S_A,S_B)=\phi_A(i,S_A) (\lambda_{i,j}+ \sum_{\ell\in \NN_B\setminus S_B }\lambda_{i,\ell}\phi_B(j,S_B|\ell))$ is given by Equation \eqref{eq:2-prob}.

Recall that for general ranking-based choice models, even the single-category unconstrained assortment problem is NP-hard to approximate within any non-trivial factor \citep{aouad2018approximability}. We therefore focus on the setting where both $\phi_A$ and $\phi_B$ are MC models, for which the single-purchase assortment optimization problem is polynomial-time solvable \citep{blanchet2016markov}. Extending to two categories is generally harder: the optimal assortment in one category can depend on the other. For previous MNL-based models that incorporate cross-category complementarity, the unconstrained optimization problem is NP-hard \citep{ke2022cross,chen2022assortmentoptimizationmultivariatemnl,tulabandula2023}. In contrast, our framework admits a polynomial-time algorithm even when both categories follow MC models.

Recall that $\phi_A$ and $\phi_B$ are defined by $(\bfpsi^A,\bfrho^A)$ and $(\bfpsi^B,\bfrho^B)$, respectively, and that $\bfl^i$ denotes the transition probabilities from product $i\in \NN_A^+$ to category $B$. Lemma~\ref{lem:choice-mc} shows that conditional on $i$, the choice process in $B$ is equivalent to an MC model $(\bfl^i,\bfrho^B)$. For notational simplicity, we denote this conditional model as $\phi_{B|i}$.

The key to solving the two-category assortment problem is to show that the joint assortment problem can be decoupled.
First, observe that if the optimal assortment in $B$ were known, then computing the optimal assortment for $A$ would reduce to a standard single-category assortment problem under MC with adjusted prices. Specifically, the price of any product $i$ in $A$ is augmented by the expected revenue from category $B$ given that $i$ is chosen from $A$. This idea of using adjusted prices has been commonly used in multi-purchase models \citep{ke2022cross,chen2024assortment}. However, computing the optimal assortment for $B$ seems to require knowing the optimal assortment for $A$, since the choice in $A$ influences the initial attraction in $B$. Remarkably, our model breaks this dependency when $\phi_B$ is an MC model. The solution hinges on a crucial property of the MC choice model: for the unconstrained assortment problem, there exists an optimal assortment that is invariant to the initial arrival distribution. This `invariant' optimal assortment depends on only the transition matrix ($\bfrho$) and product prices ($\bfr$).  Such an invariant optimal assortment can be computed using the algorithms of \cite{blanchet2016markov} or \cite{desir2020constrained}; see Lemma~\ref{lem:mc-invariance} in Appendix~\ref{app:mc-lemma}. In our framework, the choice $i$ from category $A$ affects only the arrival distribution in $B$ (via $\bfl^i$), while leaving the within-category transition matrix $\bfrho^B$ unchanged. Consequently, the invariant optimal assortment for the MC model in category $B$ remains optimal for the joint problem, regardless of the assortments in $A$. Thus, we fully decouple the problem: the optimal assortment for $B$ can be computed first, and the solution for $A$ then follows by solving a single-category MC problem with adjusted prices. We now formalize this procedure in Algorithm \ref{alg:mc_assortment}. To simplify the algorithm description, we assign the outside option $0_A$ a price $r_{0_A}^A=0$.

The algorithm first computes the invariant optimal assortment for $B$, denoted by $S_B^{\ALG}$. Since this assortment is independent of the initial arrival probabilities, we can use arbitrary initial arrival probabilities; for concreteness, the algorithm uses $\bfl^1$ in $\phi_{B|1}$. The invariant optimal assortment can be obtained in polynomial time using the iterative algorithms of \citet{blanchet2016markov} or \citet{desir2020constrained} (Algorithm 1 in each paper). The output of these algorithms is invariant to the arrival probabilities, even though not all optimal unconstrained assortments of an MC model are invariant (see Example \ref{eg:no-mc} in Appendix \ref{app:no-mc}).

\begin{algorithm}[H]
\caption{Unconstrained Two-category Assortment Optimization}
\label{alg:mc_assortment}
\SetAlgoNoLine
\textbf{Input:} Prices $(r_i^A)_{i\in \NN_A^+}$ and $(r_j^B)_{j\in \NN_B}$, choice models $\phi_A$ and $\phi_B$, transition probability vectors $(\bfl^i)_{i\in \NN_A^+}$;\\
Compute the invariant optimal unconstrained assortment $S_B^{\ALG}$ for category $B$ under $\phi_{B|1}$\;
\For{\(i\in\NN_A^+\)}{
    Compute the expected revenue from $B$:  $R_i(S_B^{\ALG})= \sum_{j\in S_B^{\ALG}}r_j^B\phi_{B|i}(j,S_B^{\ALG})$\;
    Define adjusted price for product $i$ in $A$: $r_i'=r_i^A+R_i(S_B^{\ALG})$\;
}
Compute the optimal unconstrained assortment $S_A^{\ALG}$ under $\phi_A$ using $(r_i')_{i\in\NN_A^+}$ as prices\;
\textbf{Return:} $S_A^{\ALG}$ and $S_B^{\ALG}$.
\end{algorithm}

Given $S_B^{\ALG}$, the algorithm then computes an adjusted price $r_i'$ for each $i \in \NN_A^+$ by adding the base price $r_i^A$ to the expected revenue from category~$B$, $R_i(S_B^{\ALG})$. Finally, it solves the unconstrained assortment optimization for category~$A$ under $\phi_A$ with adjusted prices, yielding $S_A^{\ALG}$. A subtlety arises in this step: the adjusted price of the outside option $0_A$, denoted $r_{0_A}'$, may be strictly positive. This differs from the standard assortment optimization setting, where the outside option has a price of zero. We show that this can be easily addressed by applying a simple price-shifting transformation that subtracts $r_{0_A}'$ from all adjusted prices to set the price of the outside option to zero. In Lemma~\ref{prop:3-opt-equ} (Appendix~\ref{app:assortment-shift}), we show that this transformation does not change the optimal solution. Thus, for computing the optimal assortment in $A$, we can use any polynomial-time MC assortment optimization algorithm. We now state the main result.
\begin{theorem}\label{thm:opt}
For any Markov chain choice models $\phi_A$ and $\phi_B$, Algorithm \ref{alg:mc_assortment} computes an optimal solution to the problem \eqref{eq:opt} in polynomial time.
\end{theorem}
\proof{Proof.}
Let $R(S_A,S_B)$ denote the total expected revenue for assortments $(S_A,S_B)$ as given in~\eqref{eq:opt}. Let $(S_A^*,S_B^*)$ be an optimal solution to~\eqref{eq:opt}, and $(S_A^{\ALG},S_B^{\ALG})$ denote the output of Algorithm~\ref{alg:mc_assortment}. We first show the optimality of $(S_A^{\ALG},S_B^{\ALG})$ and then show that the algorithm is polynomial-time.

First, consider the subproblem for category $B$, assuming a fixed optimal assortment $S_A=S_A^*$ for category $A$. In this case, from \eqref{eq:opt}, $S_B^*$ is the optimal solution of the following problem:
\[
\vspace{-5pt}
\max_{S_B\subseteq \NN_B}\; \sum_{j\in S_B^+}\sum_{i\in S_A^{*+}} \phi_A(i,S_A^*)\phi_{B|i}(j,S_B) r_j^B.
\]
Given $S_A^{*+}$, we define an aggregated choice model for category $B$ as
\(
\phi_{B|S_A^*}(j, S_B) = \sum_{i \in S_A^{*+}} \phi_A(i, S_A^{*})\phi_{B|i}(j,S_B)\).
Since each $\phi_{B|i}$ is an MC model with the same transition matrix $\bfrho^B$ but different initial arrival probabilities, their convex combination is also an MC model with transition matrix $\bfrho^B$ and a new initial arrival distribution given by $\boldsymbol{\psi}^{S_A^*} = \sum_{i \in S_A^{*+}} \phi_A(i ,S_A^*) \boldsymbol{\lambda}^i$. Thus, the subproblem is a standard unconstrained assortment optimization for MC with parameters $(\boldsymbol{\psi}^{S_A^*}, \bfrho^B)$. Algorithm \ref{alg:mc_assortment} computes the invariant assortment $S_B^{\ALG}$ that is optimal for any initial arrival probability; thus, it must also be optimal for the specific initial distribution $\boldsymbol{\psi}^{S_A^*}$ induced by $S_A^*$. This implies its expected revenue is at least as high as that of $S_B^*$: $\sum_{j\in S_B^{\ALG}} \sum_{i\in S_A^{*+}}\phi_A(i,S_A^{*})\phi_{B|i}(j,S_B^{\ALG}) r_j^B\ge \sum_{j\in S_B^*} \sum_{i\in S_A^{*+}}\phi_A(i,S_A^{*})\phi_{B|i}(j,S_B^*) r_j^B$. Therefore, considering the total expected revenue, we have our first key inequality:
\begin{equation}
\label{eq:prf-1}
    R(S_A^*,S_B^{\text{ALG}})\ge R(S_A^*,S_B^*).
\end{equation}

Next, we fix the assortment for $B$ to be $S_B^{\ALG}$. The problem of computing the assortment for $A$ then reduces to:
\(
\max_{S_A\subseteq \NN_A} \sum_{i\in S_A^+}(r_i^A+R_i(S_B^{\text{ALG}}))\phi_A(i,S_A)\),
where $R_i(S_B^{\ALG})=\sum_{j\in S_B^{\ALG}}r_j^B\phi_{B|i}(j,S_B^{\ALG})$ is the expected revenue from $B$ if product $i$ is chosen in $ A$. This is an unconstrained assortment optimization problem under $\phi_A$ with adjusted prices $r_i^A+R_i(S_B^{\ALG})$ for $i\in \NN_A^+$. Since $S_A^{\text{ALG}}$ is the optimal solution for this problem, we have $R(S_A^{\text{ALG}},S_B^{\text{ALG}})\ge R(S_A^*,S_B^{\text{ALG}})$. Combining with \eqref{eq:prf-1} yields $R(S_A^{\text{ALG}},S_B^{\text{ALG}})\ge R(S_A^*,S_B^*)$. Since $(S_A^*,S_B^*)$ is a joint optimal assortment, this proves that $(S_A^{\text{ALG}},S_B^{\text{ALG}})$ is also optimal for the joint problem \eqref{eq:opt}.

Finally, we analyze the runtime. The algorithm runs polynomial-time unconstrained MC assortment optimization for $A$ and $B$. Besides this,
the algorithm performs $O(n_A)$ iterations of polynomial-time revenue calculations. Hence, the overall runtime is polynomial.
 \hfill
\Halmos
\endproof
Because the MNL model is a special case of the MC model, Theorem~\ref{thm:opt} also applies to the \textsc{MarkovMNL} special case used in our empirical comparisons.
\begin{remark}[Generalization beyond MC models]
The optimality result in Theorem \ref{thm:opt} does not require category~$A$ to follow an MC model. It is sufficient that the unconstrained assortment problem under $\phi_A$ can be solved, or $\alpha$-approximated in polynomial time (for example, with the nested logit model), and then the same guarantee carries over to the two-category problem. However, if the choice model for category $B$ is a non-Markovian model such as the nested logit model \citep{selena,ma2023assortment}, our invariance-based decoupling argument does not apply, and the tractability of the joint optimization problem remains an open question. The framework nevertheless extends
when category~$B$ admits a fixed assortment that is approximately
optimal uniformly over all initial arrival distributions. In this
case, the decoupling argument yields a corresponding approximation
guarantee; Appendix~\ref{sec:approx-invariance} formalizes this
extension.
\end{remark}
\begin{remark}[Extension to the personalization setting]
Our framework also extends to the \emph{personalization} setting, where the assortment for category $B$ is determined after observing the customer's choice from category $A$. In this setting, the assortment problem in the two categories decouples naturally. Specifically, for each potential choice $i \in \NN_A^+$, we solve for the optimal conditional assortment in category~$B$, denoted $S_{B,i}^{\ALG}$. Then we solve category~$A$ with adjusted prices $r_i^A+\sum_{j\in S_{B,i}^{\ALG}}r_j^B\phi_{B|i}(j,S_{B,i}^{\ALG})$. Notably, the optimality of the procedure only requires that the single-category assortment optimization problems under $\phi_A$ and $\phi_B$ are tractable.
\end{remark}

\subsection{Assortment Optimization for a General DAG}
The Markovian structure of the dependencies between categories enables our algorithm to extend to any general DAG. Because the choice process in each category $W$ depends only on the selection made in its immediate predecessors, we do not need to track the customer's full purchase history. This memoryless property allows us to solve the problem via backward induction, starting from the terminal categories.

We formalize this backward induction using a dynamic programming approach over the dependency graph $G=(\mathcal{V},\mathcal{E})$. The categories are processed in reverse topological order (i.e., each category is visited only after all its successors have been processed), ensuring that when we solve the subproblem for any category, the invariant optimal assortments for all its descendant categories have already been computed. For simplicity, we again first present the recursion for the case in which every non-root category has a single parent. The maximum expected revenue from choosing option $i_U\in \NN_U^+$ and continuing optimally through all of its descendant categories serves as the product's adjusted price. This value is calculated recursively using the following Bellman equation:
\begin{equation} \label{eq:bellman}
r'_{i_U}=r_{i_U}^U+\sum_{W\in \text{Children}(U)}\Bigl(\max_{S_W\subseteq \NN_W} \sum_{i_W\in S_W^+}\phi_{W|{i_U}}(i_W,S_W)r'_{i_W} \Bigr),
\end{equation}
where $\text{Children}(U)$ is the set of children of $U$, and any terminal category $K$ (a category with no children) has $r'_{i_K}=r_{i_K}^K$ for all $i_K\in \NN_K^+$. {
For a category $W$ with multiple parents, the Bellman recursion is modified by weighting the continuation value associated with each incoming edge $(U,W)$ by the corresponding pooling weight $\omega_{U,W}$. The formal recursion for the general case is provided in Appendix~\ref{app:multi-parent}. }

Crucially, this dynamic program is tractable when each choice model $\phi_W$ is an MC model. The inner maximization in the Bellman equation is a single-category assortment problem for category $W$ with a given set of adjusted prices $\{r_{i_W}'\}$. As established previously, there exists an optimal assortment $S_W^*$ for this subproblem that is invariant to the initial arrival distribution. Since the preceding choice $i_U$ only influences the initial arrivals for $\phi_{W|i_U}$, we can find $S_W^*$ and compute the value of the subproblem without needing to know the choice probabilities within category $U$. This decoupling at each step ensures the Bellman recursion can be solved efficiently. The complete algorithm and the formal result are included in Appendix~\ref{app:alg-dag} (see Algorithm~\ref{alg:dag_assortment}).

\begin{remark}[Generalization beyond MC models] Let $\mathcal{V}_r$ be the set of root categories, i.e., nodes with no incoming edges. The MC assumption can be relaxed on these nodes: if all non-root categories follow MC models and each root-category assortment problem is polynomial-time solvable, the backward induction remains polynomial-time. If root category $U$ admits only an $\alpha_U$-approximation, the overall approximation factor is $\min_{U\in\mathcal{V}_r}\alpha_U$.
\end{remark}
 \subsection{Cardinality Constrained Assortment Optimization}\label{sec:cardinality-constr} We now analyze the assortment optimization problem under cardinality constraints, which impose an upper bound on the number of distinct products that can be offered in an assortment. We consider the case where all non-root categories follow an MC model, while the root categories may follow more general choice models. In this section, we show that the problem remains tractable if constraints are limited to the root categories, but becomes computationally hard if they are applied throughout the graph.
\subsubsection{Constraints on Root Categories}
We first consider the setting where a cardinality constraint applies only to the root categories of the DAG. This setting models many practical scenarios, such as offering a limited set of flights (the primary, root products) with a full range of ancillary services (downstream add-ons). In this case, our backward induction approach remains effective. The algorithm is modified only at each root category, where a cardinality-constrained subproblem is solved using the adjusted prices propagated from downstream. If the choice models in the root categories admit a polynomial-time constant-factor approximation algorithm, we can still obtain a constant-factor approximation for the overall problem.

\begin{theorem}\label{prop:car}
    Consider the multi-category assortment problem with cardinality constraints applied only to the root categories $\mathcal{V}_r$. Suppose that for each root category $U\in \mathcal{V}_r$, the cardinality-constrained single-category assortment problem under its choice model admits a polynomial-time $\alpha_U$-approximation algorithm. Then the overall problem can be solved in polynomial time with an approximation factor of $\min_{U\in \mathcal{V}_r}\alpha_U$.
\end{theorem}
The proof is provided in Appendix~\ref{app:root-card-proof}.
As an immediate corollary, if each root category’s cardinality-constrained assortment problem is polynomial-time solvable (i.e., $\alpha_U=1$), then the overall problem can also be solved optimally in polynomial time. This is the case, for example, if all root categories follow the MNL model, for which the cardinality-constrained assortment problem is known to be solvable in polynomial time \citep{rusmevichientong2010dynamic}.

\subsubsection{Constraints on All Categories}
The problem becomes significantly more difficult when cardinality constraints are imposed on all categories. In this fully constrained setting, the key invariance property of the MC model breaks down: the optimal assortment for a downstream category now depends on the initial attraction probabilities, which are determined by choices in parent categories. This interdependence prevents the use of our efficient backward induction algorithm and leads to computational intractability. {Our hardness result is stated under the Exponential Time Hypothesis (ETH), the standard conjecture that 3-SAT cannot be solved in subexponential time \citep{impagliazzo2001complexity,impagliazzo2001problems}. ETH is commonly used to rule out polynomial-time approximation algorithms for dense-subgraph problems, including the Bipartite Densest $\kappa$-Subgraph problem used in our reduction. We show the following result, with proof in Appendix \ref{app:prop-car-hard}. The reduction uses an idea similar to that of Theorem 8 of \citet{chen2022assortmentoptimizationmultivariatemnl}.}
\begin{theorem}\label{prop:car-hard}
There exists a constant $c > 0$ such that, assuming the Exponential Time Hypothesis (ETH),
no polynomial-time algorithm can achieve an $\Omega\!\left(g^{-1/(\log \log g)^c}\right)$-approximation
for the multi-category assortment problem with cardinality constraints on all categories,
where $g$ denotes the total number of products.
\end{theorem}

Theorem~\ref{prop:car-hard} establishes that the fully cardinality-constrained problem is computationally intractable. Given the practical importance of such constraints, we develop a scalable alternating optimization heuristic, \textsc{AltOpt}, for the two-category \textsc{MarkovMC} model. Fixing the assortment in one category reduces the other category's problem to a single-category cardinality-constrained MC assortment problem: fixing $S_A$ determines the arrival distribution in $B$, whereas fixing $S_B$ determines the adjusted prices in $A$. \textsc{AltOpt} therefore alternates between solving the cardinality-constrained assortment optimization problems for $A$ and $B$. Each block update uses the polynomial-time $(1/2-\epsilon)$-approximation algorithm of \citet{desir2020constrained}. Algorithm~\ref{alg:alt_opt} gives the procedure, where $R(S_A,S_B)$ denotes the joint expected revenue in~\eqref{eq:opt}.

\begin{algorithm}[H]
\caption{Alternating Optimization for Cardinality-Constrained Assortments (\textsc{AltOpt})}
\label{alg:alt_opt}
\SetAlgoNoLine
\textbf{Input:} MC model parameters and prices, cardinality limits $K_A,K_B$, approximation parameter $\epsilon\in(0,1/2)$, iteration limit $T\ge 1$, and revenue tolerance $\eta\ge 0$\;
Initialize $S_A^{(0)}\gets\NN_A$ and $R^{\mathrm{best}}\gets-\infty$\;
\For{$t=1,\ldots,T$}{
    Use the single-category MC approximation routine to obtain $S_B^{(t)}$ for
    $\displaystyle\max_{S_B\subseteq\NN_B:\,|S_B|\le K_B}R(S_A^{(t-1)},S_B)$\;
    Use the same routine to obtain $S_A^{(t)}$ for
    $\displaystyle\max_{S_A\subseteq\NN_A:\,|S_A|\le K_A}R(S_A,S_B^{(t)})$\;
    Evaluate the exact joint revenue $R^{(t)}\gets R(S_A^{(t)},S_B^{(t)})$\;
    \If{$R^{(t)}>R^{\mathrm{best}}+\eta$}{
        Set $(S_A^{\mathrm{best}},S_B^{\mathrm{best}})\gets(S_A^{(t)},S_B^{(t)})$ and $R^{\mathrm{best}}\gets R^{(t)}$\;
    }
    \If{$S_A^{(t)}=S_A^{(t-1)}$}{\textbf{break}\;}
}
\textbf{Return:} $(S_A^{\mathrm{best}},S_B^{\mathrm{best}})$.
\end{algorithm}
The algorithm returns the best feasible pair encountered, replacing the incumbent only when the joint revenue improves by more than $\eta$, and stops when two successive category-$A$ assortments coincide or the iteration limit is reached. The numerical results in Section~\ref{sec:cardinality_experiments} show that this algorithm is empirically near-optimal relative to exact MILP solutions and scales to problem sizes where the MILP does not terminate within the time limit.

\section{Parameter Estimation}\label{sec:estimation}
In this section, we develop an estimation procedure for our model. We begin with the case of two categories where substitution within each category is captured by an MNL model, as this case preserves the core cross-category complementarity and highlights the main estimation challenges.

Let $\Theta = (\bfl, \bfv^A, \bfv^B)$ denote the full set of model parameters. Let $\mathcal{H}=\{(S_A^t, S_B^t, a^t, b^t)\}_{t\in [T]} $ be the purchase history from $T$ customers, where $(S_A^t,S_B^t)$ are the offered assortments and $(a^t,b^t)$ are the corresponding choices from $S_A^{t,+}$ and $S_B^{t,+}$. Based on the choice probabilities from Lemma~\ref{lem:choice_prob}, define the observed-data log-likelihood $\mathcal{L}(\Theta):=\log L(\Theta;\mathcal H)$ as
\begin{align}
\mathcal{L}(\Theta)
&= \sum_{t=1}^{T} \Biggl[
 \log\left(
            \frac{v_{a^t}^{A}}{V^{A}(S_A^{t}) + 1}
        \right) + \log\left(
            \lambda_{a^t, b^t} +
            \frac{v_{b^t}^{B}}{V^{B}(S_B^{t}) + 1}
            \sum_{m \in \NN_B \setminus S_B^{t}}
                \lambda_{a^t ,m}
        \right)
\Biggr].
\label{eq:obs_log_likelihood}
\end{align}
The terms involving $\bfv^A$ are decoupled from $(\bfl,\bfv^B)$ and correspond to a standard MNL log-likelihood; thus, $\bfv^A$ can be estimated separately by maximizing that log-likelihood. However, the log-likelihood is not jointly concave in the remaining parameters $(\bfl, \bfv^B)$, as shown in Appendix~\ref{app:estimation_details}. This non-concavity prevents the use of standard convex optimization methods. We therefore develop an Expectation-Maximization (EM) algorithm,
which has been widely used to estimate choice models such as the Markov chain choice model or the general ranking-based choice models \citep{van2017expectation,csimcsek2018expectation}.

 The EM algorithm introduces latent variables for the customer's unobserved initial interest in category $B$. Let $X_{m}^{t}$ be an indicator that is $1$ if the initial downstream state in transaction $t$ is outcome $m \in \NN_B^+$, and $0$ otherwise, and let $\mathbf X$ collect these latent variables. Given the observed downstream choice $b^t$, the compatible latent states are
\(
\mathcal M_t:=\{b^t\}\cup(\NN_B\setminus S_B^t).
\)
The resulting M-step separates into subproblems that are concave in $\bfl$ and in the log-weights of the two MNL models.

We now derive this complete-data log-likelihood.
After a customer chooses $a^t$ from $S_A^{t,+}$, their initial interest transitions to an outcome $m \in \NN_B^+$ with probability $\lambda_{a^t ,m}$. If we knew $m$, the logic would be straightforward: if $m$ is available ($m\in S_B^{t,+}$), the final choice must be $m$; if $m$ is unavailable ($m\notin S_B^{t,+}$), the customer substitutes and chooses outcome $j$ with probability $\frac{v_j^B}{V^B(S_B^t)+1}$, based on the derivation from Lemma \ref{lem:choice_prob}.
The complete-data log-likelihood, $\log L_C(\Theta; \mathcal{H}, \mathbf{X})$, is therefore:
\begin{align*}
\log L_C(\Theta;\mathcal{H}, \mathbf{X})
&= \sum_{t=1}^T \Bigl\{
\log\left(
            \frac{v_{a^t}^{A}}{V^{A}(S_A^{t}) + 1}
        \right) + \\
        & \sum_{m \in \mathcal M_t}  X_m^t \log\left( \lambda_{a^t ,m} \left(
\1\{m = b^t\} + \1\{m \in \NN_B\setminus S_B^t \} \frac{v_{b^t}^B}{V^B(S_B^t) + 1}
\right) \right)
\Bigr\}.
\end{align*}
The EM algorithm iteratively alternates between an Expectation (E) step and a Maximization (M) step. For the convergence result below, fix constants
$0<\varepsilon_\lambda<1/(n_B+1)$ and $0<\varepsilon_v<1<C$.
With the no-purchase weights normalized to $v_{0_A}^A=v_{0_B}^B=1$,
define the compact estimation domain $\mathcal P := \{ (\bfl,\bfv^A,\bfv^B): \sum_{m\in\NN_B^+}\lambda_{i,m}=1,
\lambda_{i,m}\ge\varepsilon_\lambda,
\forall i\in\NN_A^+,\ m\in\NN_B^+; \bfv^A\in[\varepsilon_v,C]^{n_A},
\bfv^B\in[\varepsilon_v,C]^{n_B}\}$. This restriction is imposed only as a regularity condition for estimation. Starting from $\Theta^{(0)}\in\mathcal P$, the algorithm proceeds as follows:

\noindent \emph{Step 1 (Initialization).} Choose initial estimates $\Theta^{(0)}=( \bfl^{(0)}, \bfv^{A,(0)}, \bfv^{B,(0)}) \in \mathcal{P}$ and set the iteration counter $\ell=0$.

\noindent \emph{Step 2 (E-step).} For each transaction $t$ and compatible latent state $m\in\mathcal M_t$, compute the conditional expectation of the latent variable given the observed data and current parameter estimates: $\hat{X}_{m}^{t,(\ell)} = \mathbb{E}[X_m^t | \mathcal{H}, \Theta^{(\ell)}]$.

\noindent \emph{Step 3 (M-step).} Let
$Q(\Theta\mid\Theta^{(\ell)}):=
\mathbb E_{\mathbf X\mid\mathcal H,\Theta^{(\ell)}}
[\log L_C(\Theta;\mathcal H,\mathbf X)]$ and choose
$\Theta^{(\ell+1)}\in\arg\max_{\Theta\in\mathcal P}
Q(\Theta\mid\Theta^{(\ell)})$. Increment $\ell\leftarrow\ell+1$ and return to Step 2. The process is repeated until the change in the observed log-likelihood or the parameters falls below a predefined tolerance.

We provide the detailed derivations for the E-step and M-step in Appendix \ref{app:em-alg} and show that each step can be computed efficiently. The sequence of log-likelihood values generated by the EM algorithm is non-decreasing and converges to a finite limit. Any limit point of the
parameter sequence is a constrained stationary point of $\mathcal{L}$ on $\mathcal P$, although it need not be a global maximizer.
The result is stated in the following theorem, with proof in Appendix \ref{app:thm-em-conv}.
\begin{theorem}
    \label{thm:em-conv}
Suppose $\Theta^{(0)}\in\mathcal P$. Then
\(
\mathcal{L}(\Theta^{(\ell+1)})\ge \mathcal{L}(\Theta^{(\ell)})\) for any \(\ell\ge0\),
the sequence $\{\mathcal{L}(\Theta^{(\ell)})\}_{\ell\ge0}$ converges to a finite
limit, and any limit point of $\{\Theta^{(\ell)}\}_{\ell\ge0}$ is a
constrained stationary point of $\mathcal{L}$ on $\mathcal P$.
\end{theorem}

\noindent\textbf{Identifiability.}
In the first category, the estimation problem reduces to a standard single-category choice model: MNL for \textsc{MarkovMNL} and MC for \textsc{MarkovMC}. Hence, the usual support conditions for these models apply \citep{blanchet2016markov,gupta2020parameter}. For category $B$, a simple sufficient condition for identifying $\boldsymbol{\lambda}$ is that, conditional on each first-stage choice $a=i$, the full downstream assortment $S_B=\NN_B$ is observed. In this case, no downstream substitution occurs, and
\(
\mathbb P(b=j\mid a=i,S_B=\NN_B)=\lambda_{i,j}.
\)
Once $\boldsymbol{\lambda}$ is fixed, the remaining category-$B$ problem reduces
to a standard single-category identification problem again. Thus, the two-category
construction does not introduce a fundamentally new identification difficulty.
For MC models, all-but-one assortments give a simple sufficient condition under
the standard MC choice-model argument of \citet{blanchet2016markov};
\citet{gupta2020parameter} provide weaker sufficient conditions. Details are
in Appendix~\ref{app:identifiability}.

\noindent \textbf{Generalizations of the estimation framework. }
Our EM-based estimation framework can be generalized beyond the MNL models and two categories. Its modular design handles more general choice models and multi-category structures arranged in a DAG, as we detail in Appendix \ref{app:est-dag}. The key idea is to introduce edge-specific latent variables and leverage the likelihood's modularity to decouple the M-step updates for each choice model and transition probability.

\section{Synthetic Data Experiments}\label{sec:synthetic_experiments}

We now evaluate the performance of our proposed model using synthetic data. Our objective is to compare our model with established benchmarks in a controlled environment where the ground truth is known. Specifically, we focus on assessing the model's ability to capture varying degrees of cross-category complementarity and whether improved modeling accuracy translates into better assortment decisions.
We complement this section with real-world case studies in Section \ref{sec:real_data}.

\subsection{Experimental Design}
We begin by describing the data generation process, the models used for comparison, and the evaluation metrics.

\subsubsection{Ground Truth and Data Generation}
We first describe the ground-truth customer choice model, which captures complex dependencies, and then outline the data generation process.

\noindent \textbf{Choice model construction. }
\begingroup
Our ground-truth model is based on the ranking-based choice model, which is a general choice model and has been widely used as the ground-truth model in the literature \citep{aouad2023exponomial}. We use this model to simulate complex real-world choice behavior. We consider two product categories, $A$ and $B$, with $n_A = 12$ and $n_B=12$ products, respectively.

For category $A$, customer preferences are characterized by $m_A = 10$ customer classes, and each class $k$ has an arrival probability $\alpha_k^A$, generated by sampling values $\beta_k$ independently from a uniform distribution $U[0,1]$ and then normalizing: $\alpha_k^A = \beta_k / \sum_{\kappa=1}^{m_A} \beta_\kappa$. Product rankings for each customer class $k$ in category $A$ are generated as a random permutation of all products, including the no-purchase option. To naturally simulate non-consideration, the permutation is truncated at the no-purchase option, meaning any products ranked below it are deemed unacceptable and discarded from the customer's ranking. For category~$B$, we generate $m_B = 10$ baseline rankings $\{\sigma_k^0\}_{k=1}^{m_B}$ and arrival probabilities $\{\psi_k^B\}_{k=1}^{m_B}$ identically. To introduce cross-category complementarity, we conditionally perturb these baseline rankings based on the chosen product $i \in \NN_A$. For each product $i$, we define a subset of complementary goods $C_i \subseteq \NN_B$ by including each product $j \in \NN_B$ independently with probability $0.3$.
Next, we derive the perturbed rank of product $j$. The baseline utility of product $j$ under rank $k$ is defined as $u_{kj}^0 = n_B - \text{rank}_k(j)$. For complementary pairs ($j \in C_i$), we sample a shared utility shock $\delta_{ij} \sim \text{Exp}(1)$ and define an augmented utility $u_{kij} = u_{kj}^0 + \theta \delta_{ij}$, where $\theta \ge 0$ controls the strength of complementarity. For non-complementary products ($j \notin C_i$), or if the no-purchase option was chosen in category~$A$, no utility boost is applied. Consequently, when $\theta = 0$, choices in category~$B$ are strictly independent of the choice in category~$A$. As $\theta$ increases, the conditional ranking increasingly favors the complementary products. The probability of choosing product $j \in \NN_B^+$ from an offered assortment $S_B$, given the prior choice $i$, is then determined by these augmented conditional rankings $\sigma_k^i$:
$
\sum_{k=1}^{m_B} \psi_k^B \1\!\left\{j = \arg \min_{\ell \in S_B^+} \sigma_k^i(\ell)\right\}$.

\endgroup

\noindent \textbf{Transaction data simulation. } {To ensure our results are robust, we perform $20$ independent replications of the entire model configuration generation process. Each replication consists of the following steps. First, we create a new ground-truth model for category $A$ using the procedure described above. Then, keeping the category $A$ model fixed, we generate $10$ corresponding models for category $B$ by setting the complementarity parameter $\theta$ to $10$ distinct values, ranging from $0.0$ (no complementarity) to $9.0$ (strong complementarity). This procedure results in a total of $200$ unique ground-truth model configurations.}
For each model configuration, we simulate a dataset of $14{,}400$ transactions. In each transaction, an assortment for each category is drawn uniformly at random from the power set of its products ($\mathcal{P}(\NN_A)$ and $\mathcal{P}(\NN_B)$). A customer's choice is then generated from the corresponding ground-truth model. This process yields $200$ distinct datasets.

\subsubsection{Models for Comparison}
\leavevmode%
\begingroup
For the synthetic experiments, we instantiate our framework with either an MC or MNL model for both $\phi_A$ and $\phi_B$, denoted by \textsc{MarkovMC} and \textsc{MarkovMNL}, respectively, as detailed in Section~\ref{sec:model}. Although the ground-truth data are generated from a more general ranking-based model, MNL and MC models are the most commonly used tractable specifications. We compare our model against the following benchmarks:
\endgroup
\begin{itemize}
    \item \textsc{MultiMNL}: The two-category MNL model from \cite{ke2022cross}, which was introduced in Section \ref{sec:model-two}. Recall that in this model, choices in category $A$ induce a mixture of MNLs on category $B$, making it a highly expressive model. Finding the optimal assortment for this model requires solving a mixed-integer linear program (MILP).
    \item \textsc{IndMNL}: A baseline using two independent MNL models, one for each category. This model assumes no complementarity across categories, and its optimal assortment can be computed efficiently with the revenue-ordered optimal structure \citep{talluri2004revenue}. Recall that \textsc{IndMNL} is a special case of both \textsc{MarkovMNL} and \textsc{MultiMNL}, while \textsc{MultiMNL} and \textsc{MarkovMNL} cannot be reduced to one another.
      \item {{\textsc{IndMC}: A baseline using two independent MC models. Similar to the independent MNL model, this model assumes no complementarity across categories, and we compute the optimal assortment using the modified greedy algorithm \citep{desir2020constrained}. Note that \textsc{IndMC} is a special case of \textsc{MarkovMC}}.}
\end{itemize}
The independent model benchmarks allow us to isolate the value of modeling cross-category complementarity. In particular, comparing \textsc{MarkovMC} with \textsc{IndMC} distinguishes the gains from complementarity from those arising simply from replacing the MNL model with the richer MC model.

\noindent \textbf{Parameter estimation. } {For each of the $200$ datasets, we perform a $70\%$--$30\%$ random split to create training and test sets. On the training data, we estimate the parameters for \textsc{MarkovMNL}, \textsc{MarkovMC}, and \textsc{IndMC} using the EM algorithm from Section \ref{sec:estimation} and Appendix \ref{app:em-alg}, which terminates when the change in average log-likelihood or the maximum parameter deviation between iterations falls below a threshold of $10^{-4}$. We run EM from five random initializations for each instance,
    drawing the initial parameters uniformly at random, and retain the estimate with
    the highest in-sample log-likelihood. The MNL-based benchmarks \textsc{MultiMNL}
    and \textsc{IndMNL} are estimated by maximizing their respective
    log-likelihood functions, which are concave in the log-weights.} {
Empirically, EM converged quickly and reliably for both \textsc{MarkovMC} and
\textsc{MarkovMNL}. The in-sample log-likelihood increased monotonically in
every run.
\textsc{MarkovMNL} converged in a median of $25$ iterations (at most $39$) and
\textsc{MarkovMC} in a median of $76$ iterations (at most $89$). The estimates were
also robust to initialization: across the five restarts, the converged in-sample
log-likelihoods differed by less than $1.2\times10^{-4}$ per observation for
\textsc{MarkovMC} and by less than $2\times10^{-8}$ for \textsc{MarkovMNL},
indicating that the runs reached essentially the same attained likelihood
value regardless of starting point.}

\subsubsection{Evaluation Metrics} We evaluate the performance of all estimated models on the test sets. {{Unless otherwise specified, each metric is calculated for a specific value of $\theta$ by averaging the outcomes from the $20$ independent replications.}}

\noindent \textbf{Model fit and predictive accuracy. } We use three metrics to evaluate how well each model explains customer behavior:
\begin{itemize}
    \item {{\textbf{Log-likelihood:} The average total log-likelihood on both the training set and test set.}}
    \item {{\textbf{Top-$K$ hit rate:} For each transaction, we check whether the realized choice in category~$A$ is among the top-$K$ predicted alternatives for~$A$, and whether the realized choice in category~$B$ is among the top-$K$ predicted alternatives under the conditional distribution given the realized choice in~$A$. We then report the average of the two hit rates. This metric evaluates the model's ability to predict both stages of the purchase decision while avoiding a mechanical dominance of the larger joint choice space.}}
    \item {\textbf{Mean rank:} For each transaction, we rank the realized
choice in category~$A$ under its predicted distribution and the realized
choice in category~$B$ under its predicted distribution conditional on
the realized choice in~$A$. We report the average rank across the two
categories and all transactions. Lower values indicate better performance.}
\end{itemize}
These metrics are commonly used in related literature and highlight distinct aspects of model fit and predictive accuracy: log-likelihood measures overall model fit, the top-$K$ hit rate offers a practical measure of predictive power, and mean rank provides a nuanced view of how well the model prioritizes the actual chosen product.

\noindent \textbf{Assortment optimization performance. } To assess the economic value of each model, we evaluate the quality of its assortment decisions. For each estimated model, we solve for its optimal unconstrained assortment and then evaluate the expected revenue of the optimal assortment under the ground-truth choice model. This provides a direct measure of each model's practical impact.

\subsection{Results}
For conciseness, we present the key experimental trends graphically here. The detailed numerical results corresponding to Figure~\ref{fig:model_comparison} are tabulated in Appendix~\ref{app:syn}.

\noindent \textbf{Model fit evaluation. }
\begingroup
We evaluate model fit using in-sample and out-of-sample log-likelihoods. As illustrated in Figures \ref{fig:model_comparison}(a) and \ref{fig:model_comparison}(b), \textsc{MarkovMC} consistently achieves the highest log-likelihood across all $\theta > 0$ in both the training and test data. When there is no complementarity ($\theta=0$), \textsc{MarkovMC} and \textsc{IndMC} perform nearly identically, as expected.

The comparison with independent models highlights the value of capturing complementarity: the gap between \textsc{MarkovMC} and \textsc{IndMC} widens as $\theta$ increases, and \textsc{MarkovMNL} similarly pulls ahead of \textsc{IndMNL}. Moreover, the results show that richer within-category substitution also improves model performance. In terms of both in-sample and out-of-sample log-likelihood, \textsc{IndMC} outperforms all MNL-based models when complementarity is relatively weak ($\theta\leq 3.0$), while \textsc{MarkovMC} substantially outperforms every benchmark for all $\theta>0$. Taken together, these comparisons show that the strong performance of \textsc{MarkovMC} comes from combining Markovian cross-category transitions with expressive MC within-category substitution, rather than from richer within-category substitution alone. This finding is consistent with the modeling discussion in Section~\ref{sec:model}: when the within-category choice model is sufficiently expressive, the first-category choice can meaningfully affect the downstream substitution process, allowing our Markovian framework to capture richer forms of cross-category dependence.

Among the MNL-based models, the two non-nested specifications, \textsc{MarkovMNL} and \textsc{MultiMNL}, perform comparably overall. At $\theta=0$, the data-generating process contains no cross-category complementarity. Although \textsc{IndMNL} is a special case of \textsc{MultiMNL} \citep{ke2022cross}, the more flexible \textsc{MultiMNL} specification can fit finite-sample variation by allowing the downstream MNL parameters to depend on the first-stage choice. As a result, \textsc{MultiMNL} improves in-sample fit relative to \textsc{IndMNL} by $0.25\%$, but performs $0.47\%$ worse out of sample. By contrast, \textsc{MarkovMNL}, which also nests \textsc{IndMNL} while allowing richer downstream substitution, improves out-of-sample fit by $0.92\%$; see Tables~\ref{tab:train_ll_synthetic} and~\ref{tab:test_ll_synthetic}. As $\theta$ increases, neither \textsc{MarkovMNL} nor \textsc{MultiMNL} uniformly dominates the other.
\endgroup

\begingroup
\begin{figure}[htbp]
    \centering
    \includegraphics[width=0.9\linewidth]{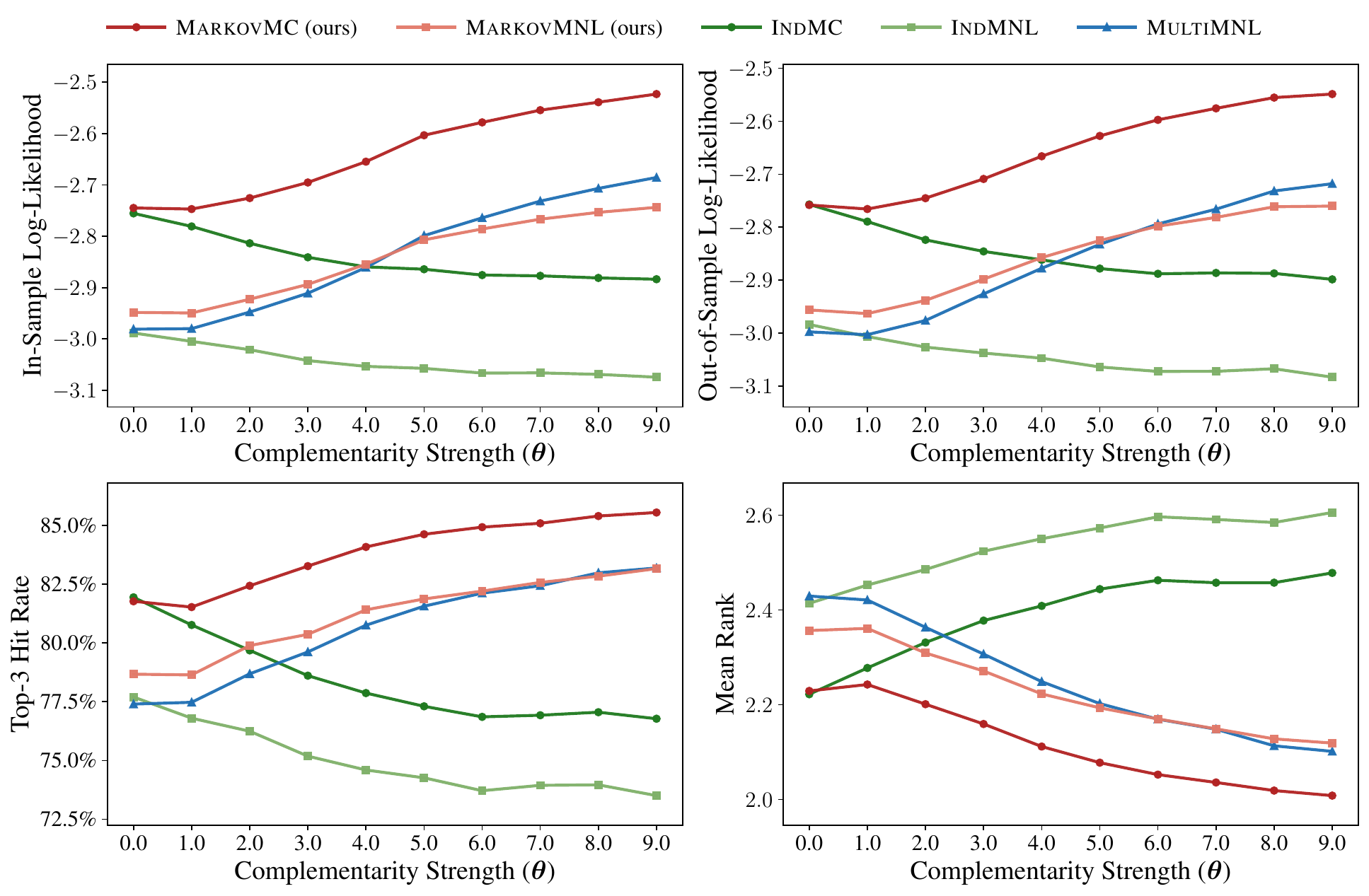}
    \caption{\textbf{Comparison of Model Performance.} The panels display (a) in-sample log-likelihood, (b) out-of-sample log-likelihood, (c) top-3 hit rate, and (d) mean rank across varying levels of complementarity strength ($\theta$).}
    \label{fig:model_comparison}
\end{figure}
\endgroup
\noindent \textbf{Prediction accuracy evaluation. }
\begingroup
We next evaluate prediction accuracy using out-of-sample top-3 hit rates and mean rank. Figures \ref{fig:model_comparison}(c) and \ref{fig:model_comparison}(d) mirror the likelihood results: independent models degrade as complementarity strengthens, the MNL-based complementarity models improve on them, and \textsc{MarkovMC} attains the highest hit rate and lowest mean rank for all $\theta>0$.
\endgroup

\begingroup
\noindent \textbf{Expected revenue from unconstrained optimal assortments. }
To evaluate the economic value of the models, we compare the expected revenue from their unconstrained optimal assortments. We draw prices from log-normal distributions with $\mu_A=2.5$, $\mu_B=2.0$, and $\sigma_A=\sigma_B=2.0$, averaging revenue over $100$ independent price trials. Figure~\ref{fig:revenue_comparison} shows that complementarity-aware models outperform independent benchmarks as $\theta$ grows; relative to \textsc{IndMC}, \textsc{MarkovMC} improves expected revenue by $4$--$7\%$ when $\theta\ge 5.0$.

\endgroup
\begingroup
\begin{figure}
    \centering      \includegraphics[width=0.5\linewidth]{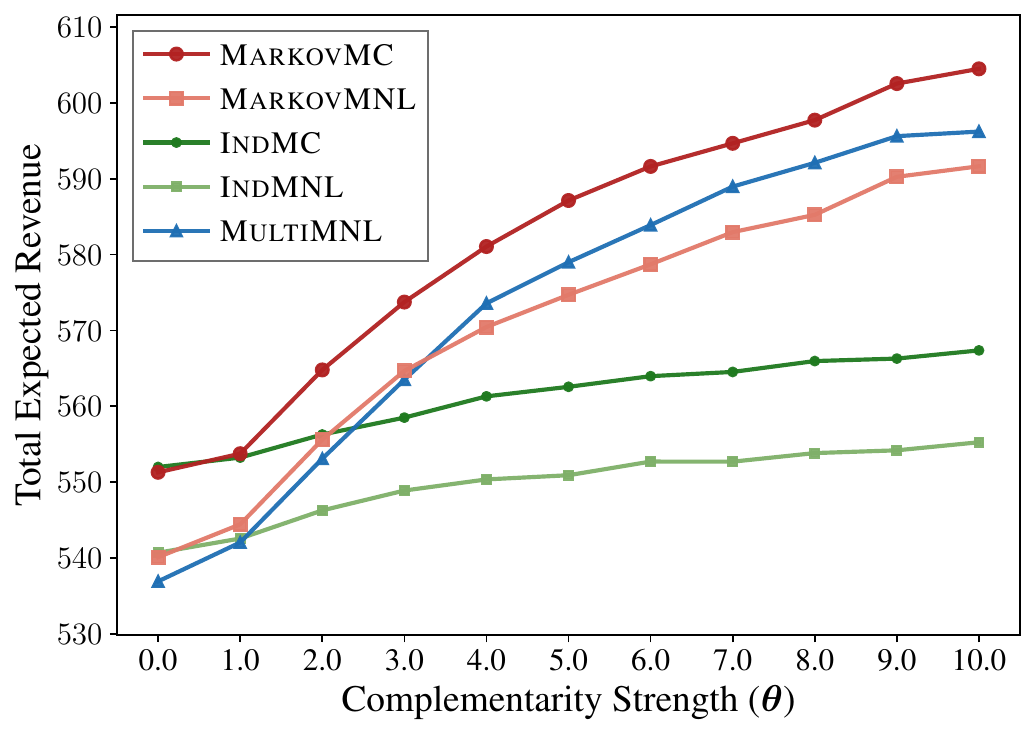}
        \caption{Comparison of Expected Revenue.}
    \label{fig:revenue_comparison}
\end{figure}
\endgroup

\section{Case Studies With Real-World Data}\label{sec:real_data}
In the previous section, we conducted synthetic experiments to evaluate our model under systematically varying complementarity strength. Here, we use real-world transaction data to address the following questions: (i) Does meaningful complementarity exist across product categories in retail settings, and if so, what is its magnitude? (ii) For category pairs exhibiting such complementarity, how does our proposed model perform against established benchmarks? (iii) Can our model uncover actionable managerial insights about product-level interactions?

To answer these questions, we first introduce a metric to quantify complementarity at the category level and use it to identify high-signal category pairs. Next, we evaluate model performance on selected pairs and compare fit and predictive accuracy. In addition, we introduce a second metric to quantify the strength of complementarity between specific products, which we use to extract interpretable, product-level insights. Finally, we demonstrate the framework's flexibility by extending the analysis to a three-category setting.
\subsection{Data and Preprocessing}
{We use a large-scale dataset from a major U.S. grocery store, which contains $5{,}968{,}528$ item-level purchase records covering $5{,}590$ products and $3{,}206$ customers over $97$ weeks; the records group into per-trip shopping baskets (transactions).} This dataset is available to researchers at Stanford and Berkeley by application and has been used in prior literature such as \cite{che2012investigating,ruiz2019shopper}; see \url{https://www.ocf.berkeley.edu/~sberto/sberto-web-site/SGDC/} for more details. We randomly partition the data into a $70\%$ training set and a $30\%$ test set.

The raw transaction data require several preprocessing steps to prepare for choice model estimation.
First, a common challenge in retail data is that the offered assortment for each transaction is unobserved. To address this, we follow the standard approach of inferring the offered set from products purchased within a given time period, as discussed by \cite{jagabathula2019limit} and \cite{mitrofanov2024assortment}. Specifically, for each week, we define the assortment as the set of all unique products purchased by any customer during that period. Second, if there are multiple purchases in a single category, we decompose each transaction into distinct choice observations. For instance, a transaction with purchases $a_1$, $a_2$ from category $A$ and $b_1$ from category $B$ is mapped to two observations: ($a_1, b_1$) and ($a_2, b_1$). {These decomposed category-pair choices, rather than raw purchase records or baskets, are the observations used for estimation.}
Next, we narrow our analysis to transactions with at least one purchase from category $A$. This is necessary because we cannot distinguish customers who considered but rejected all products in a category from those who ignored the category entirely. With this filtering, we can naturally model the conditional choice behavior in category $B$ and observe no-purchase events in that category. While some customers with no purchases in $B$ may not be actively interested in $B$ at all (e.g., they already have enough of those products at home), we assume that all of them consider products from both categories. Finally, to ensure robust model estimation, we restrict the ground set to the set of products that appear in at least $10\%$ of transactions with purchases in their respective category.

\subsection{Metrics for Quantifying Complementarity}\label{sec:cm-metrics}
To quantify the strength and nature of cross-category effects, we develop and use two distinct metrics. The first is the \emph{Aggregate Complementarity Metric (CM)}, which is a data-driven measure used to identify categories with complementarity effects from raw transaction data. The second is the \emph{Specific Complementarity Score (SCS)}, which is a model-based measure derived from estimated parameters to provide interpretable insights into product-level interactions.
\paragraph{Aggregate Complementarity Metric (CM) for Screening.}
This metric is designed to distinguish true, heterogeneous complementarity from simple, homogeneous co-occurrence. A heterogeneous effect implies that the choice of a specific product in one category influences the choice of a specific product in another. In contrast, a homogeneous effect, such as in the classic `beer and diapers' case, reflects a correlation driven by shopper demographics rather than product interactions (i.e., the choice of a specific beer brand does not affect the preference for a specific diaper brand).

To capture this distinction, the CM metric measures how much the purchase distribution in category $B$ deviates when conditioned on the purchase of a specific product $i \in \NN_A$. Let $c(i,j)$ be the count of transactions where product $i \in \NN_A^+$ and product $j \in \NN_B^+$ are purchased together. For this metric, we assume all products in $\NN_B$ are shown to each customer. Thus, the empirical conditional probability of purchasing $j$ given a purchase of $i$ is:
$\hat{P}(j|i) = \frac{c(i,j)}{\sum_{k \in \NN_B^+} c(i,k)}$. We compare this to the aggregate unconditional purchase probability of $j$, calculated as:
 $\hat{P}(j) = (\sum_{i' \in \NN_A^+} c(i',j)) / (\sum_{k \in \NN_B^+} \sum_{i' \in \NN_A^+} c(i',k))$.
For each product $i \in \NN_A^+$, we measure the deviation of its conditional purchase distribution from the aggregate purchase probability distribution using the $L_1$ norm:
$d(i) = \sum_{j \in \NN_B^+} |\hat{P}(j|i) - \hat{P}(j)|$. The final metric is the weighted average of these deviations, where $f_i$ is the total number of transactions in which product $i\in \NN_A^+$ was selected:
\[\text{CM} = \sum_{i \in \NN_A^+} \frac{f_i}{\sum_{k \in \NN_A^+} f_k} d(i). \]

The CM score is bounded between $0$ and $2$, as it represents a weighted $L_1$ distance between probability distributions. A score of $\text{CM} = 0$ occurs when $\hat{P}(j|i) = \hat{P}(j)$ for all $i\in \NN_A^+,j\in \NN_B^+$, indicating that the choice in category $B$ is completely decoupled from the specific choice in category $A$. In such a scenario, our model would reduce to two independent choice models. Larger values indicate stronger heterogeneous cross-category dependence
in the observed data. We use CM as an empirical screening measure to filter out the false identifications of co-occurrence as complementarity. For example, in the `beer and diapers' case, the co-purchase is driven by shopper demographics rather than product-level interactions. This homogeneous effect would yield a low CM score.

\begin{figure}[htbp]
  \centering
  \begin{minipage}[t]{0.46\textwidth}
    \centering
    \includegraphics[width=0.88\linewidth]{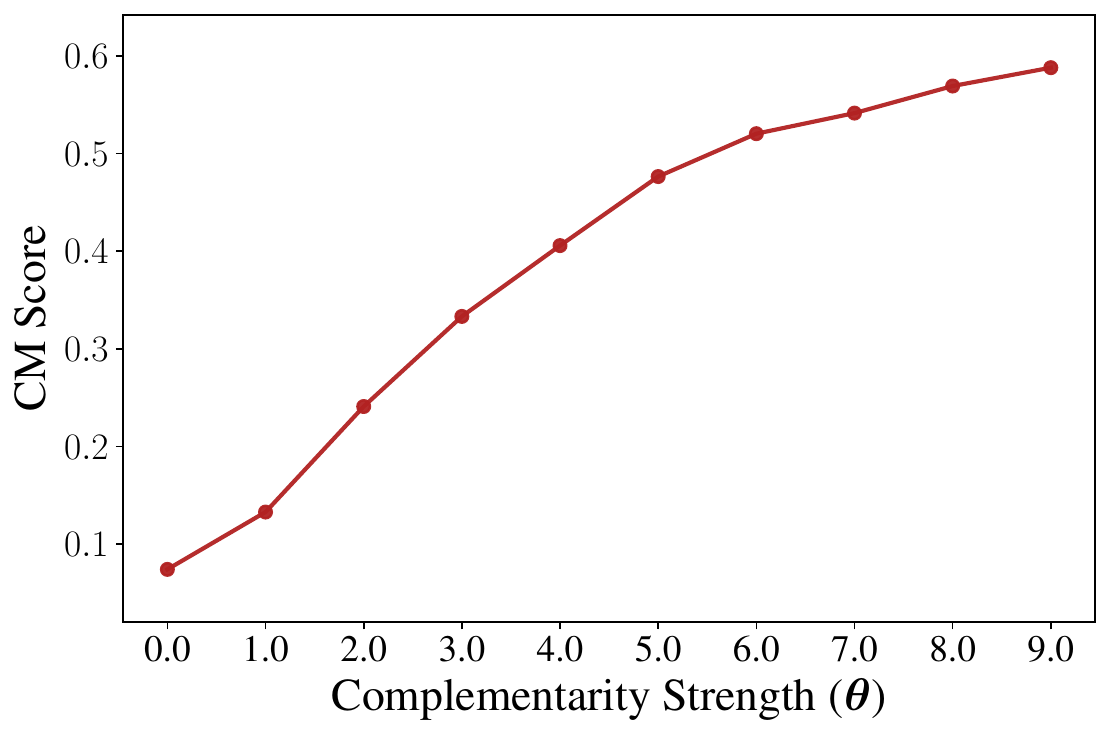}
    \captionof{figure}{CM versus ground-truth complementarity parameter~$\theta$ on synthetic data.}
    \label{fig:cm_score_vs_theta}
  \end{minipage}%
  \hfill
  \begin{minipage}[t]{0.4\textwidth}
  \vspace{-20ex}
    \centering
    \smallXnoCE
    \begin{tabular}{lll}
      \toprule
      \textbf{Category A} & \textbf{Category B} & \textbf{CM} \\
      \midrule
            Diaper         & Beer          & 0.0712 \\
      Coffee         & Creamer       & 0.0738 \\
      Milk/yogurt    & Granola       & 0.1303 \\
      Pasta          & Pasta sauce   & 0.2593 \\
      Patties/franks & Buns          & 0.4202 \\
      Cake mix       & Cake frosting & 0.5844 \\
      \bottomrule
    \end{tabular}
    \vspace{2.2ex}
    \captionof{table}{CM scores for selected category pairs in real data.}
    \label{tab:cm}
  \end{minipage}
   \vspace{-3ex}
\end{figure}

To validate this, we compute CM scores on our synthetic data and real data. Figure~\ref{fig:cm_score_vs_theta} shows that the CM score increases monotonically with our ground-truth complementarity parameter, $\theta$. This result confirms that CM is a valid tool for measuring the strength of cross-category complementarity. Table~\ref{tab:cm} shows the CM scores for several category pairs, illustrating a low score for pairs with homogeneous co-occurrence (e.g., diapers and beer) and progressively higher scores for those with stronger interactions (e.g., cake mix and cake frosting). Subsequently, we use this metric to screen for category pairs with strong complementarity effects in the real-world dataset.

\paragraph{Specific Complementarity Score (SCS) for Interpretation.}
We measure product-level interactions by proposing another metric,  Specific Complementarity Score (SCS), which is computed from the fitted transition probabilities of \textsc{MarkovMC}, and defined as
\(
\tilde{\lambda}_{ij}
=\lambda_{i,j}-\frac{1}{n_A}\sum_{h\in\NN_A}\lambda_{h,j}\) for $i\in\NN_A, j\in\NN_B$.
The baseline is the average initial-attraction probability for product $j$ across all category-$A$ products, assigning equal weight to each product and excluding the no-purchase option. A positive (negative) SCS indicates that initial attraction to $j$ conditional on purchasing $i$ is more (less) likely than this average; a zero score indicates equality with the average. All scores therefore vanish when the purchases across categories are independent, in which case $\lambda_{i,j}$ does not vary with $i$. The SCS captures product-level complementarity in the fitted model through differences in initial attraction across purchases of different category-$A$ products. This interpretation is descriptive, and we do not claim a causal effect of purchasing $i$. In the results below, we first compute the SCS for each product pair and then average these scores within each brand or subcategory pair.

\subsection{Model Evaluation and Interpretability on Real Data}
We now evaluate our model on several combinations of categories with high CM scores. Our goal is to compare model performance and illustrate the interpretability of model estimates.
For evaluation, we use the same metrics as in the synthetic data experiments: log-likelihood, top-3 hit rates, and mean rank.

\paragraph{Example 1: Cake Mixes and Cake Frostings.}
We first consider cake mixes as category $A$ ($39$ products) and cake frosting as category $B$ ($31$ products). The CM score for this pair is $0.5844$, indicating a high degree of observed cross-category dependence. This score is higher than those of the other pairs in Table~\ref{tab:cm} and also aligns with the data for  $\theta\approx 3.0$ in the synthetic experiments (Figure \ref{fig:cm_score_vs_theta}). This pair is also known to exhibit strong complementarity, as noted in previous literature \citep{manchanda1999shopping,sinitsyn2012coordination}.

\begingroup
Table~\ref{tab:model_comparison_cake} shows that accounting for complementarity leads to clear gains in product categories with strong complementarity effects: all complementarity-aware models substantially outperform \textsc{IndMNL} and \textsc{IndMC}. Among these models, \textsc{MarkovMC} achieves the best performance across all metrics. Relative to \textsc{IndMNL}, it improves out-of-sample log-likelihood by 5.52\%, increases the top-3 hit rate by
3.57 p.p., and reduces the mean rank by
6.15\%. These results suggest that \textsc{MarkovMC} captures both complementarity and substitution in the real data very well, which is consistent with the results of synthetic data. \textsc{MultiMNL} and \textsc{MarkovMNL} are comparable in those metrics. Recall that \textsc{MarkovMNL} allows tractable assortment optimization, whereas assortment optimization under \textsc{MultiMNL} is NP-hard.
\endgroup
\begingroup
\begin{table}[htbp]
\smallXnoCE
\centering
\caption{Model fit and prediction accuracy for cake mixes and cake frostings.}
\label{tab:model_comparison_cake}
\resizebox{0.9\textwidth}{!}{%
\begin{tabular}{@{}lllll@{}}
\toprule
Model & In-Sample LL & Out-of-Sample LL & Top-3 HR & Mean Rank \\
\midrule
\textsc{IndMNL}    & -120198.30           & -51739.69           & 44.97\%              & 17.41 \\
\textsc{IndMC}     & -118742.89 (+1.21\%) & -51239.62 (+0.97\%) & 46.45\% (+1.48 p.p.) & 16.90 (-2.93\%) \\
\textsc{MultiMNL}    & -113431.65 (+5.63\%) & -49316.55 (+4.68\%) & 47.39\% (+2.42 p.p.) & 16.81 (-3.45\%) \\
\textsc{MarkovMNL} & -113346.47 (+5.70\%) & -49461.20 (+4.40\%) & 47.44\% (+2.47 p.p.) & 16.81 (-3.45\%) \\
\textsc{MarkovMC}  & -111749.02 (+7.03\%) & -48883.94 (+5.52\%) & 48.54\% (+3.57 p.p.) & 16.34 (-6.15\%) \\
\bottomrule
\end{tabular}%
}
\end{table}
\endgroup

\begin{figure}[htbp]
  \centering
  \subfigure[SCS by cake-mix and frosting brand.]{%
    \includegraphics[width=0.49\linewidth]{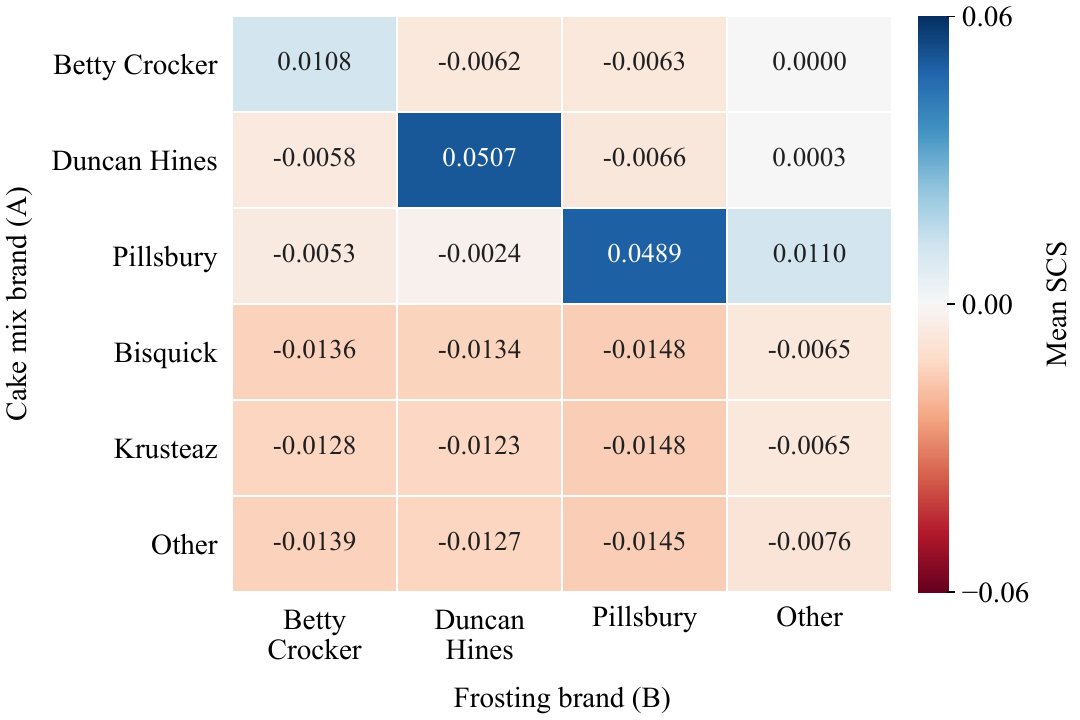}%
\label{fig:brand_complementarity_heatmap}%
  }\hfill
  \subfigure[SCS by meat and bakery subcategory.]{%
    \includegraphics[width=0.49\linewidth]{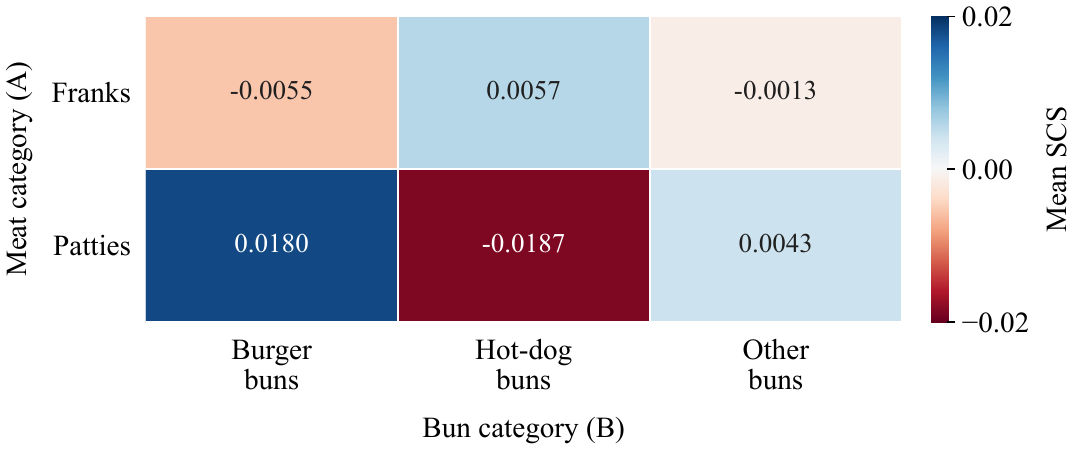}%
    \label{fig:meats_buns_heatmap}%
  }
  \caption{Mean SCS across product pairs within each brand or subcategory pair, computed from the fitted \textsc{MarkovMC} transition probabilities. Scores use a uniform baseline over category-$A$ products; each color scale is centered at zero.}
  \label{fig:combined_heatmaps}
\end{figure}

Figure~\ref{fig:brand_complementarity_heatmap} reports the mean SCS for each brand pair. The scores indicate a preference for matching brands: the within-brand means are positive for Betty Crocker ($0.0108$), Duncan Hines ($0.0507$), and Pillsbury ($0.0489$), whereas all six cross-brand means among these three brands are negative. Duncan Hines has the largest within-brand mean. These estimates are consistent with brand matching and identify potential candidates for coordinated promotions.
\paragraph{Example 2: Meats and Buns.}
Our second case study investigates the classic pairing of meats (category $A$: beef patties, franks, etc.; $30$ products) and buns (category $B$: hamburger, hot dog, etc.; $23$ products, including two English muffins grouped with other buns). This pair has a CM score of $0.4202$, indicating substantial, though slightly less pronounced, complementarity compared to our first example. \begingroup
Table~\ref{tab:model_comparison_meats} shows that modeling complementarity again yields clear gains for this pair. \textsc{MarkovMC} performs best across all metrics: relative to \textsc{IndMNL}, it improves the out-of-sample log-likelihood by $3.18\%$, increases the top-3 hit rate by $2.64$ p.p., and reduces the mean rank by $3.91\%$.
\endgroup

\begingroup
\begin{table}[htbp]

\smallXnoCE
\centering
\caption{Model fit and prediction accuracy for meats and buns. }
\label{tab:model_comparison_meats}
\resizebox{0.9\textwidth}{!}{%
\begin{tabular}{@{}lllll@{}}
\toprule
Model & In-Sample LL & Out-of-Sample LL & Top-3 HR & Mean Rank \\
\midrule
\textsc{IndMNL}    & -132468.40           & -57326.84           & 46.08\%              & 14.83 \\
\textsc{IndMC}     & -131138.79 (+1.00\%) & -56792.61 (+0.93\%) & 46.90\% (+0.82 p.p.) & 14.70 (-0.88\%) \\
\textsc{MultiMNL}    & -128811.21 (+2.76\%) & -55944.46 (+2.41\%) & 47.89\% (+1.81 p.p.) & 14.38 (-3.03\%) \\
\textsc{MarkovMNL} & -128713.70 (+2.83\%) & -56038.91 (+2.25\%) & 47.89\% (+1.81 p.p.) & 14.38 (-3.03\%) \\
\textsc{MarkovMC}  & -127376.52 (+3.84\%) & -55503.71 (+3.18\%) & 48.72\% (+2.64 p.p.) & 14.25 (-3.91\%) \\
\bottomrule
\end{tabular}%
}
\end{table}
\endgroup
 Figure~\ref{fig:meats_buns_heatmap} reports mean SCS scores by subcategory. Franks have a positive mean SCS with hot dog buns ($0.0057$), and patties have a positive mean SCS with burger buns ($0.0180$). In contrast, the mean SCS for patties and hot dog buns is negative ($-0.0187$). Thus, relative to the uniform average across meat products, the estimates associate franks with greater initial attraction to hot dog buns and patties with greater initial attraction to burger buns and lower initial attraction to hot dog buns.

\begingroup
\begin{remark}[Scalability to multiple categories]
Our framework's scalability is demonstrated with a three-category `meats, buns, and condiments' case study. We model the choice process with a tree dependency structure ($B\leftarrow A\rightarrow C$), as this reflects the consumer behavior where the primary choice of meat ($A$) is the main driver for the selection of both the bun ($B$) and the condiment ($C$). The results show that \textsc{MarkovMC} performs best among the models, and \textsc{MarkovMNL} and \textsc{MultiMNL} dominate independent MNL models. The detailed experimental setup and full results are presented in Appendix~\ref{app:three_group_extension}.
\end{remark}
\endgroup
\section{Extensions}
In this section, we present two additional results.
First, we evaluate the computational performance of
\textsc{AltOpt} for cardinality-constrained joint assortment optimization (Algorithm \ref{alg:alt_opt}).
Second, we extend our model to capture bidirectional complementarity
between two categories and study the resulting estimation and
assortment optimization problems.
\subsection{Computational Performance of the Alternating Optimization Heuristic} \label{sec:cardinality_experiments}
We evaluate the alternating optimization heuristic, \textsc{AltOpt}, introduced in Algorithm~\ref{alg:alt_opt} for the cardinality-constrained two-category \textsc{MarkovMC} model. Implementation details and the random instance-generation procedure are provided in Appendix~\ref{sec:heuristic}.

We test \textsc{AltOpt} on $50$ random instances for each problem size $N\in\{10,20,40,80,160,320,640\}$, with $n_A=n_B=N$. When Gurobi certifies a MILP solution as optimal to a relative tolerance of $10^{-4}$ within two hours ($N\le 160$), we use that solution as the numerical benchmark. For larger instances, where the MILP does not terminate within two hours, we compare against the LP relaxation upper bound. Table~\ref{tab:cardinality-altopt} reports the average and worst-case performance ratios.

\begin{table}[htbp]
\centering
\smallXnoCE
\caption{Performance of the alternating optimization heuristic (\textsc{AltOpt}) for cardinality-constrained two-category assortment optimization.}
\label{tab:cardinality-altopt}
\begin{tabular}{@{}rccccc@{}}
\toprule
$N$ & Benchmark & Avg. Ratio & Worst-case Ratio & \textsc{AltOpt} Time (s) & MILP Time (s) \\
\midrule
10  & MILP  & 99.46\% & 92.80\% & 0.0164 & 0.1719 \\
20  & MILP  & 99.83\% & 96.83\% & 0.0229 & 0.6771 \\
40  & MILP  & 99.81\% & 94.54\% & 0.0565 & 6.0129 \\
80  & MILP  & 99.96\% & 99.02\% & 0.4038 & 36.1994 \\
160 & MILP  & 99.99\% & 99.65\% & 1.3730 & 248.3148 \\
320 & LP UB & 88.24\% & 68.98\% & 5.7400 & -- \\
640 & LP UB & 88.60\% & 71.44\% & 35.3900 & -- \\
\bottomrule
\end{tabular}
\end{table}

The results show that \textsc{AltOpt} is near optimal relative to the MILP benchmark whenever Gurobi certifies optimality to the $10^{-4}$ relative tolerance: its average performance ratio is at least $99.46\%$, and its worst-case ratio is at least $92.80\%$. It also scales substantially better than the exact formulation. At $N=320$ and $640$, the MILP does not solve within two hours, while \textsc{AltOpt} terminates in $5.74$ and $35.39$ seconds, respectively. The lower ratios against the LP upper bound for these larger instances should be interpreted conservatively, since the relaxation is only an upper bound rather than the true integer optimum.
\subsection{A Bidirectional Extension for Parallel Categories}
\label{sec:bidirectional}

The sequential model in Section~\ref{sec:model-two} assumes
that one category serves as the primary category and influences
the customer's subsequent choice in another category.
However, in some settings, the two categories may be parallel,
without a clear primary category.
To accommodate such settings, we extend our model by allowing
customers to follow either direction of cross-category dependence.
Specifically, let $D\in\{A\to B,B\to A\}$ denote a latent
direction, with
\(
\mathbb P(D=A\to B)=\gamma,
\mathbb P(D=B\to A)=1-\gamma.
\)
Conditional on $D=A\to B$, the customer follows the sequential
model introduced in Section~\ref{sec:model-two}.
Conditional on $D=B\to A$, the same construction applies
with the roles of the two categories reversed.

Let $\lambda^{AB}_{i,\ell}$ denote the transition probability
from product $i\in\NN_A^+$ to initial state $\ell\in\NN_B^+$,
and let $\lambda^{BA}_{j,k}$ denote the transition probability
from product $j\in\NN_B^+$ to initial state $k\in\NN_A^+$.
For $i\in S_A^+$ and $j\in S_B^+$, define
\[
\begin{aligned}
P^{AB}(i,j,S_A,S_B)
&=
\phi_A(i,S_A)
\left(
\lambda^{AB}_{i,j}
+
\sum_{\ell\in\NN_B\setminus S_B}
\lambda^{AB}_{i,\ell}
\phi_B(j,S_B\mid\ell)
\right),\\
P^{BA}(i,j,S_A,S_B)
&=
\phi_B(j,S_B)
\left(
\lambda^{BA}_{j,i}
+
\sum_{k\in\NN_A\setminus S_A}
\lambda^{BA}_{j,k}
\phi_A(i,S_A\mid k)
\right).
\end{aligned}
\]
The joint choice probability under the bidirectional model is
\(
P^{\mathrm{bi}}(i,j,S_A,S_B)
=
\gamma P^{AB}(i,j,S_A,S_B)
+
(1-\gamma)P^{BA}(i,j,S_A,S_B).
\) Note that the direction $D$ does not necessarily represent the chronological
order of purchases, but rather the category that anchors the
customer's preferences. Thus, the bidirectional model allows
cross-category complementarity to operate in both directions
while preserving the Markovian structure within each direction.
The original sequential model is recovered when $\gamma=1$,
while $\gamma=0$ corresponds to the reversed model.

We next consider assortment optimization.
When the within-category choice models are MC models,
the assortment optimization problem for each direction
can be solved optimally in polynomial time by
Theorem~\ref{thm:opt}.
Although the bidirectional objective generally does not
admit the same decoupling, we can obtain a polynomial-time
$1/2$-approximation by solving the two directional problems
separately and returning the assortment pair with the larger
bidirectional expected revenue.
The formal result and proof are provided in
Appendix~\ref{app:bi-opt}.

For estimation, we extend the EM framework by introducing the
direction $D$ as an additional latent variable. We provide the
detailed algorithm for the MNL specification in
Appendix~\ref{app:bi-opt}. The MC specification can be handled
by an analogous extension of the EM algorithm in
Appendix~\ref{app:markovmc-em}.
\section{Conclusion}\label{sec:conclusion}
In this paper, we introduce a sequential multi-purchase choice model that captures cross-category complementarity through Markovian dependencies on a DAG, while allowing a general within-category choice model to govern substitution. By conditioning each downstream choice on the product selected upstream, the model preserves the tractability retailers need. When within-category choices follow a Markov chain model, we obtain a polynomial-time algorithm for the unconstrained assortment problem that extends from two categories to arbitrary DAGs, together with a scalable EM procedure for estimation. The large-scale grocery experiments reveal interpretable structure such as the brand-matching effect between cake mix and frosting. These results suggest that complementarity matters even for routine purchases, and that managers can use the framework to identify high-complementarity category pairs and product-level interactions that inform assortment, cross-selling, promotions, and shelf placement. Several interesting questions remain open, including stronger guarantees or exact optimization under bidirectional complementarity effects, and the joint optimization of assortments and prices.
\bibliographystyle{informs2014} %
\bibliography{ref.bib} %

@article{ma2023assortment,
  title={When is assortment optimization optimal?},
  author={Ma, Will},
  journal={Management Science},
  volume={69},
  number={4},
  pages={2088--2105},
  year={2023},
  publisher={INFORMS}
}

@article{abdallah2024multi,
  title={Multi-purchase assortment optimization under a general random utility model},
  author={Abdallah, Tarek and Braverman, Anton and Gu, Wenhao},
  journal={Available at SSRN 4842012},
  year={2024}
}

@article{selena,
  title={The non-Markovian nature of nested logit choice},
  author={Li, Selena and Udwani, Rajan},
  journal={Operations Research Letters},
  volume={51},
  number={4},
  pages={461--467},
  year={2023},
  publisher={Elsevier}
}

@article{gupta2020parameter,
  title={Parameter identification in Markov chain choice models},
  author={Gupta, Arushi and Hsu, Daniel},
  journal={Theoretical Computer Science},
  volume={808},
  pages={99--107},
  year={2020},
  publisher={Elsevier}
}

@article{csimcsek2018expectation,
  title={An expectation-maximization algorithm to estimate the parameters of the Markov chain choice model},
  author={{\c{S}}im{\c{s}}ek, A Serdar and Topaloglu, Huseyin},
  journal={Operations Research},
  volume={66},
  number={3},
  pages={748--760},
  year={2018},
  publisher={INFORMS}
}

@article{blanchet2016markov,
  title={A Markov chain approximation to choice modeling},
  author={Blanchet, Jose and Gallego, Guillermo and Goyal, Vineet},
  journal={Operations Research},
  volume={64},
  number={4},
  pages={886--905},
  year={2016},
  publisher={INFORMS}
}

@misc{ruiz2019shopper,
      title={SHOPPER: A Probabilistic Model of Consumer Choice with Substitutes and Complements}, 
      author={Francisco J. R. Ruiz and Susan Athey and David M. Blei},
      year={2019},
      eprint={1711.03560},
      archivePrefix={arXiv},
      primaryClass={stat.ML}
}

@article{che2012investigating,
  title={Investigating effects of out-of-stock on consumer stockkeeping unit choice},
  author={Che, Hai and Chen, Xinlei and Chen, Yuxin},
  journal={Journal of Marketing Research},
  volume={49},
  number={4},
  pages={502--513},
  year={2012},
  publisher={SAGE Publications Sage CA: Los Angeles, CA}
}

@article{mitrofanov2024assortment,
  title={Assortment Optimization with Replacement Options for Retail Platforms with Stockout Risk},
  author={Mitrofanov, Dmitry and Topaloglu, Huseyin and Wang, Yuheng},
  journal={Available at SSRN 4754632},
  year={2024}
}

@article{jagabathula2019limit,
  title={The limit of rationality in choice modeling: Formulation, computation, and implications},
  author={Jagabathula, Srikanth and Rusmevichientong, Paat},
  journal={Management Science},
  volume={65},
  number={5},
  pages={2196--2215},
  year={2019},
  publisher={INFORMS}
}

@article{flores2019assortment,
  title={Assortment optimization under the sequential multinomial logit model},
  author={Flores, Alvaro and Berbeglia, Gerardo and Van Hentenryck, Pascal},
  journal={European Journal of Operational Research},
  volume={273},
  number={3},
  pages={1052--1064},
  year={2019},
  publisher={Elsevier}
}

@article{feng2018substitutability,
  title={On substitutability and complementarity in discrete choice models},
  author={Feng, Guiyun and Li, Xiaobo and Wang, Zizhuo},
  journal={Operations Research Letters},
  volume={46},
  number={1},
  pages={141--146},
  year={2018},
  publisher={Elsevier}
}

@article{gao2021assortment,
  title={Assortment optimization and pricing under the multinomial logit model with impatient customers: Sequential recommendation and selection},
  author={Gao, Pin and Ma, Yuhang and Chen, Ningyuan and Gallego, Guillermo and Li, Anran and Rusmevichientong, Paat and Topaloglu, Huseyin},
  journal={Operations research},
  volume={69},
  number={5},
  pages={1509--1532},
  year={2021},
  publisher={INFORMS}
}

@article{liu2020assortment,
  title={Assortment optimization under the multinomial logit model with sequential offerings},
  author={Liu, Nan and Ma, Yuhang and Topaloglu, Huseyin},
  journal={INFORMS Journal on Computing},
  volume={32},
  number={3},
  pages={835--853},
  year={2020},
  publisher={INFORMS}
}

@article{bai2024multi,
  title={Assortment optimization under the multi-purchase multinomial logit choice model},
  author={Bai, Yicheng and Feldman, Jacob and Segev, Danny and Topaloglu, Huseyin and Wagner, Laura},
  journal={Operations Research},
  volume={72},
  number={6},
  pages={2631--2664},
  year={2024},
  publisher={INFORMS}
}

@article{talluri2004revenue,
  title={Revenue management under a general discrete choice model of consumer behavior},
  author={Talluri, Kalyan and Van Ryzin, Garrett},
  journal={Management Science},
  volume={50},
  number={1},
  pages={15--33},
  year={2004},
  publisher={INFORMS}
}

@article{Bront2009,
author = {Bront, Juan Jos\'{e} Miranda and M\'{e}ndez-D\'{\i}az, Isabel and Vulcano, Gustavo},
title = {A Column Generation Algorithm for Choice-Based Network Revenue Management},
journal = {Operations Research},
volume = {57},
number = {3},
pages = {769-784},
year = {2009}}

@article{udwani2025submodular,
  title={Submodular Order Functions and Assortment Optimization},
  author={Udwani, Rajan},
  journal={Management Science},
  volume={71},
  number={1},
  pages={202--218},
  year={2025},
  publisher={INFORMS}
}

@article{Davis2014,
author = {Davis, James M. and Gallego, Guillermo and Topaloglu, Huseyin},
title = {Assortment Optimization Under Variants of the Nested Logit Model},
journal = {Operations Research},
volume = {62},
number = {2},
pages = {250-273},
year = {2014},
}

@article{strauss2018,
title = {A review of choice-based revenue management: Theory and methods},
journal = {European Journal of Operational Research},
volume = {271},
number = {2},
pages = {375-387},
year = {2018},
issn = {0377-2217},
author = {Arne K. Strauss and Robert Klein and Claudius Steinhardt},
}

@article{gallego2004managing,
  title={Managing flexible products on a network},
  author={Gallego, Guillermo and Iyengar, Garud and Phillips, Robert and Dubey, Abhay},
  journal={Available at SSRN 3567371},
  year={2004}
}

@article{van2017expectation,
  title={An expectation-maximization method to estimate a rank-based choice model of demand},
  author={van Ryzin, Garrett and Vulcano, Gustavo},
  journal={Operations Research},
  volume={65},
  number={2},
  pages={396--407},
  year={2017},
  publisher={Informs}
}

@article{Aurier2014,
author = {Philippe Aurier and Victor Mejía},
title ={Multivariate Logit and Probit models for simultaneous purchases: Presentation, uses, appeal and limitations},

journal = {Recherche et Applications en Marketing (English Edition)},
volume = {29},
number = {2},
pages = {75-94},
year = {2014},
}

@article{Feldman2017,
author = {Feldman, Jacob B. and Topaloglu, Huseyin},
title = {Revenue Management Under the Markov Chain Choice Model},
journal = {Operations Research},
volume = {65},
number = {5},
pages = {1322-1342},
year = {2017}
}

@article{nettleton1999convergence,
  title={Convergence properties of the EM algorithm in constrained parameter spaces},
  author={Nettleton, Dan},
  journal={Canadian Journal of Statistics},
  volume={27},
  number={3},
  pages={639--648},
  year={1999},
  publisher={Wiley Online Library}
}

@article{rusmevichientong2010dynamic,
  title={Dynamic assortment optimization with a multinomial logit choice model and capacity constraint},
  author={Rusmevichientong, Paat and Shen, Zuo-Jun Max and Shmoys, David B},
  journal={Operations research},
  volume={58},
  number={6},
  pages={1666--1680},
  year={2010},
  publisher={INFORMS}
}

@article{sinitsyn2012coordination,
  title={Coordination of price promotions in complementary categories},
  author={Sinitsyn, Maxim},
  journal={Management Science},
  volume={58},
  number={11},
  pages={2076--2094},
  year={2012},
  publisher={INFORMS}
}

@article{mulhern1991implicit,
  title={Implicit price bundling of retail products: A multiproduct approach to maximizing store profitability},
  author={Mulhern, Francis J and Leone, Robert P},
  journal={Journal of Marketing},
  volume={55},
  number={4},
  pages={63--76},
  year={1991},
  publisher={SAGE Publications Sage CA: Los Angeles, CA}
}

@article{rodriguez2011assortment,
  title={Assortment selection and pricing for configurable products under demand uncertainty},
  author={Rodr{\'\i}guez, Betzab{\'e} and Ayd{\i}n, G{\"o}ker},
  journal={European Journal of Operational Research},
  volume={210},
  number={3},
  pages={635--646},
  year={2011},
  publisher={Elsevier}
}

@article{lo2019,
title = {Assortment optimization under the multinomial logit model with product synergies},
journal = {Operations Research Letters},
volume = {47},
number = {6},
pages = {546-552},
year = {2019},
issn = {0167-6377},
author = {Venus Lo and Huseyin Topaloglu}
}

@article{manchanda1999shopping,
  title={The “shopping basket”: A model for multicategory purchase incidence decisions},
  author={Manchanda, Puneet and Ansari, Asim and Gupta, Sunil},
  journal={Marketing science},
  volume={18},
  number={2},
  pages={95--114},
  year={1999},
  publisher={INFORMS}
}

@article{ghoniem2016optimizing,
  title={Optimizing assortment and pricing of multiple retail categories with cross-selling},
  author={Ghoniem, Ahmed and Maddah, Bacel and Ibrahim, Ameera},
  journal={Journal of Global Optimization},
  volume={66},
  pages={291--309},
  year={2016},
  blisher={Springer}
}

@article{block1959random,
  title={Random orderings and stochastic theories of response},
  author={Block, Henry David and Marschak, Jacob},
  year={1959}
}

@article{farias2013nonparametric,
  title={A nonparametric approach to modeling choice with limited data},
  author={Farias, Vivek F and Jagabathula, Srikanth and Shah, Devavrat},
  journal={Management science},
  volume={59},
  number={2},
  pages={305--322},
  year={2013},
  publisher={INFORMS}
}

@article{kahn1962,
author = {Kahn, A. B.},
title = {Topological sorting of large networks},
year = {1962},
issue_date = {Nov. 1962},
publisher = {Association for Computing Machinery},
address = {New York, NY, USA},
volume = {5},
number = {11},
issn = {0001-0782},
journal = {Commun. ACM},
month = nov,
pages = {558–562},
numpages = {5}
}

@article{li2019intertemporal,
  title={Intertemporal price discrimination with complementary products: E-books and e-readers},
  author={Li, Hui},
  journal={Management Science},
  volume={65},
  number={6},
  pages={2665--2694},
  year={2019},
  publisher={INFORMS}
}

@article{walters1991assessing,
  title={Assessing the impact of retail price promotions on product substitution, complementary purchase, and interstore sales displacement},
  author={Walters, Rockney G},
  journal={Journal of marketing},
  volume={55},
  number={2},
  pages={17--28},
  year={1991},
  publisher={SAGE Publications Sage CA: Los Angeles, CA}
}

@article{KWAK201519,
title = {An Analysis of Assortment Choice in Grocery Retailing},
journal = {Journal of Retailing},
volume = {91},
number = {1},
pages = {19-33},
year = {2015},
issn = {0022-4359},
author = {Kyuseop Kwak and Sri Devi Duvvuri and Gary J. Russell}
}

@article{seetharaman2005models,
  title={Models of multi-category choice behavior},
  author={Seetharaman, PB and Chib, Siddhartha and Ainslie, Andrew and Boatwright, Peter and Chan, Tat and Gupta, Sachin and Mehta, Nitin and Rao, Vithala and Strijnev, Andrei},
  journal={Marketing letters},
  volume={16},
  pages={239--254},
  year={2005},
  publisher={Springer}
}

@article{lee2013direct,
  title={A direct utility model for asymmetric complements},
  author={Lee, Sanghak and Kim, Jaehwan and Allenby, Greg M},
  journal={Marketing Science},
  volume={32},
  number={3},
  pages={454--470},
  year={2013},
  publisher={INFORMS}
}

@article{feldman2021assortment,
  title={Assortment optimization under the multi-purchase multinomial logit choice model},
  author={Feldman, Jacob and Segev, Danny and Topaloglu, Huseyin and Wagner, Laura and Bai, Yicheng},
  journal={Available at SSRN 3866734},
  year={2021}
}

@article{heger2024assortment,
  title={Assortment optimization: a systematic literature review},
  author={Heger, Julia and Klein, Robert},
  journal={OR Spectrum},
  volume={46},
  number={4},
  pages={1099--1161},
  year={2024},
  publisher={Springer}
}

@article{chen2024assortment,
  title={Assortment planning for recommendations at checkout under inventory constraints},
  author={Chen, Xi and Ma, Will and Simchi-Levi, David and Xin, Linwei},
  journal={Mathematics of Operations Research},
  volume={49},
  number={1},
  pages={297--325},
  year={2024},
  publisher={Informs}
}

@inproceedings{benson2018,
author = {Benson, Austin R. and Kumar, Ravi and Tomkins, Andrew},
title = {A Discrete Choice Model for Subset Selection},
year = {2018},
isbn = {9781450355810},
publisher = {Association for Computing Machinery},
address = {New York, NY, USA},
booktitle = {Proceedings of the Eleventh ACM International Conference on Web Search and Data Mining},
pages = {37–45},
numpages = {9},
location = {Marina Del Rey, CA, USA},
series = {WSDM '18}
}

@article{gallego2019threshold,
  title={Threshold utility model with applications to retailing and discrete choice models},
  author={Gallego, Guillermo and Tang, Zhuodong and Wang, Ruxian},
  journal={Available at SSRN 3420155},
  year={2019}
}

@article{immorlica2021combinatorial,
  title={Combinatorial assortment optimization},
  author={Immorlica, Nicole and Lucier, Brendan and Mao, Jieming and Syrgkanis, Vasilis and Tzamos, Christos},
  journal={ACM Transactions on Economics and Computation (TEAC)},
  volume={9},
  number={1},
  pages={1--34},
  year={2021},
  publisher={ACM New York, NY, USA}
}

@article{kok2015assortment,
  title={Assortment planning: Review of literature and industry practice},
  author={K{\"o}k, A G{\"u}rhan and Fisher, Marshall L and Vaidyanathan, Ramnath},
  journal={Retail supply chain management: Quantitative models and empirical studies},
  pages={175--236},
  year={2015},
  publisher={Springer}
}

@article{tanusondjaja2016understanding,
  title={Understanding shopper transaction data: how to identify cross-category purchasing patterns using the duplication coefficient},
  author={Tanusondjaja, Arry and Nenycz-Thiel, Magda and Kennedy, Rachel},
  journal={International Journal of Market Research},
  volume={58},
  number={3},
  pages={401--419},
  year={2016},
  publisher={SAGE Publications Sage UK: London, England}
}

@article{Jasin2024,
author = {Jasin, Stefanus and Lyu, Chengyi and Najafi, Sajjad and Zhang, Huanan},
title = {Assortment Optimization with Multi-Item Basket Purchase Under Multivariate MNL Model},
journal = {Manufacturing \& Service Operations Management},
volume = {26},
number = {1},
pages = {215-232},
year = {2024},
}

@article{tulabandula2023,
author = {Tulabandhula, Theja and Sinha, Deeksha and Karra, Saketh Reddy and Patidar, Prasoon},
title = {Multi-purchase Behavior: Modeling, Estimation, and Optimization},
journal = {Manufacturing \& Service Operations Management},
volume = {25},
number = {6},
pages = {2298-2313},
year = {2023}
}

@article{aouad2023exponomial,
  title={The exponomial choice model for assortment optimization: An alternative to the MNL model?},
  author={Aouad, Ali and Feldman, Jacob and Segev, Danny},
  journal={Management Science},
  volume={69},
  number={5},
  pages={2814--2832},
  year={2023},
  publisher={INFORMS}
}

@article{ban2024selling,
  title={Selling personalized substitutes and co-purchases in online grocery retail},
  author={Ban, Gah-Yi and BENBITOUR, Mohammed Hichame and Chen, Boxiao},
  journal={Available at SSRN},
  year={2024}
}

@article{ke2022cross,
  title={Cross-category retailing management: Substitution and complementarity},
  author={Ke, Chenxu and Wang, Ruxian},
  journal={Manufacturing \& service operations management},
  volume={24},
  number={2},
  pages={1128--1145},
  year={2022},
  publisher={INFORMS}
}

@article{desir2020constrained,
  title={Constrained assortment optimization under the Markov chain--based choice model},
  author={D{\'e}sir, Antoine and Goyal, Vineet and Segev, Danny and Ye, Chun},
  journal={Management Science},
  volume={66},
  number={2},
  pages={698--721},
  year={2020},
  publisher={INFORMS}
}

@article{cao2023revenue,
  title={Revenue management under a mixture of independent demand and multinomial logit models},
  author={Cao, Yufeng and Rusmevichientong, Paat and Topaloglu, Huseyin},
  journal={Operations Research},
  volume={71},
  number={2},
  pages={603--625},
  year={2023},
  publisher={INFORMS}
}

@misc{chen2022assortmentoptimizationmultivariatemnl,
      title={Assortment Optimization Under the Multivariate MNL Model}, 
      author={Xin Chen and Jiachun Li and Menglong Li and Tiancheng Zhao and Yuan Zhou},
      year={2022},
      eprint={2209.15220},
      archivePrefix={arXiv},
      primaryClass={cs.DS},
      url={https://arxiv.org/abs/2209.15220}, 
}

@article{aouad2018approximability,
  title={The approximability of assortment optimization under ranking preferences},
  author={Aouad, Ali and Farias, Vivek and Levi, Retsef and Segev, Danny},
  journal={Operations Research},
  volume={66},
  number={6},
  pages={1661--1669},
  year={2018},
  publisher={INFORMS}
}

@article{impagliazzo2001complexity,
  title={On the complexity of {k-SAT}},
  author={Impagliazzo, Russell and Paturi, Ramamohan},
  journal={Journal of Computer and System Sciences},
  volume={62},
  number={2},
  pages={367--375},
  year={2001},
  publisher={Elsevier}
}

@article{impagliazzo2001problems,
  title={Which problems have strongly exponential complexity?},
  author={Impagliazzo, Russell and Paturi, Ramamohan and Zane, Francis},
  journal={Journal of Computer and System Sciences},
  volume={63},
  number={4},
  pages={512--530},
  year={2001},
  publisher={Elsevier}
}

\newpage
\renewcommand{\theHsection}{app.\Alph{section}}
\renewcommand{\theHsubsection}{app.\Alph{section}.\arabic{subsection}}
\renewcommand{\theHsubsubsection}{app.\Alph{section}.\arabic{subsection}.\arabic{subsubsection}}
\begin{APPENDICES}
\section{Omitted Details and Proofs in Section \ref{sec:model}}
\subsection{Invariance of MNL Probabilities Under Conditioning}\label{app:rcm-spec}
\begin{lemma}\label{prop:rcm-mnl}
Let $\NN$ be a ground set and $S\subseteq\NN$ an assortment.
Assume a ranking-based choice model $\phi$ is equivalent to an MNL model with preference weights $\{v_j\}_{j\in\NN}$.
If we condition on the ranking starting with an unavailable product $k\in \NN\setminus S$, then the conditional choice probability $\phi(i,S\mid k)$ is independent of $k$, and is given by
\[
\phi(i,S\mid k)\;=\;\frac{v_i}{\,1+\sum_{j\in S}v_j\,}\;=\;\phi(i,S).
\]
\end{lemma}

\proof{Proof of Lemma \ref{prop:rcm-mnl}.}
The proof proceeds in two parts. First, we show that the standard MNL model is equivalent to a specific ranking-based model. Second, we show that under this ranking model, conditioning on an unavailable product being ranked first does not alter the choice probabilities.

We start by constructing the equivalent ranking-based choice model. Let $v_0=1$ be the weight for the no-purchase option, and set $\NN^+ := \NN\cup\{0\}$. We generate the rankings over $\NN^+$ sequentially. Given a prefix $P\subseteq \NN^+$ of already-ranked products, the next product $x$ is drawn from the set of unranked products $\NN^+\setminus P$ with probability
\[
\Pr[\text{next product is }x \mid \text{prefix is }P]
=\frac{v_x}{\sum_{m\in \NN^+\setminus P} v_m}.
\]
A customer's choice is the first element of $S^+ := S\cup\{0\}$ that appears in the ranking. This occurs after a (possibly empty) prefix $P$ of products from outside the assortment, i.e., $P \subseteq \NN\setminus S$.
We now derive the choice probabilities $\phi(i,S)$ for $i\in S^+$. Consider any step in the process where a choice has not yet been made, meaning the prefix $P$  contains only items from $\NN\setminus S$. The probability of choosing $i\in S^+$ next, relative to choosing the outside option next, is:
\[
\frac{\Pr[\text{next is } i \mid \text{prefix is } P]}{\Pr[\text{next is } 0 \mid \text{prefix is } P]}
=\frac{v_i/\sum_{m\in \NN^+\setminus P} v_m}{1/\sum_{m\in \NN^+\setminus P} v_m}
= v_i.
\]

This ratio is constant and independent of the set of unranked products. Since the final choice probabilities are aggregated over all sequences leading to a choice, this constant odds ratio is preserved. Specifically, for any $i\in S$, we have
\(
\phi(i,S)/\phi(0,S) =v_i/v_0=v_i
\), which implies $\phi(i,S)=v_i\phi(0,S)$.
Moreover, the choice probabilities for all products must sum to $1$: $\sum_{i\in S^+}\phi(i,S)=1$. Thus, \[
\phi(0,S)+\sum_{i\in S}\phi(i,S)=\phi(0,S)+\sum_{i\in S} v_i\,\phi(0,S)=1.
\]
This implies $\phi(0,S)=\frac{1}{1+\sum_{j\in S} v_j}$ and thus $\phi(i,S)=\frac{v_i}{1+\sum_{j\in S} v_j}$. This confirms that the model is equivalent to the MNL model.

Finally, we consider the conditional probability $\phi(i,S|k)$, where the ranking is conditioned to start with an unavailable product $k\in \NN\setminus S$. The key observation is that the sequential generation process is memoryless. After $k$ is ranked first, the subsequent ranking of the remaining items $\NN^+\setminus\{k\}$ follows the same process with the same relative weights. Since $k\notin S^+$, the set of offered choices $S^+$ remains unchanged, and the relative weights $\{v_j\}_{j\in S^+}$ are also unchanged.

Therefore, the derivation of choice probabilities for $S^+$ on the remaining items is identical to the unconditional case. The same odds argument holds, giving
\[
\phi(i,S\mid k)=v_i\,\phi(0,S\mid k)\qquad\text{for all }i\in S.
\]
The normalization remains the same, yielding $\phi(0,S|k)=\frac{1}{1+\sum_{j\in S}v_j}$. It follows that
for any $i\in S$, $\phi(i,S\mid k)=\frac{v_i}{1+\sum_{j\in S} v_j}=\phi(i,S).$ This shows the conditional probability is independent of $k$ and equal to the original MNL probability.
\hfill
\Halmos
\endproof

\subsection{Relationship Between Two-Category Choice Models} \label{app:connect}
In this section, we compare our two-category model with the related models in \cite{ke2022cross} and \cite{chen2022assortmentoptimizationmultivariatemnl}. Recall that we use \textsc{MarkovMC} and \textsc{MarkovMNL} to denote our model with MC and MNL as special cases, respectively. Figure \ref{fig:venn} visualizes the relationships between all related models. We first formally define the choice models in \cite{ke2022cross} and \cite{chen2022assortmentoptimizationmultivariatemnl}.
\begin{figure}[htbp]
    \centering
\includegraphics[width=0.5\linewidth]{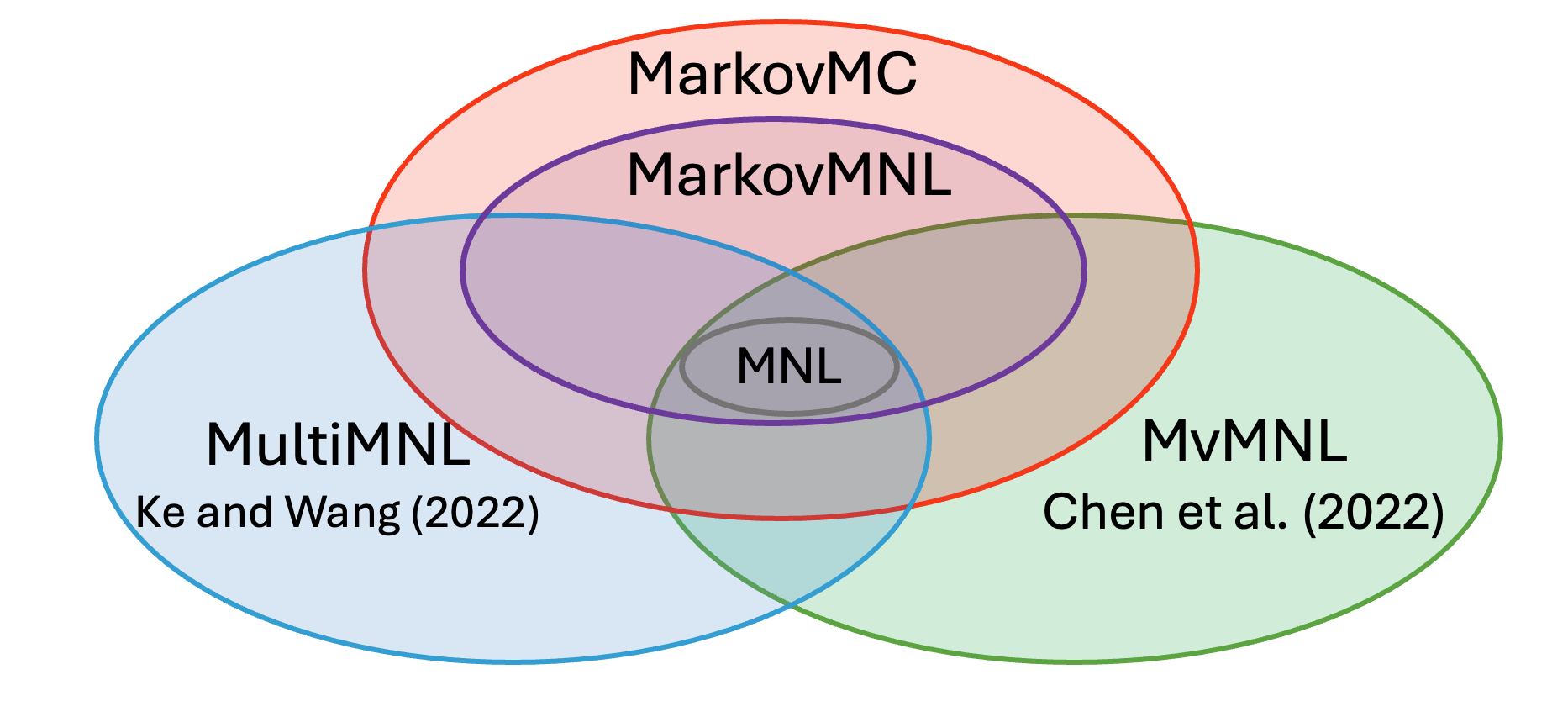}
    \caption{A Venn diagram illustrating the relationships between the choice models.}
    \label{fig:venn}
\end{figure}
We use \textsc{MultiMNL} to denote the sequential MNL choice model from \citet{ke2022cross}. In this model, a customer first chooses a product $i$ from the first-category assortment $S_A$ according to an MNL model with preference weights $\{w_i\}_{i\in \NN_A}$. Conditional on choosing $i$, the customer chooses from $S_B$ according to a distinct MNL model with parameters $\{w_{i,j}\}_{j \in \NN_B}$. A key feature is that each choice $i \in \NN_A^+$ induces its own set of preference parameters for the second category. Without loss of generality, set $w_{0_A}=w_{i,0_B}=1$ for all $i\in \NN_A^+$. The probability of choosing $i\in S_A^+$ and $j\in S_B^+$ is:
\[
P_{\textsc{MultiMNL}}(i,j, S_A, S_B) = \frac{w_i}{\sum_{k\in S_A}w_k+1}\frac{w_{i,j}}{\sum_{\ell \in S_B} w_{i,\ell}+1}.
\]

We use \textsc{MvMNL} to denote the multivariate MNL model from \citet{chen2022assortmentoptimizationmultivariatemnl}. In this model, each bundle $(i,j)$ for $i \in \NN_A^+$ and $j \in \NN_B^+$ has a preference weight $u_{i,j}$. The preference weight of the no-purchase option is set to $u_{0_A,0_B}=1$. The probability of choosing bundle $(i,j)\in S_A^+\times S_B^+$ is:
\[
P_{\textsc{MvMNL}}(i,j , S_A, S_B) = \frac{u_{i,j}}{\sum_{k \in S_A^+} \sum_{\ell \in S_B^+} u_{k,\ell}}.
\]
A core property of this model is that the odds of choosing bundle $(i,j)$ over $(k,\ell)$, given by $u_{i,j}/u_{k,\ell}$, are independent of the offered assortment, provided that all of these products are included in the assortments.

For each model considered here, let $P_{\textsc{model}}(i,j,S_A,S_B)$
denote the probability of choosing outcomes $i\in S_A^+$ and $j\in S_B^+$.
Define the conditional probability
\[
P_{\textsc{model}}(j,S_B\mid i,S_A)
=\frac{P_{\textsc{model}}(i,j,S_A,S_B)}
{\sum_{\ell\in S_B^+}P_{\textsc{model}}(i,\ell,S_A,S_B)}.
\]
For the models compared below, this conditional probability is independent
of $S_A$ once $i$ is given, so we write $P_{\textsc{model}}(j,S_B\mid i)$.

\subsubsection*{Comparison with \cite{ke2022cross}}

\begin{lemma}\label{prop:STC-multimnl}
The \textsc{MarkovMNL} model and the \textsc{MultiMNL} model are distinct; neither model class contains the other.
\end{lemma}

\proof{Proof of Lemma \ref{prop:STC-multimnl}. }We prove this by constructing counterexamples in both directions.

(\emph{\textsc{MarkovMNL} $\not\subseteq$ \textsc{MultiMNL}}).
Consider a universe with $\NN_A=\{1\}$ and $\NN_B=\{2,3\}$. Let the \textsc{MarkovMNL} parameters be $\lambda_{1,2}=\lambda_{1,3}=1/3$ and $v_{0_B}^B=v_2^B=1, v_3^B=2$. (We omit parameters related to the first-stage no-purchase option $0_A$ as they are not needed here.)

Conditional on purchasing product $1$ from category $A$, the choice probabilities for products in category $B$ are:
\begin{itemize}
    \item For assortment $S_B = \{2,3\}$:
    $P_{\textsc{MarkovMNL}}(2, S_B|1) = \lambda_{1,2} = \frac{1}{3}$ and
    $P_{\textsc{MarkovMNL}}(3, S_B|1) = \lambda_{1,3} = \frac{1}{3}$.
    \item For assortment $S_B = \{2\}$:
    $P_{\textsc{MarkovMNL}}(2, S_B|1) = \lambda_{1,2} + \lambda_{1,3} \frac{v^B_2}{v^B_2+1} = \frac{1}{3} + \frac{1}{3}  \frac{1}{2} = \frac{1}{2}$.
    \item For assortment $S_B = \{3\}$:
    $P_{\textsc{MarkovMNL}}(3, S_B|1) = \lambda_{1,3} + \lambda_{1,2} \frac{v^B_3}{v^B_3+1} = \frac{1}{3} + \frac{1}{3}  \frac{2}{3} = \frac{5}{9}$.
\end{itemize}
For a \textsc{MultiMNL} model to match these probabilities, its parameters $\{w_{1,j}\}$ must satisfy:
\begin{itemize}
    \item From $S_B = \{2,3\}$: $\frac{w_{1,2}}{w_{1,2}+w_{1,3}+1} = \frac{1}{3}$ and $\frac{w_{1,3}}{w_{1,2}+w_{1,3}+1} = \frac{1}{3}$, which implies $w_{1,2}=w_{1,3}$.
    \item From $S_B = \{2\}$: $\frac{w_{1,2}}{w_{1,2}+1} = \frac{1}{2}$, which gives $w_{1,2}=1$.
    \item From $S_B = \{3\}$: $\frac{w_{1,3}}{w_{1,3}+1} = \frac{5}{9}$, which gives $w_{1,3}=\frac{5}{4}$.
\end{itemize}
These conditions imply both $w_{1,2}=w_{1,3}$ and $w_{1,2}\ne w_{1,3}$, a contradiction. Thus, no \textsc{MultiMNL} model can represent this \textsc{MarkovMNL} instance.

	(\emph{\textsc{MultiMNL} $\not\subseteq$ \textsc{MarkovMNL}}).
    Consider a universe with $\NN_A=\{1,2\}$ and $\NN_B=\{3,4\}$. Let the \textsc{MultiMNL} parameters be $w_{1,3}=w_{1,4}=2$, and $w_{2,3}=w_{2,4}=w_{2,0_B}=w_{1,0_B}=1$.

To represent this instance with a \textsc{MarkovMNL} model, we first match the choice probabilities when the full assortment $S_B=\{3,4\}$ is offered. This determines the required $\bfl$ parameters by equating them with the \textsc{MultiMNL} probabilities:
\begin{align*}
    \lambda_{1,3} &= P_{\textsc{MultiMNL}}(3, \{3,4\}|1) = \frac{w_{1,3}}{w_{1,3}+w_{1,4}+1} = \frac{2}{5}, \\
    \lambda_{1,4} &= P_{\textsc{MultiMNL}}(4, \{3,4\}|1) = \frac{w_{1,4}}{w_{1,3}+w_{1,4}+1} = \frac{2}{5}, \\
    \lambda_{2,3} &= P_{\textsc{MultiMNL}}(3, \{3,4\}|2) = \frac{w_{2,3}}{w_{2,3}+w_{2,4}+1} = \frac{1}{3}, \\
    \lambda_{2,4} &= P_{\textsc{MultiMNL}}(4, \{3,4\}|2) = \frac{w_{2,4}}{w_{2,3}+w_{2,4}+1} = \frac{1}{3}.
\end{align*}
Now, consider the smaller assortment $S_B=\{3\}$. The \textsc{MultiMNL} probabilities are $P_{\textsc{MultiMNL}}(3, \{3\}|1) = \frac{w_{1,3}}{w_{1,3}+1} = \frac{2}{3}$ and $P_{\textsc{MultiMNL}}(3, \{3\}|2) = \frac{w_{2,3}}{w_{2,3}+1} = \frac{1}{2}$.
For the \textsc{MarkovMNL} model to match these, its parameters must satisfy:
\begin{align*}
    P_{\textsc{MarkovMNL}}(3, \{3\}|1) &= \lambda_{1,3} + \lambda_{1,4} \frac{v_3^B}{v_3^B+1} = \frac{2}{3}, \\
    P_{\textsc{MarkovMNL}}(3, \{3\}|2) &= \lambda_{2,3} + \lambda_{2,4} \frac{v_3^B}{v_3^B+1} = \frac{1}{2}.
\end{align*}
Substituting the $\lambda$ values we derived:
$
    \frac{2}{5} + \frac{2}{5} \frac{v_3^B}{v_3^B+1} = \frac{2}{3},
    \frac{1}{3} + \frac{1}{3} \frac{v_3^B}{v_3^B+1} = \frac{1}{2}$.
The first equation yields $v_3^B=2$, whereas the second gives $v_3^B=1$, a contradiction. \hfill
 \Halmos
\endproof

\subsubsection*{Comparison with \cite{chen2022assortmentoptimizationmultivariatemnl}}
\begin{lemma}\label{prop:STC-MNL-Mv}
The \textsc{MarkovMNL} model and the \textsc{MvMNL} model are distinct; neither model class contains the other.
\end{lemma}
\proof{Proof of Lemma \ref{prop:STC-MNL-Mv}.}
(\emph{\textsc{MarkovMNL} $\not\subseteq$ \textsc{MvMNL}}). The key intuition here is that our model can violate the IIA property for bundles in \textsc{MvMNL}. Consider a universe with $\NN_A=\{1\}$ and $\NN_B=\{2,3\}$. Set the \textsc{MarkovMNL} parameters to $\lambda_{0_A,2}=\lambda_{0_A,3}=\lambda_{1,2}=\lambda_{1,3}=\frac{1}{3}$ and $v_{0_B}^B=v_2^B=1, v_3^B=2$. Let $S_A=\{1\}$. For simplicity, let $P_1=\phi_A(1, \{1\})$ denote the choice probability of product $1$ from category $A$.

\begin{itemize}
    \item When assortment $S_B=\{2,3\}$ is offered, the choice probabilities in \textsc{MarkovMNL} models are $P_{\textsc{MarkovMNL}}(1,2,S_A,S_B)=\lambda_{1,2}P_1=\frac{1}{3}P_1$ and $P_{\textsc{MarkovMNL}}(1,3,S_A,S_B)=\lambda_{1,3}P_1=\frac{1}{3}P_1$. The odds ratio of these joint probabilities is $1$. In an \textsc{MvMNL} model, this implies $u_{1,2}/u_{1,3} = 1$. Similarly, $u_{0_A,2}/u_{0_A,3} = 1$.
    \item When $S_B=\{2\}$, the joint probability is $P_{\textsc{MarkovMNL}}(1,2,S_A,S_B)= (\lambda_{1,2} + \lambda_{1,3} \frac{v_2^B}{v_2^B+1})P_1 = \frac{1}{2}P_1$. Similarly, when $S_B=\{3\}$, the joint probability is $P_{\textsc{MarkovMNL}}(1,3,S_A,S_B)= (\lambda_{1,3} + \lambda_{1,2} \frac{v_3^B}{v_3^B+1})P_1= (\frac{1}{3} + \frac{1}{3}  \frac{2}{3})P_1 = \frac{5}{9}P_1$. For an \textsc{MvMNL} model, $u_{1,2}=u_{1,3}$ and $u_{0_A,2}=u_{0_A,3}$ imply
    \[
    P_{\textsc{MvMNL}}(1,3,S_A,\{3\})
    =P_{\textsc{MvMNL}}(1,2,S_A,\{2\}).
    \]
    The \textsc{MarkovMNL} probabilities above are unequal, so this instance has no \textsc{MvMNL} representation.
\end{itemize}

(\emph{\textsc{MvMNL} $\not\subseteq$ \textsc{MarkovMNL}}).
Consider a universe with $\NN_A=\{ 1\}$ and $\NN_B=\{ 2, 3\}$. Let the \textsc{MvMNL} bundle preference weights be $u_{1,2}=2, u_{0_A,3}=2$, and set the preference weights of all other bundles to $1$.

First, we determine the $\bfl$ parameters that a \textsc{MarkovMNL} model would require to match the conditional probabilities when $S_A=\{1\}$ and $S_B=\{2,3\}$ are offered. We have the following choice probabilities from \textsc{MvMNL}:
\begin{align*}
    P_{\textsc{MvMNL}}(2,S_B|1) &= \frac{u_{1,2}}{u_{1,2}+u_{1,3}+u_{1,0_B}}  = \frac{1}{2}, \\
    P_{\textsc{MvMNL}}(3,S_B| 1) &= \frac{u_{1,3}}{u_{1,2}+u_{1,3}+u_{1,0_B}} = \frac{1}{4}, \\
   P_{\textsc{MvMNL}}(2,S_B|0_A) &= \frac{u_{0_A,2}}{u_{0_A,2}+u_{0_A,3}+u_{0_A,0_B}}= \frac{1}{4}, \\
   P_{\textsc{MvMNL}}(3, S_B|0_A) &= \frac{u_{0_A,3}}{u_{0_A,2}+u_{0_A,3}+u_{0_A,0_B}} = \frac{1}{2}.
\end{align*}

Any \textsc{MarkovMNL} representation would therefore require
$\lambda_{1,2}=\frac{1}{2}$, $\lambda_{1,3}=\frac{1}{4}$, $\lambda_{0_A,2}=\frac{1}{4}$, $\lambda_{0_A,3}=\frac{1}{2}$.
Now consider assortment $S_B=\{2\}$. From \textsc{MvMNL}, we have:
 \begin{align*}
      &P_{\textsc{MvMNL}}(2,  S_B|1)= \frac{u_{1,2}}{u_{1,0_B}+u_{1,2}}=\frac{2}{3}, \\
           &P_{\textsc{MvMNL}}(2, S_B|0_A)= \frac{u_{0_A,2}}{u_{0_A,0_B}+u_{0_A,2}}=\frac{1}{2}.
  \end{align*}
To match the \textsc{MarkovMNL} model, we need
 \begin{align*}
      &P_{\textsc{MarkovMNL}}(2,  S_B|1)= \lambda_{1,2}+\lambda_{1,3}\frac{v_2^B}{v_2^B+1}=\frac{1}{2}+\frac{1}{4}\frac{v_2^B}{v_2^B+1}= \frac{2}{3}, \\
      & P_{\textsc{MarkovMNL}}(2, S_B| 0_A)=\lambda_{0_A,2}+\lambda_{0_A,3}\frac{v_2^B}{v_2^B+1}=\frac{1}{4}+\frac{1}{2}\frac{v_2^B}{v_2^B+1}=\frac{1}{2}.
  \end{align*}
  No value of $v_2^B$ could satisfy both equations, so this \textsc{MvMNL} model cannot be represented by any \textsc{MarkovMNL} model.  \hfill
\Halmos
\endproof

Next, we show that \textsc{MvMNL} and \textsc{MultiMNL} cannot express one another.

\begin{lemma}\label{prop:mvmnl-multi}
The \textsc{MvMNL} model and the \textsc{MultiMNL} model are distinct; neither model class contains the other.
\end{lemma}
\proof{Proof of Lemma \ref{prop:mvmnl-multi}.} We give a counterexample for each direction.

(\emph{\textsc{MvMNL} $\not\subseteq$ \textsc{MultiMNL}}).
Let $\NN_A=\{1\}$ and $\NN_B=\{2,3\}$. In \textsc{MvMNL}, set the bundle preference weights as $u_{1,3}=3$ and $u_{i,j}=1$ for all other $(i,j)\in\{0_A,1\}\times\{0_B,2,3\}$.

First, with $S_A=\{1\}$ and $S_B=\{2\}$, the available bundles are
$\{(1,2),(1,0_B),(0_A,2),(0_A,0_B)\}$, each with weight $1$. The ratio of the choice probability of $1$ to that of $0_A$ from $A$ is:
\[
\frac{P_{\textsc{MvMNL}}(\text{choose }1)}{P_{\textsc{MvMNL}}(\text{choose }0_A)}
= \frac{u_{1,2}+u_{1,0_B}}{u_{0_A,2}+u_{0_A,0_B}}=1.
\]
Hence any matching \textsc{MultiMNL} must have the same odds $w_1/w_{0_A}=1$.

Next, with $S_A=\{1\}$ and $S_B=\{3\}$, the available bundles are
$\{(1,3),(1,0_B),(0_A,3),(0_A,0_B)\}$ with weights $\{3,1,1,1\}$. Thus the ratio of the choice probability of $1$ to that of $0_A$ from $A$ is:
\[
\frac{P_{\textsc{MvMNL}}(\text{choose }1)}{P_{\textsc{MvMNL}}(\text{choose }0_A)}
= \frac{u_{1,3}+u_{1,0_B}}{u_{0_A,3}+u_{0_A,0_B}}=2.
\]
Therefore a matching \textsc{MultiMNL} would need $w_1/w_{0_A}=2$, which contradicts $w_1/w_{0_A}=1$ above.

(\emph{\textsc{MultiMNL} $\not\subseteq$ \textsc{MvMNL}}).
Again let $\NN_A=\{1\}$ and $\NN_B=\{2,3\}$. In \textsc{MultiMNL}, take
\(
w_1=w_{0_A}=w_{0_A,2}=w_{1,3}=1,\,
w_{1,2}=w_{0_A,3}=2,
\)
and normalize the relevant no-purchase weights as $w_{1,0_B}=w_{0_A,0_B}=1$.

For $S_A=\{1\}$, $S_B=\{2\}$, the bundle probabilities are
\begin{align*}
P_{\textsc{MultiMNL}}(1,0_B,S_A,S_B)&=\frac{w_1}{w_1+w_{0_A}}\frac{1}{w_{1,2}+1}
=\frac{1}{2}\frac{1}{3}=\frac{1}{6},\\
P_{\textsc{MultiMNL}}(0_A,0_B,S_A,S_B)&=\frac{w_{0_A}}{w_1+w_{0_A}}\frac{1}{w_{0_A,2}+1}
=\frac{1}{2}\frac{1}{2}=\frac{1}{4}.
\end{align*}
Any \textsc{MvMNL} representation has fixed bundle weights $\{u_{i,j}\}$ across assortments, so it must satisfy
\[
\frac{u_{1,0_B}}{u_{0_A,0_B}}
=\frac{P_{\textsc{MultiMNL}}(1,0_B,S_A,S_B)}{P_{\textsc{MultiMNL}}(0_A,0_B,S_A,S_B)}
=\frac{1/6}{1/4}=\frac{2}{3}.
\]

For $S_A=\{1\}$, $S_B=\{3\}$, we obtain
\begin{align*}
P_{\textsc{MultiMNL}}(1,0_B,S_A,S_B)&=\frac{w_1}{w_1+w_{0_A}}\frac{1}{w_{1,3}+1}
=\frac{1}{2}\frac{1}{2}=\frac{1}{4},\\
P_{\textsc{MultiMNL}}(0_A,0_B,S_A,S_B)&=\frac{w_{0_A}}{w_1+w_{0_A}}\frac{1}{w_{0_A,3}+1}
=\frac{1}{2}\frac{1}{3}=\frac{1}{6},
\end{align*}
so
\[
\frac{u_{1,0_B}}{u_{0_A,0_B}}
=\frac{P_{\textsc{MultiMNL}}(1,0_B,S_A,S_B)}{P_{\textsc{MultiMNL}}(0_A,0_B,S_A,S_B)}
=\frac{1/4}{1/6}=\frac{3}{2},
\]
which contradicts the previous ratio $2/3$. Hence, no \textsc{MvMNL} model represents this \textsc{MultiMNL} instance. \hfill\Halmos
\endproof
\begin{lemma}The \textsc{MarkovMC} model class and the \textsc{MultiMNL} model class are distinct; neither model class contains the other.
    \label{prop: STC-MC-Multi}
\end{lemma}
\proof{Proof of Lemma \ref{prop: STC-MC-Multi}.} (\emph{\textsc{MarkovMC} $\not\subseteq$ \textsc{MultiMNL}}).
By Lemma \ref{prop:STC-multimnl} we already have \textsc{MarkovMNL} $\not\subseteq$ \textsc{MultiMNL}. Since \textsc{MarkovMNL} is a special case of \textsc{MarkovMC}, it follows that \textsc{MarkovMC} $\not\subseteq$ \textsc{MultiMNL}.

(\emph{$ \textsc{MultiMNL} \not\subseteq \textsc{MarkovMC}$}).
Let $\NN_A=\{1\}$ and $\NN_B=\{2,3\}$. Consider a \textsc{MultiMNL} with weights $\{w_1,w_{0_A}\}$ from $A$ and weights $\{w_{i,2},w_{i,3}\}$ from $B$ for $i\in\{1,0_A\}$, with the second-stage no-purchase weight normalized to $1$.

We consider two assortments in category $B$.

\noindent\underline{Assortment $S_B=\{2,3\}$.}
Under \textsc{MultiMNL},
\[
P_{\textsc{MultiMNL}}(2,S_B\mid i)=\frac{w_{i,2}}{w_{i,2}+w_{i,3}+1},
\qquad
P_{\textsc{MultiMNL}}(3,S_B\mid i)=\frac{w_{i,3}}{w_{i,2}+w_{i,3}+1}.
\]
In a \textsc{MarkovMC}, offered products are absorbing. Hence, with $S_B=\{2,3\}$, the probability of choosing $2$ equals the initial arrival mass on $2$, denoted $\lambda_{i,2}$, and similarly for $3$. Matching the two models gives
\[
\lambda_{i,2}=\frac{w_{i,2}}{w_{i,2}+w_{i,3}+1},
\qquad
\lambda_{i,3}=\frac{w_{i,3}}{w_{i,2}+w_{i,3}+1}.
\]
\noindent\underline{Assortment $S_B=\{2\}$.}
Under \textsc{MultiMNL},
\(
P_{\textsc{MultiMNL}}(2,S_B\mid i)=\frac{w_{i,2}}{w_{i,2}+1}.
\)
In a \textsc{MarkovMC}, with only $2$ offered, the probability of choosing $2$ is
\(
P_{\textsc{MarkovMC}}(2,S_B\mid i)=\lambda_{i,2}+\lambda_{i,3}\, h,
\)
where $h\in[0,1]$ is the probability, determined solely by $\bfrho$, that the chain starting from state $3$ hits $2$ before absorption at $0_B$. Crucially, $h$ does not depend on $i$.

Equating the two models and substituting the $\lambda$’s identified from $S_B=\{2,3\}$ yields, for each $i\in\{1,0_A\}$,
\[
\frac{w_{i,2}}{w_{i,2}+1}
=\frac{w_{i,2}}{w_{i,2}+w_{i,3}+1}
+ \frac{w_{i,3}}{w_{i,2}+w_{i,3}+1}\,h.
\]
Solving for $h$ gives
\(
h
=\frac{w_{i,2}}{w_{i,2}+1}\).
Thus $h$ must equal $w_{i,2}/(w_{i,2}+1)$ for both $i=1$ and $i=0_A$, which forces
\[
\frac{w_{1,2}}{w_{1,2}+1}=\frac{w_{0_A,2}}{w_{0_A,2}+1}
\quad\Longrightarrow\quad w_{1,2}=w_{0_A,2}.
\]
Choosing any \textsc{MultiMNL} with $w_{1,2}\neq w_{0_A,2}$ (and $w_{i,3}>0$) yields a contradiction. Hence \textsc{MultiMNL} is not contained in \textsc{MarkovMC}. \hfill\Halmos
\endproof
\begin{lemma}
The \textsc{MarkovMC} and \textsc{MvMNL} model classes are distinct; neither contains the other.
\label{prop:STC-MC-Mv}
\end{lemma}
\proof{Proof of Lemma \ref{prop:STC-MC-Mv}.}
(\emph{\textsc{MarkovMC} $\not\subseteq$ \textsc{MvMNL}}).
By Lemma \ref{prop:STC-MNL-Mv}, \textsc{MarkovMNL} $\not\subseteq$ \textsc{MvMNL}. Since \textsc{MarkovMNL} is a special case of \textsc{MarkovMC}, we conclude \textsc{MarkovMC} $\not\subseteq$ \textsc{MvMNL}.

(\emph{\textsc{MvMNL} $\not\subseteq$ \textsc{MarkovMC}}).
Let $\NN_A=\{1\}$ and $\NN_B=\{2,3\}$ with no-purchase options $0_A,0_B$. Consider an \textsc{MvMNL} with bundle weights
\(\{u_{i,j}\}_{i\in\{1,0_A\},\,j\in\{2,3,0_B\}}\) and \(u_{1,0_B}=u_{0_A,0_B}=1\).

\noindent \underline{Assortment \(S_B=\{2,3\}\).} In \textsc{MarkovMC}, absorption at offered items gives
\(
P_{\textsc{MarkovMC}}(2,S_B\mid i)=\lambda_{i,2},
P_{\textsc{MarkovMC}}(3,S_B\mid i)=\lambda_{i,3}.
\)
Matching \textsc{MvMNL}, we have:
\[
P_{\textsc{MvMNL}}(2,S_B\mid i)=\frac{u_{i,2}}{u_{i,2}+u_{i,3}+1},\quad
P_{\textsc{MvMNL}}(3,S_B\mid i)=\frac{u_{i,3}}{u_{i,2}+u_{i,3}+1},
\]
which forces
\[
\lambda_{i,2}=\frac{u_{i,2}}{u_{i,2}+u_{i,3}+1},\quad
\lambda_{i,3}=\frac{u_{i,3}}{u_{i,2}+u_{i,3}+1}.
\]
\noindent \underline{Assortment \(S_B=\{2\}\).}
Now
\(
P_{\textsc{MarkovMC}}(2,S_B\mid i)=\lambda_{i,2}+\lambda_{i,3}\,h\), \(
P_{\textsc{MvMNL}}(2,S_B\mid i)=\frac{u_{i,2}}{u_{i,2}+1}.
\)
Substituting the \(\lambda\)’s and solving for \(h\) yields, for each \(i\in\{1,0_A\}\),
\[
\frac{u_{i,2}}{u_{i,2}+1}
=\frac{u_{i,2}}{u_{i,2}+u_{i,3}+1}+\frac{u_{i,3}}{u_{i,2}+u_{i,3}+1}\,h
\;\Longrightarrow\;
h=\frac{u_{i,2}}{u_{i,2}+1}.
\]
Thus \(h\) must simultaneously equal \(u_{1,2}/(u_{1,2}+1)\) and \(u_{0_A,2}/(u_{0_A,2}+1)\). Choosing \(u_{1,2}\neq u_{0_A,2}\) (with \(u_{i,3}>0\)) yields a contradiction. Hence \textsc{MvMNL} $\not\subseteq$ \textsc{MarkovMC}. \hfill\Halmos
\endproof

\begin{lemma}\label{prop:MarkovMNL-MNL}
   \textsc{MarkovMNL} includes the independent MNL models as a special case.
\end{lemma}
\proof{Proof of Lemma \ref{prop:MarkovMNL-MNL}.}
We prove this in two parts: (i) any two independent MNL models can be represented by a \textsc{MarkovMNL} model; (ii) the inclusion is strict, since some \textsc{MarkovMNL} models cannot be represented by any two independent MNL models.

For part (i), consider independent MNL models for categories $A$ and $B$, with preference weights $\{u^A_i\}_{i\in \NN_A}$ and $\{u^B_j\}_{j\in \NN_B}$ and normalized no-purchase weights $u^A_{0_A}=u^B_{0_B}=1$. We construct an equivalent \textsc{MarkovMNL} model by setting $v^A_i=u^A_i$ for all $i\in \NN_A$ and $v^B_j=u^B_j$ for all $j\in \NN_B$. The key step is to define the transition probabilities as the MNL probabilities under the full assortment, independent of the first-stage choice $i$:
\[
\lambda_{i,j}=\frac{u_j^B}{1+\sum_{k\in\NN_B}u_k^B},
\qquad i\in\NN_A^+,\quad j\in\NN_B^+.
\]
Then, for any assortment $S_B\subseteq \NN_B$, the probability of choosing $j\in S_B$ in the \textsc{MarkovMNL} model, conditional on first-stage choice $i$, is
\begin{align*}
P_{\textsc{MarkovMNL}}(j,S_B\mid i)
 &= \lambda_{i,j} + \sum_{k\in \NN_B\setminus S_B}\lambda_{i,k}\frac{v^B_j}{1+\sum_{q\in S_B} v^B_q} \\
 &= \frac{u^B_j}{1+\sum_{m\in \NN_B} u^B_m} + \sum_{k\in \NN_B\setminus S_B}\frac{u^B_k}{1+\sum_{m\in \NN_B} u^B_m}\frac{u^B_j}{1+\sum_{q\in S_B} u^B_q} \\
 &= \frac{u^B_j}{1+\sum_{q\in S_B} u^B_q}.
\end{align*}
This matches the standard MNL probability.

For part (ii), we give a counterexample. The key insight is that \textsc{MarkovMNL} can violate the IIA property of standard MNL models. Consider $\NN_A=\{1\}$, $\NN_B=\{2,3\}$, with parameters $\lambda_{1,2}=\lambda_{1,3}=1/3$ and $v^B_{0_B}=v^B_2=1$, $v^B_3=2$. (Parameters related to the first-stage no-purchase option $0_A$ are omitted, as they are not needed here.)

Conditional on choosing product $1$ from $A$, the choice probabilities in $B$ are:
\begin{itemize}
    \item $S_B=\{2,3\}$: $P_{\textsc{MarkovMNL}}(2,S_B\mid 1)=\frac{1}{3}$, $P_{\textsc{MarkovMNL}}(3,S_B\mid 1)=\frac{1}{3}$.
    \item $S_B=\{2\}$: $P_{\textsc{MarkovMNL}}(2,S_B\mid 1)=\frac{1}{2}$.
    \item $S_B=\{3\}$: $P_{\textsc{MarkovMNL}}(3,S_B\mid 1)=\frac{5}{9}$.
\end{itemize}
For a standard MNL model to match these, its parameters $\{u^B_j\}_{j\in \NN_B}$ must satisfy:
\begin{itemize}
    \item $S_B=\{2,3\}$: $\frac{u^B_2}{u^B_2+u^B_3+1}=\frac{1}{3}$ and $\frac{u^B_3}{u^B_2+u^B_3+1}=\frac{1}{3}$, implying $u^B_2=u^B_3$.
    \item $S_B=\{2\}$: $\frac{u^B_2}{u^B_2+1}=\frac{1}{2}$, giving $u^B_2=1$.
    \item $S_B=\{3\}$: $\frac{u^B_3}{u^B_3+1}=\frac{5}{9}$, giving $u^B_3=\frac{5}{4}$.
\end{itemize}
These conditions imply both $u^B_2=u^B_3$ and $u^B_2\neq u^B_3$, a contradiction. Hence, no standard MNL model can represent this \textsc{MarkovMNL} instance. \hfill\Halmos
\endproof

\subsection[Connection to the Model of Cao et al. (2023)]{Connection to the Model of \cite{cao2023revenue}} \label{app:cao}
In our model, for any given assortment $S_A$, the induced choice model in category $B$ has the following choice probabilities for any product $j\in S_B$: $\sum_{i\in S_A^+}\frac{v_i^A}{V^A(S_A)+1}\left(\lambda_{i,j}+ \frac{v_j^B }{V^B(S_B) +1} \sum_{m\in \NN_B \setminus S_B}\lambda_{i,m}\right)$.
\cite{cao2023revenue} propose a single-purchase model that is a mixture of an independent demand model and an MNL model. In their framework, a customer belongs to one of two segments. With probability $q\in (0,1)$, a customer belongs to the first segment and chooses product $i\in S$ with a fixed probability $\theta_i$. With probability $1-q$, the customer is in the second segment and chooses according to a standard MNL model with weights $(v_j)_{j\in \NN_B}$. The resulting choice probability for a product $i\in S$ is
\(
\phi_{\textsc{mix}}(i,S)=q \theta_i+(1-q)\frac{v_i}{\sum_{k\in S}v_k+1}.
\)

The two constructions use the MNL component differently.
In our model, its coefficient is the probability of initially
preferring an unavailable product, and therefore depends on
the offered assortment. In the model of
\citet{cao2023revenue}, the coefficient is a fixed
customer-segment weight.

\subsection{Multi-Parent Categories: Pooled Arrivals}\label{app:multi-parent}
In this section, we formalize the choice process in a category with multiple parent categories. Let $\text{Pa}(W)$ denote the set of parents of category $W$, and let $i_U\in\NN_U^+$ denote the realized outcome in parent $U$ (a product or the no-purchase option). Categories are processed in a topological order, so all parental outcomes are realized before $W$ is reached. Each incoming edge $(U,W)$ carries a fixed weight $\omega_{U,W}\ge 0$ with $\sum_{U\in\text{Pa}(W)}\omega_{U,W}=1$; unless stated otherwise, we take $\omega_{U,W}=1/|\text{Pa}(W)|$, and $\omega_{U,W}=1$ when $U$ is the only parent. Given the parental outcomes $\mathbf{i}_{\text{Pa}(W)}=(i_U)_{U\in\text{Pa}(W)}$, the initial attraction in $W$ is drawn from the pooled arrival distribution
\begin{equation}\label{eq:pooled-arrival}
\psi^W_\ell\bigl(\mathbf{i}_{\text{Pa}(W)}\bigr) \;=\; \sum_{U\in\text{Pa}(W)} \omega_{U,W}\,\lambda^{U,W}_{i_U,\ell}, \qquad \ell\in\NN_W^+.
\end{equation}
The customer then makes exactly one choice in $W$, following the same mechanism as in the two-category model: if the attracted option $\ell$ is available ($\ell\in S_W^+$), it is chosen; otherwise, the customer substitutes according to $\phi_W(\cdot,S_W\mid\ell)$. The resulting choice probability of $j\in S_W^+$ is
\begin{equation}\label{eq:pooled-choice}
P_W\bigl(j, S_W \mid \mathbf{i}_{\text{Pa}(W)}\bigr) \;=\; \psi^W_j\bigl(\mathbf{i}_{\text{Pa}(W)}\bigr) + \sum_{\ell\in\NN_W\setminus S_W} \psi^W_\ell\bigl(\mathbf{i}_{\text{Pa}(W)}\bigr)\,\phi_W(j,S_W\mid \ell),
\end{equation}
which generalizes Equation~\eqref{eq:2-prob}. Since the weights and each transition row sum to one, \eqref{eq:pooled-choice} is a valid probability distribution over $S_W^+$: parental outcomes influence $W$ through the pooled arrival distribution, the customer makes exactly one choice in $W$ (possibly $0_W$), and no product can be purchased twice through different incoming edges. Equivalently, $\omega_{U,W}$ can be interpreted as the probability that the customer \emph{anchors} on parent $U$: a latent anchor parent is drawn with probabilities $(\omega_{U,W})_{U\in\text{Pa}(W)}$, and the attraction is then formed from the anchor's realized outcome, in the same spirit as the latent-direction mixture of Section~\ref{sec:bidirectional}. When $\phi_W$ is an MC model, Lemma~\ref{lem:choice-mc} extends verbatim: conditional on $\mathbf{i}_{\text{Pa}(W)}$, the choice process in $W$ is an MC model with arrival distribution $\psi^W(\mathbf{i}_{\text{Pa}(W)})$ and transition matrix $\bfrho^W$.

The weights can be fixed a priori or estimated: treating the anchor parent as an additional latent variable, the EM framework of Appendix~\ref{app:est-dag} yields closed-form updates for the weights from posterior responsibilities. We note that the weights are identified jointly with the edge transition matrices only up to transformations that leave the pooled arrival distributions \eqref{eq:pooled-arrival} unchanged; fixing the weights (e.g., uniform) avoids this indeterminacy. If two parameterizations induce the same pooled arrival distribution for every parental-outcome tuple in the full Cartesian set $\prod_{U\in\operatorname{Pa}(W)}\NN_U^+$, then each parent's edge contribution to its adjusted prices changes by a common constant across that parent's outcomes. The corresponding optimal assortment is therefore unchanged by Lemma~\ref{prop:3-opt-equ}. Under more limited support, this common-constant conclusion need not hold.

\section{Omitted Details and Proofs in Section \ref{sec:3-opt}}

\subsection{Existence of an Invariant Unconstrained Optimal Assortment Under an MC Model} \label{app:mc-lemma}
\begin{lemma}
    \label{lem:mc-invariance} For any Markov chain choice model $\phi$ defined by a price vector $\mathbf{r}$, transition matrix $\bfrho$, and initial arrival probabilities
$\bfpsi$, there exists an optimal unconstrained assortment under $\phi$ that depends only on $(\mathbf{r}, \bfrho)$ and is independent of $\bfpsi$.
\end{lemma}

\proof{Proof of Lemma \ref{lem:mc-invariance}.}
\cite{blanchet2016markov} show that the optimal revenue-maximizing assortment for a Markov chain choice model can be found by solving a linear program (LP) (see Theorem 5.1 and Lemma 5.2 in \cite{blanchet2016markov}). The objective function and constraints of this LP are constructed using only the price vector $\mathbf{r}$ and the transition matrix $\bfrho$.
Since the initial arrival probabilities $\bfpsi$ do not appear as parameters in this LP formulation, the solution, and thus the optimal assortment derived from it, is independent of $\bfpsi$.

Moreover, the iterative algorithm given by \cite{desir2020constrained} also yields an optimal assortment that is independent of the initial arrival probabilities (see Algorithm 1 and Theorem 2 in \cite{desir2020constrained}).
\hfill \Halmos
\endproof

\subsection{Existence of Non-invariant Optimal Assortments for MC}\label{app:no-mc}
\begin{example}[Non-invariant optimal assortments]\label{eg:no-mc}Not all optimal assortments for an MC model are invariant across different arrival distributions.
Consider a category with two products $\{1,2\}$. Let $0$ denote the outside option. The transition probabilities are
$\rho_{1,0}=\rho_{2,0}=1$, and the revenues are $r_1=r_2=1$.

If the arrival probabilities are $\psi_1=1,\;\psi_2=0$, an optimal assortment is
$S^*=\{1\}$. However, if the arrival probabilities switch to
$\psi_1=0,\;\psi_2=1$, then $\{1\}$ is no longer optimal. In this case,
 both $S^*=\{2\}$ and $S^*=\{1,2\}$ are optimal.

This example shows that the optimality of some assortments depends on the arrival distribution.
Nevertheless, the maximal assortment $S^*=\{1,2\}$ is optimal for both arrival probability
vectors, and is an invariant optimal assortment for this transition
matrix $\boldsymbol{\rho}$.
\end{example}
\subsection{Equivalence of Assortment Problems with a Non-Zero Price Outside Option}
\label{app:assortment-shift}In this section, we show that an unconstrained assortment optimization problem where the outside option has a non-zero price can be transformed into an equivalent problem where the outside option has a zero price. This equivalence allows standard assortment optimization algorithms, which often assume a zero-price outside option, to be applied directly.

Consider a ground set $\NN$ and an outside option denoted by $0$.
Let $\phi(i,S)$ be the probability of choosing option $i$ from $S^+:= S\cup\{0\}$. Let $(r_i)_{i \in \NN^+}$ be a set of prices, where the price of the outside option, $r_0$, may be positive. The objective is to find the assortment $S^* \subseteq \NN$ that solves:
\begin{equation} \max_{S \subseteq \NN} R(S) = \sum_{i \in S^+} r_i  \phi(i,S). \label{eq:non-zero-1}
\end{equation}

\begin{lemma} \label{prop:3-opt-equ}
Problems \eqref{eq:non-zero-1} and \eqref{eq:non-zero-2} have the same set of optimal assortments, where $r'_i=r_i-r_0$ for every $i\in\NN^+$:
\begin{equation}  \max_{S \subseteq \NN} R'(S) = \sum_{i \in S^+} r'_i  \phi(i,S) .
 \label{eq:non-zero-2}\end{equation}
 \end{lemma}
\proof{Proof of Lemma \ref{prop:3-opt-equ}.}
Let $c = r_0$ be the price of the outside option in the original problem. The shifted prices are defined as $r'_i = r_i - c$ for all $i \in \NN^+$. Note that for the outside option, the shifted price is $r'_0 = r_0 - c = r_0 - r_0 = 0$.

We can rewrite the original objective as:
\begin{align*}
R(S) &= \sum_{i \in S^+} r_i  \phi(i,S) = \sum_{i \in S^+} (r'_i + c)  \phi(i,S) = \sum_{i \in S^+} r'_i  \phi(i,S) + \sum_{i \in S^+} c  \phi(i,S) = R'(S) + c \sum_{i \in S^+} \phi(i,S) \\
& = R'(S)+c
\end{align*}
The last equation holds since $\phi(\cdot,S)$ is a probability distribution over the options in $S^+$, thus $\sum_{i \in S^+} \phi(i,S)=1$. Since $R(S)$ and $R'(S)$ differ only by the constant $c$, which is independent of the assortment $S$, both problems have the same set of optimal assortments.

The problem \eqref{eq:non-zero-2} is a standard assortment optimization problem where the outside option has zero price. The only difference is that now the products could have negative prices. For any substitution-based choice model, products with negative prices can be ignored without loss of generality. Therefore, after filtering out products where $r_i'<0$, problem \eqref{eq:non-zero-2} can be solved by existing algorithms. \hfill
\Halmos
\endproof
\subsection{Algorithms for DAGs}\label{app:alg-dag}
We first present the algorithm for a general DAG. The core idea is to use dynamic programming, processing the categories in a reverse topological order. As in the two-category case, for each category, we compute an adjusted price for each of its options, including the no-purchase option, whose adjusted price captures the expected continuation revenue when the customer leaves the category without a purchase. This adjusted price represents the option's intrinsic revenue plus its weighted expected future revenue from all downstream categories, where the pooling weights of Appendix~\ref{app:multi-parent} allocate the revenue of a multi-parent child exactly once across its parents.
\label{app:alg-dag-opt}
\begin{algorithm}[h!]\small
\caption{Unconstrained Assortment Optimization in a General DAG}
\label{alg:dag_assortment}
\SetAlgoNoLine
\textbf{Input:} A DAG $G=(\mathcal{V}, \mathcal{E})$ with pooling weights $(\omega_{U,W})_{(U,W)\in \mathcal{E}}$; for each category $U\in \mathcal{V}$, a ground set of products $\NN_U$, prices $r_{i_U}^U$ for each $i_U\in \NN_U$ (and $r_{0_U}^U=0$), and a choice model $\phi_U$\;
Let $(U_1, U_2, \dots, U_{|\mathcal{V}|})$ be a reverse topological ordering of the categories in $\mathcal{V}$\;
\For{$n \leftarrow 1$ \KwTo $|\mathcal{V}|$}
{
    Let $U \leftarrow U_n$ be the current category\;
    Define adjusted prices $(r'_{i_U})_{i_U \in \NN_U^+}$ for all options in $U$ (including $0_U$) as:
    $$r'_{i_U} \leftarrow r_{i_U}^U + \sum_{W \in \text{Children}(U)} \omega_{U,W} \Bigl( \sum_{i_W \in S_W^{*+}} \phi_{W|{i_U}}(i_W, S_W^*) r'_{i_W} \Bigr);$$
    Compute the optimal invariant unconstrained assortment $S_U^*$ for $U$ under $\phi_U$ using the shifted adjusted prices $(r'_{i_U}-r'_{0_U})_{i_U\in \NN_U^+}$ (cf.\ Lemma~\ref{prop:3-opt-equ})\;
}
\textbf{Return:} The set of optimal assortments $\{S_U^*\}_{U \in \mathcal{V}}$.
\end{algorithm}

To illustrate, consider a three-category chain $A \to B \to C$. The procedure begins at the terminal category $C$, where $S_C^*$ is computed using its original prices (so $r'_{i_C}=r_{i_C}^C$). It then moves to category $B$. The adjusted price for each product $i_B$ equals its original price plus the expected revenue from the fixed optimal assortment $S_C^*$, conditional on choosing $i_B$:
\[
r'_{i_B}=r_{i_B}^B+\sum_{i_C\in S_C^*} \phi_{C|i_B}(i_C,S_C^*)\,r_{i_C}^C.
\]
The optimal assortment $S_B^*$ is then computed using these forward-looking prices. Finally, at category $A$, each product’s adjusted price incorporates both its intrinsic price and the expected value of the entire $B\to C$ path, computed from the adjusted prices $r'_{i_B}$.

\begin{theorem}\label{thm:dag}
For a customer choice process on a DAG $G=(\mathcal{V},\mathcal{E})$, if the choice model $\phi_U$ for every category $U\in \mathcal{V}$ is a Markov chain choice model, then Algorithm \ref{alg:dag_assortment} computes an optimal set of assortments $\{S_U^*\}_{U\in \mathcal{V}}$ in polynomial time.
\end{theorem}

\proof{Proof of Theorem \ref{thm:dag}.}The proof has two parts: establishing optimality via backward induction and analyzing the computational complexity.

\emph{Edge-allocated decomposition.} Fix any assortments $\{S_V\}_{V\in \mathcal{V}}$ and let $i_V\in S_V^+$ denote the realized outcome in category $V$. The total expected revenue is $R(\{S_V\})=\sum_{V\in \mathcal{V}}\EE[r^V_{i_V}]$. Consider a non-root category $W$. Conditional on the realized parental outcomes, $i_W$ is generated by the MC model $\phi_W$ with the pooled arrival distribution $\psi^W(\mathbf{i}_{\text{Pa}(W)})$ of Equation~\eqref{eq:pooled-arrival}. Since the choice probabilities of an MC model are linear in its arrival distribution,
\[
\EE[r^W_{i_W}]=\sum_{U\in \text{Pa}(W)}\omega_{U,W}\,\EE\Bigl[\sum_{i_W\in S_W^+}\phi_{W|i_U}(i_W,S_W)\,r^W_{i_W}\Bigr],
\]
where $\phi_{W|i_U}$ denotes the MC model in $W$ with arrival distribution $\boldsymbol{\lambda}^{U,W}_{i_U}$ and the outer expectation is over the realized outcome $i_U$ in parent $U$. Because $\sum_{U\in \text{Pa}(W)}\omega_{U,W}=1$, the expected revenue of every category is allocated exactly once across its incoming edges. For the fixed assortments $\{S_V\}$, define value functions by the recursion
\[
v_{i_V}=r^V_{i_V}+\sum_{W\in \text{Children}(V)}\omega_{V,W}\sum_{i_W\in S_W^+}\phi_{W|i_V}(i_W,S_W)\,v_{i_W},\qquad i_V\in \NN_V^+,
\]
with $v_{i_K}=r^K_{i_K}$ for terminal categories $K$. Applying the allocation identity recursively down the DAG yields
\begin{equation}\label{eq:edge-alloc}
R(\{S_V\})=\sum_{U\in \mathcal{V}_r}\sum_{i_U\in S_U^+}\phi_U(i_U,S_U)\,v_{i_U},
\end{equation}
where $\mathcal{V}_r$ denotes the set of root categories.

\emph{Optimality.} We claim that for every category $V$, every $i_V\in \NN_V^+$, and arbitrary assortments of the descendants of $V$, we have $v_{i_V}\le r'_{i_V}$, where $r'_{i_V}$ is the adjusted price computed by Algorithm~\ref{alg:dag_assortment}, and that equality holds when every descendant $W$ uses $S_W^*$. We proceed by induction in reverse topological order. For terminal categories, $v_{i_K}=r^K_{i_K}=r'_{i_K}$. For the inductive step, consider a category $V$ whose children have all been processed, and fix any $i_V\in \NN_V^+$. Using the inductive hypothesis $v_{i_W}\le r'_{i_W}$ together with the nonnegativity of the weights and choice probabilities,
\[
v_{i_V}\le r^V_{i_V}+\sum_{W\in \text{Children}(V)}\omega_{V,W}\sum_{i_W\in S_W^+}\phi_{W|i_V}(i_W,S_W)\,r'_{i_W}.
\]
For each child $W$, the inner sum is the objective of the single-category MC assortment problem in $W$ with prices $(r'_{i_W})_{i_W\in \NN_W^+}$ and arrival distribution $\boldsymbol{\lambda}^{V,W}_{i_V}$. Shifting all of these prices by the constant $r'_{0_W}$ changes the objective of every assortment by exactly $r'_{0_W}$, because choice probabilities sum to one; hence, by Lemma~\ref{prop:3-opt-equ}, the shift affects neither the comparison between assortments nor the optimal assortment. By Lemma~\ref{lem:mc-invariance}, the invariant assortment $S_W^*$ computed by the algorithm is optimal for the shifted problem simultaneously for every arrival distribution, in particular for $\boldsymbol{\lambda}^{V,W}_{i_V}$ for every $i_V\in\NN_V^+$. Therefore,
\[
\sum_{i_W\in S_W^+}\phi_{W|i_V}(i_W,S_W)\,r'_{i_W}\le \sum_{i_W\in S_W^{*+}}\phi_{W|i_V}(i_W,S_W^*)\,r'_{i_W},
\]
and consequently $v_{i_V}\le r'_{i_V}$, with equality when all descendants use the algorithm's assortments. Combining the claim with \eqref{eq:edge-alloc}, for any assortments $\{S_V\}$,
\[
R(\{S_V\})\le \sum_{U\in \mathcal{V}_r}\sum_{i_U\in S_U^+}\phi_U(i_U,S_U)\,r'_{i_U}\le \sum_{U\in \mathcal{V}_r}\sum_{i_U\in S_U^{*+}}\phi_U(i_U,S_U^*)\,r'_{i_U}=R(\{S_V^*\}),
\]
where the second inequality holds because the algorithm solves each root category's unconstrained problem optimally (again using the price shift of Lemma~\ref{prop:3-opt-equ}), and the final equality applies the claim with equality. Thus, the returned set $\{S_U^*\}_{U \in \mathcal{V}}$ is globally optimal.

We next analyze the runtime of the algorithm.
\begin{enumerate}
    \item A reverse topological sort of the DAG can be performed in $O(|\mathcal{V}| + |\mathcal{E}|)$ time, for example, using Kahn's algorithm \citep{kahn1962}.
    \item The main loop iterates $|\mathcal{V}|$ times, once for each category. Inside the loop for a category $U$, computing the adjusted prices for all its products takes time proportional to $|\NN_U| \sum_{W \in \text{Children}(U)} |\NN_W|$. The final step inside the loop is solving a single-category assortment problem for an MC model. This problem can be solved in polynomial time in the number of products $|\NN_U|$ (e.g., via the linear program formulation in \cite{blanchet2016markov} or the iterative algorithm in \cite{desir2020constrained}).
\end{enumerate}
Since each step is polynomial in the input size (number of categories, edges, and products per category), the total runtime of the algorithm is polynomial.
\hfill
\Halmos
\endproof

\subsection{Proof of Theorem \ref{prop:car}}\label{app:root-card-proof}
\proof{Proof of Theorem \ref{prop:car}.}
Let $K_U$ denote the cardinality upper bound for each root category
$U\in\mathcal V_r$. Consider the following algorithm. Process the categories in
reverse topological order, as in Algorithm~\ref{alg:dag_assortment}. For every
non-root category $W\notin\mathcal V_r$, compute the invariant optimal
unconstrained assortment $S_W^{\ALG}$ under the MC choice model using the
current adjusted prices (shifted so that the no-purchase option has price
zero, cf.\ Lemma~\ref{prop:3-opt-equ}). For each root category $U\in\mathcal V_r$, after all
of its descendants have been processed, define the adjusted price of each
$i_U\in\NN_U^+$ by
\[
r'_{i_U}
=
r^U_{i_U}
+
\sum_{W\in \mathrm{Children}(U)}
\omega_{U,W}
\left(
\sum_{i_W\in S_W^{\ALG,+}}
\phi_{W|i_U}(i_W,S_W^{\ALG})\,r'_{i_W}
\right).
\]
Then solve the cardinality-constrained single-category problem
\[
\max_{S_U\subseteq \NN_U:\ |S_U|\le K_U}
\sum_{i_U\in S_U^+}
\phi_U(i_U,S_U)\,r'_{i_U}
\]
using the assumed polynomial-time $\alpha_U$-approximation algorithm. Let
$S_U^{\ALG}$ denote the assortment returned for root category $U$. The algorithm
returns $\{S_U^{\ALG}\}_{U\in\mathcal V_r}$ for the root categories and
$\{S_W^{\ALG}\}_{W\notin\mathcal V_r}$ for the non-root categories.

We first show that, once the root assortments are fixed, the assortments chosen
by the algorithm for the non-root categories are optimal. This follows from the
same backward-induction argument as in Theorem~\ref{thm:dag}. For any terminal
non-root category, the adjusted prices coincide with the original prices, so the
computed invariant MC assortment is optimal. Now suppose that all descendants of
a non-root category $W$ have already been assigned optimal continuation
assortments and correct adjusted prices. Conditional on any realized outcomes in its
parent categories, the choice in $W$ changes only through the pooled initial arrival
distribution \eqref{eq:pooled-arrival} into the MC model of $W$; the transition matrix and adjusted prices
remain fixed. By Lemma~\ref{lem:mc-invariance}, there exists an optimal
unconstrained assortment for $W$ that is invariant to the initial arrival
distribution. Hence the assortment $S_W^{\ALG}$ computed by the algorithm is
optimal for every possible parent choice. Therefore, by induction, the
algorithm computes optimal continuation assortments and correct adjusted prices
for every non-root category.

It remains to analyze the root categories. By the preceding paragraph, after
the non-root categories have been processed, the entire downstream value
generated by choosing product $i_U$ in root category $U$ is exactly captured by
the adjusted price $r'_{i_U}$. Thus the remaining optimization over each root
category is precisely the single-category cardinality-constrained problem
\(
\operatorname{OPT}_U
=
\max_{S_U\subseteq \NN_U:\ |S_U|\le K_U}
\sum_{i_U\in S_U^+}
\phi_U(i_U,S_U)\,r'_{i_U} .
\)
By the edge-allocated decomposition \eqref{eq:edge-alloc} established in the proof of
Theorem~\ref{thm:dag}, each category's expected revenue is allocated exactly once
across its incoming edges because $\sum_{U\in \text{Pa}(W)}\omega_{U,W}=1$, and the
accounts of distinct root categories are optimized independently---the total optimal value
of the multi-category problem with cardinality constraints only on the roots is
therefore
\(
\operatorname{OPT}
=
\sum_{U\in\mathcal V_r} \operatorname{OPT}_U .
\)
For each root category $U$, write $F_U(S_U)=\sum_{i_U\in S_U^+}\phi_U(i_U,S_U)\,r'_{i_U}$,
so that $\operatorname{OPT}_U=\max_{S_U\subseteq\NN_U:\,|S_U|\le K_U}F_U(S_U)$. The assumed
approximation algorithm returns $S_U^{\ALG}$ satisfying
\(
F_U(S_U^{\ALG}) \ge \alpha_U \operatorname{OPT}_U .
\)
Let $\alpha=\min_{U\in\mathcal V_r}\alpha_U$. The total revenue obtained by the
algorithm is
\[
R^{\ALG}
=
\sum_{U\in\mathcal V_r} \sum_{i_U\in S_U^{\ALG,+}}
\phi_U(i_U,S_U^{\ALG})\,r'_{i_U}
\ge
\sum_{U\in\mathcal V_r} \alpha_U \operatorname{OPT}_U
\ge
\alpha \sum_{U\in\mathcal V_r} \operatorname{OPT}_U
=
\alpha \operatorname{OPT} .
\]
Hence the algorithm achieves approximation factor
$\min_{U\in\mathcal V_r}\alpha_U$.

Finally, the algorithm runs in polynomial time. The reverse topological order
can be computed in polynomial time. For every non-root category, the
unconstrained MC assortment problem is solved in polynomial time by
Lemma~\ref{lem:mc-invariance} and the algorithms cited there. For every root
category, the assumed $\alpha_U$-approximation algorithm runs in polynomial
time. The adjusted-price computations require only polynomially many evaluations
of choice probabilities and revenue sums. Therefore, the overall algorithm is
polynomial-time.
\Halmos
\endproof
\subsection{Proof of Theorem \ref{prop:car-hard}}\label{app:prop-car-hard}
To show the inapproximability, we perform a reduction from the Bipartite Densest $\kappa$-Subgraph (BD$\kappa$S) problem. We first define the problem and state a known hardness result from \cite{chen2022assortmentoptimizationmultivariatemnl}.

\begin{definition}
    Given an undirected bipartite graph $G=(\NN,\MM,\bfE)$ with vertex sets $\NN$ and $\MM$ and edge set $\bfE\subseteq \NN\times \MM$, the BD$\kappa$S problem is to find subsets $\NN_1\subseteq \NN$ and $\MM_1\subseteq \MM$ that solve:
    \begin{align*}
        & \max_{\NN_1,\MM_1} |E(\NN_1\times \MM_1)| \\
        & \text{s.t.} \quad |\NN_1|=|\MM_1| = \kappa\,, \quad \NN_1\subseteq \NN\,, \quad \MM_1\subseteq \MM,
    \end{align*}
    where $E(\NN_1\times \MM_1)$ is the set of edges with one endpoint in $\NN_1$ and the other in $\MM_1$.
\end{definition}

 \begin{lemma}[Lemma 6 in \cite{chen2022assortmentoptimizationmultivariatemnl}]
    There is a constant $c > 0$ such that, assuming the Exponential Time Hypothesis (ETH), no polynomial-time algorithm can achieve an $\Omega\!\left(g^{-1/(\log \log g)^c}\right)$-approximation
for BD$\kappa$S on bipartite graphs with $g$ vertices.
\end{lemma}

Next, we prove the hardness of our cardinality-constrained problem.

\proof{Proof of Theorem \ref{prop:car-hard}.}
We reduce an arbitrary instance of BD$\kappa$S to our cardinality-constrained assortment problem. Consider a bipartite graph $G=(\NN_A, \NN_B, \bfE)$ and a cardinality $\kappa$. Let $m=|\NN_A|+|\NN_B|$.

We construct a two-category assortment instance where categories $A$ and $B$ correspond to the vertex sets $\NN_A$ and $\NN_B$. We use the MNL choice model for both categories with the following parameters:
\begin{itemize}
    \item For category $A$: $r_i^A=0$ and $v_i^A=1/m$ for all $i\in \NN_A$.
    \item For category $B$: $r_j^B = 1$ and $v_j^B=0$ for all $j\in \NN_B$.
    \item Transition probabilities from $i \in \NN_A$: $\lambda_{i,j} = \frac{1}{m}\1_{(i,j)\in \bfE}$ for $j \in \NN_B$, and the no-purchase transition is $\lambda_{i,0_B}=1-\sum_{j\in \NN_B} \lambda_{i,j} > 0$.
    \item Following no purchase in category $A$, set $\lambda_{0_A,0_B}=1$ and $\lambda_{0_A,j}=0$ for every $j\in\NN_B$.
\end{itemize}
The cardinality constraints for the assortment problem are set to $K_A=K_B=\kappa$. The expected revenue $R(S_A,S_B)$ for assortments $S_A \subseteq \NN_A$ and $S_B \subseteq \NN_B$ simplifies to:
\[
R(S_A,S_B) = \sum_{i \in S_A} \frac{v_i^A}{1 + \sum_{k \in S_A} v_k^A} \left( \sum_{j \in S_B} \lambda_{i,j} r_j^B \right) = \frac{\frac{1}{m^2} |E(S_A\times S_B)|}{1 + |S_A|/m}.
\]
This revenue function is bounded as follows for any feasible $S_A, S_B$ (where $|S_A|\le\kappa, |S_B|\le\kappa$):
\begin{align}
    R(S_A,S_B) &\le \frac{1}{m^2} |E(S_A\times S_B)| \label{prop1:eq2} ,\\
    R(S_A,S_B) &\ge \frac{1}{2m^2} |E(S_A\times S_B)| \label{prop1:eq1}.
\end{align}
Inequality \eqref{prop1:eq2} holds because the denominator $1+|S_A|/m \ge 1$. Inequality \eqref{prop1:eq1} holds because $|S_A| \le \kappa \le |\NN_A| < m$, which implies the denominator $1+|S_A|/m < 2$.

Now, suppose we have an $\alpha$-approximation algorithm for the assortment problem, which returns a solution $(\tilde{S}_A, \tilde{S}_B)$. We have:
\begin{align*}
|E(\tilde{S}_A\times \tilde{S}_B)| &= m^2 \left(1 + \frac{|\tilde{S}_A|}{m}\right) R(\tilde{S}_A, \tilde{S}_B) \ge m^2 R(\tilde{S}_A, \tilde{S}_B) \\
&\ge \alpha m^2 \max_{S_A, S_B} R(S_A,S_B) \\
&\ge \alpha m^2 \left( \frac{1}{2m^2} \max_{S_A, S_B} |E(S_A\times S_B)| \right) = \frac{\alpha}{2} \max_{S_A, S_B} |E(S_A\times S_B)|.
\end{align*}
The first inequality holds because $1 + |\tilde{S}_A|/m \ge 1$. The final inequality uses the bound from \eqref{prop1:eq1} applied to the optimal solution.

The solution $(\tilde{S}_A, \tilde{S}_B)$ is feasible for the assortment problem, meaning $|\tilde{S}_A| \le \kappa$ and $|\tilde{S}_B| \le \kappa$. The BD$\kappa$S problem requires sets of size exactly $\kappa$. Since adding vertices to $\tilde{S}_A$ or $\tilde{S}_B$ can only increase the edge count, we can arbitrarily add vertices to $\tilde{S}_A$ and $\tilde{S}_B$ to form sets $\tilde{S}'_A$ and $\tilde{S}'_B$ with $|\tilde{S}'_A| = |\tilde{S}'_B| = \kappa$. This new solution $(\tilde{S}'_A, \tilde{S}'_B)$ is feasible for BD$\kappa$S and satisfies:
\[
|E(\tilde{S}'_A\times \tilde{S}'_B)| \ge |E(\tilde{S}_A\times \tilde{S}_B)| \ge \frac{\alpha}{2} \max_{|S_A|=\kappa, |S_B|=\kappa} |E(S_A\times S_B)|.
\]
Thus an $\alpha$-approximation for the assortment problem yields an
$(\alpha/2)$-approximation for BD$\kappa$S. The reduction preserves the
size parameter: the total number of products is $g=m$. Therefore, an
$\Omega(g^{-1/(\log\log g)^c})$-approximation for the assortment problem
would contradict the stated hardness of BD$\kappa$S under ETH.
\hfill \Halmos
  \hfill
\endproof
\section{Omitted Details and Proofs in Section \ref{sec:estimation}}\label{app:estimation_details}
\subsection{Example of Non-Concavity of the Log-Likelihood}

The observed-data log-likelihood function in Equation~\eqref{eq:obs_log_likelihood} is generally not concave. We show this using a simple example. Consider ground sets $\NN_A = \{i\}$ and $\NN_B=\{j,m\}$.
Suppose the purchase history $\mathcal{H}$ consists of a single observation: at time $t=1$, the offered assortments were $S_A^1=\{i\}$ and $S_B^1=\{j\}$, and the customer purchased $i$ and $j$. Then the log-likelihood function is
$
\log L(\bfl, \bfv^A, \bfv^B;\mathcal{H})
= \log \frac{v_i^A}{v_i^A+1}
+ \log\!\left(\lambda_{i,j}+ \lambda_{i,m}\frac{v_j^B}{v_j^B+1}\right).
$
The first term corresponds to the choice of $i$ from $S_A^1$, and the second term corresponds to the subsequent purchase of $j$ from $S_B^1$.

If we reparameterize $v_i^A = \exp(\mu_i)$, the first term, $\log(\frac{\exp(\mu_i)}{1+\exp(\mu_i)}) = \mu_i - \log(1+\exp(\mu_i))$, is a standard logistic log-likelihood, which is concave in $\mu_i$. Since this term is independent of the second, the overall function's concavity depends entirely on the second term.

Let us analyze the function associated with the second term. For notational simplicity, write the parameter vector as $\boldsymbol{\theta} = (\lambda_j, \lambda_m, v_j)$ and define the function:
\(
f(\boldsymbol{\theta}) = f(\lambda_j, \lambda_m, v_j) = \log\!\left(\lambda_{j}+\lambda_{m}\frac{v_j}{1+v_j}\right).
\)
The following numerical example shows that $f$ is not concave. Let the two parameter vectors be:
\(
\boldsymbol{\theta}_1 = (0.2611, 0.5875, 7.5261)\) and \(
\boldsymbol{\theta}_2 = (0.6689, 0.0783, 2.4141).
\)
Evaluating the function at these points gives:
\begin{align*}
    f(\boldsymbol{\theta}_1) &= \log\!\left(0.2611 + 0.5875  \frac{7.5261}{1+7.5261}\right) \approx -0.2489, \\
    f(\boldsymbol{\theta}_2) &= \log\!\left(0.6689 + 0.0783  \frac{2.4141}{1+2.4141}\right) \approx -0.3226.
\end{align*}
The average of the endpoint values is $\tfrac{1}{2}[f(\boldsymbol{\theta}_1)+f(\boldsymbol{\theta}_2)]\approx -0.2857$.

Now, we evaluate the function at the midpoint of the parameter vectors:
\(
\boldsymbol{\theta}' = \tfrac{1}{2}(\boldsymbol{\theta}_1+\boldsymbol{\theta}_2) = (0.4650, 0.3329, 4.9701).
\)
This yields:
\[
f(\boldsymbol{\theta}') = \log\!\left(0.4650 + 0.3329  \frac{4.9701}{1+4.9701}\right) \approx -0.2982.
\]
Since $f(\boldsymbol{\theta}')\approx -0.2982<-0.2857
\approx \tfrac12[f(\boldsymbol{\theta}_1)+f(\boldsymbol{\theta}_2)]$,
the concavity inequality is violated. This demonstrates that the observed-data log-likelihood is not, in general, a concave function of the model parameters.

\subsection{EM Algorithm Details}\label{app:em-alg}

\noindent \textbf{E-step: Computing the expectation of latent variables.}
In the expectation step of iteration $\ell$, for each $m\in\mathcal M_t$ we compute the conditional expectation of the latent variable given the observed data $\mathcal{H}$ and current parameter estimates $\Theta^{(\ell)}$:
\[
\hat{X}_{m}^{t,(\ell)} := \EE[X_m^t \mid \mathcal{H}, \Theta^{(\ell)}]
= P(X_{m}^{t}=1 \mid a^t, b^t, S_A^t, S_B^t, \bfl^{(\ell)}, \bfv^{B,(\ell)}).
\]
Using the definition of conditional probability, this can be written as:
\[
\hat{X}_{m}^{t,(\ell)}
= \frac{P(X_{m}^{t}=1 \text{ and the downstream choice is } b^t \mid a^t, S_A^t, S_B^t, \Theta^{(\ell)})}
      {P(\text{the downstream choice is } b^t \mid a^t, S_A^t, S_B^t, \Theta^{(\ell)})}.
\]
The denominator is the total probability of choosing outcome $b^t$ after choosing $a^t$, which, from Lemma~\ref{lem:choice_prob}, is
\[
P(b^t \mid a^t, S_A^t, S_B^t , \bfl^{(\ell)}, \bfv^{B,(\ell)})
= \lambda_{a^t,b^t}^{(\ell)}
  + \frac{v_{b^t}^{B,(\ell)}}{1 + \sum_{k\in S_B^t} v_k^{B,(\ell)}}
    \sum_{j\in \NN_B \setminus S_B^t} \lambda_{a^t,j}^{(\ell)}.
\]
The numerator is the probability that the initial downstream state is $m$ and the final downstream choice is $b^t$. This occurs in two mutually exclusive ways:
\begin{enumerate}
    \item The customer directly chooses $m$ (which requires $m=b^t$), with probability $\lambda_{a^t,m}^{(\ell)}$.
    \item The customer is initially attracted to $m \in \NN_B \setminus S_B^t$, and the customer then chooses outcome $b^t$ from the offered assortment $S_B^t$. This occurs with probability $\lambda_{a^t,m}^{(\ell)}  \frac{v_{b^t}^{B,(\ell)}}{1+\sum_{k\in S_B^t} v_k^{B,(\ell)}}$.
\end{enumerate}
Combining these cases into a single expression for the numerator and dividing by the denominator gives the final update rule:
\begin{equation} \label{eq:e_step_Xim_appendix_final_revised}
\hat{X}_{m}^{t,(\ell)} =
\frac{\lambda_{a^t,m}^{(\ell)} \left( \1\{b^t=m\}
+ \1\{m \in \NN_B \setminus S_B^t\} \frac{v_{b^t}^{B,(\ell)}}{1+\sum_{k \in S_B^t} v_k^{B,(\ell)}} \right)}
{\lambda_{a^t,b^t}^{(\ell)}
+ \frac{v_{b^t}^{B,(\ell)}}{1+\sum_{k\in S_B^t} v_k^{B,(\ell)}}
  \sum_{j\in \NN_B \setminus S_B^t} \lambda_{a^t,j}^{(\ell)}}.
\end{equation}
The denominator is strictly positive on $\mathcal P$ because
$\lambda_{a^t,b^t}^{(\ell)}\ge\varepsilon_\lambda$.

\noindent \textbf{M-step: Maximization.}
In the maximization step, we maximize the expectation of the complete-data log-likelihood, $Q(\Theta|\Theta^{(\ell)}) = \EE_{\mathbf{X}|\mathcal{H},\Theta^{(\ell)}}[\log L_C(\Theta; \mathcal{H},\mathbf{X})]$. This objective decomposes into three separate terms:
\[
Q(\Theta|\Theta^{(\ell)}) = Q_1 (\bfv^A|\Theta^{(\ell)}) + Q_2 (\bfl|\Theta^{(\ell)}) + Q_3 (\bfv^B|\Theta^{(\ell)}),
\]
where, after substituting the expectations $\hat{X}_{m}^{t,(\ell)}$ for the latent variables $X_m^t$:
\begin{align*}
Q_1 (\bfv^A|\Theta^{(\ell)}) &= \sum_{t=1}^T
\log\!\left( \frac{v_{a^t}^A}{1 + V^A(S_A^t)} \right), \\
Q_2 (\bfl|\Theta^{(\ell)}) &= \sum_{t=1}^T \sum_{m \in \mathcal M_t} \hat{X}_{m}^{t,(\ell)} \log \lambda_{a^t ,m}, \\
Q_3 (\bfv^B|\Theta^{(\ell)}) &= \sum_{t=1}^T\sum_{m \in \mathcal M_t} \hat{X}_{m}^{t,(\ell)}
\log \!\left(
\1\{m = b^t\}
+ \1\{m \in \NN_B\setminus S_B^t\}
\frac{v_{b^t}^B}{1+V^B(S_B^t)}
\right).
\end{align*}

These three subproblems can be maximized independently. Write
$v_i^A=e^{\alpha_i}$ and $v_j^B=e^{\beta_j}$. Then $Q_1$ and $Q_3$ are
concave in $\boldsymbol\alpha$ and $\boldsymbol\beta$, respectively, over
the boxes $[\log\varepsilon_v,\log C]^{n_A}$ and
$[\log\varepsilon_v,\log C]^{n_B}$, while $Q_2$ is concave in $\bfl$ over the parameter space.
Thus, each subproblem can be solved efficiently by convex optimization.

The $\bfl$-subproblem also has an explicit row-wise solution. Define the
posterior counts
\(
c_{i,m}^{(\ell)}
:=
\sum_{\substack{t:a^t=i\\m\in\mathcal M_t}}
\hat X_m^{t,(\ell)}.
\)
If $|\{t:a^t=i\}|>0$, the exact constrained update is
\(
\lambda_{i,m}^{(\ell+1)}
=
\max\left\{\varepsilon_\lambda,
\frac{c_{i,m}^{(\ell)}}{\nu_i}\right\}, m\in\NN_B^+\),
where $\nu_i>0$ is chosen so that
$\sum_{m\in\NN_B^+}\lambda_{i,m}^{(\ell+1)}=1$. If no transaction has
$a^t=i$, the likelihood does not depend on row $\boldsymbol\lambda^i$, and we
retain its previous feasible value.

\subsection{Proof of Theorem \ref{thm:em-conv}}\label{app:thm-em-conv}

\proof{Proof of Theorem \ref{thm:em-conv}.}
The set $\mathcal P$ is compact and lies within the strictly positive
parameter space. Consequently, every observed choice has positive probability,
and $\mathcal{L}(\Theta)$ is finite and continuously differentiable on an open
neighborhood of $\mathcal P$. For each $m\in\mathcal M_t$, the posterior
expectation $\hat X_m^t$ is continuous in the current parameter vector because
its denominator is bounded away from zero. Hence
$Q(\Theta\mid\widetilde\Theta)$ is continuous in both arguments on
$\mathcal P\times\mathcal P$.

Each iteration exactly maximizes $Q$ over $\mathcal P$ by assumption, so the
standard EM inequality gives
$\mathcal{L}(\Theta^{(\ell+1)})\ge\mathcal{L}(\Theta^{(\ell)})$. Since
$\mathcal{L}(\Theta)\le0$, the log-likelihood sequence converges to a finite limit.
The constrained-parameter-space EM convergence result of
\citet{nettleton1999convergence} then implies that every accumulation point of
the parameter sequence is a constrained stationary point of $\mathcal{L}$ on
$\mathcal P$.\hfill \Halmos
\endproof

\subsection{Identifiability} \label{app:identifiability}
In category~$A$, the model reduces to a standard
single-category choice model: MNL for \textsc{MarkovMNL} and MC for
\textsc{MarkovMC}. Hence, the usual normalization and support conditions for
these models apply.

For category~$B$, the relevant object is the conditional downstream choice
probability
\(
q_i(j,S_B):=\mathbb P(b=j\mid a=i,S_B), j\in S_B^+,
\)
for each first-stage outcome $i\in\NN_A^+$ that is observed with positive
probability. Under our model,
\[
q_i(j,S_B)
=
\lambda_{i,j}
+
\sum_{\ell\in\NN_B\setminus S_B}
\lambda_{i,\ell}\phi_B(j,S_B\mid \ell).
\]
Thus, the first-stage choice $i$ determines the initial arrival distribution
$\boldsymbol{\lambda}^i$ into category~$B$, while the downstream substitution
model $\phi_B$ is shared across first-stage choices.

A simple sufficient condition for identifying $\boldsymbol{\lambda}$ is that
the full downstream assortment $S_B=\NN_B$ is observed for each relevant
first-stage choice $i$. In that case, no downstream substitution occurs, and
therefore
\(
q_i(j,\NN_B)=\lambda_{i,j}, j\in\NN_B^+.
\)
Hence, the $i$th row of $\boldsymbol{\lambda}$ is directly identified from the
conditional distribution of category~$B$ choices given $a=i$.

Once $\boldsymbol{\lambda}$ is identified, the remaining identification problem
in category~$B$ reduces to a standard single-category identification problem.
For \textsc{MarkovMNL}, the downstream model is an MNL model with the usual
normalization $v_{0_B}^B=1$. In particular, for any reduced assortment $S_B$ such
that
\(
L_i(S_B):=\sum_{\ell\in\NN_B\setminus S_B}\lambda_{i,\ell}>0,
\)
we have, for each $j\in S_B$,
\(
q_i(j,S_B)-\lambda_{i,j}
=
\frac{v_j^B}{1+\sum_{k\in S_B}v_k^B}L_i(S_B),
\)
and
\[
q_i(0_B,S_B)-\lambda_{i,0_B}
=
\frac{1}{1+\sum_{k\in S_B}v_k^B}L_i(S_B).
\]
Taking ratios yields
\(
v_j^B
=
\frac{
q_i(j,S_B)-\lambda_{i,j}
}{
q_i(0_B,S_B)-\lambda_{i,0_B}
}.
\)
Therefore, the downstream MNL weights are identified if, for every product
$j\in\NN_B$, there exists a pair $(i,S_B)$ such that
$q_i(\cdot,S_B)$ is observed, $j\in S_B$, and $L_i(S_B)>0$.

For \textsc{MarkovMC}, after $\boldsymbol{\lambda}$ is fixed, category~$B$
becomes a standard MC choice model with common transition matrix
$\boldsymbol{\rho}^B$. A sufficient condition is that, for every
$m\in\NN_B$, there exists an upstream outcome $i$ with $\lambda_{i,m}>0$
for which the conditional choice probabilities under
$S_B=\NN_B\setminus\{m\}$ are observed. For any such pair $(i,m)$, under
the standard normalization $\rho^B_{m,m}=0$,
\(
q_i(j,\NN_B\setminus\{m\})
=
\lambda_{i,j}
+
\lambda_{i,m}\rho^B_{m,j},
j\in\NN_B^+\setminus\{m\}.
\)
Thus,
\[
\rho^B_{m,j}
=
\frac{
q_i(j,\NN_B\setminus\{m\})-\lambda_{i,j}
}{
\lambda_{i,m}
}.
\]
Repeating this argument over $m$ identifies the transition matrix
$\boldsymbol{\rho}^B$. This is the same identification logic as in the standard
MC choice model of \citet{blanchet2016markov}. More recent work by
\citet{gupta2020parameter} further shows that MC choice parameters can be
identified from substantially smaller assortments, namely assortments of sizes
two and three. Therefore, our downstream identification requirement can be
viewed as a direct application of standard MC identification conditions to the
conditional choice probabilities $q_i(j,S_B)$.

If these assortment support conditions fail, the data may identify only the
induced probabilities $q_i(j,S_B)$ over the observed assortments, rather than a
unique primitive decomposition into cross-category arrival probabilities and
downstream substitution parameters. This is a support limitation inherited from
the underlying single-category choice model, not a separate issue created by
our cross-category construction. For estimation, our EM procedure follows the
standard missing-data logic used for MC choice models, where the unobserved
Markov path is treated as latent; see \citet{csimcsek2018expectation}.

\subsection{EM Algorithm for the Two-Category MarkovMC Model}
\label{app:markovmc-em}

In this section, we provide the full EM algorithm for the two-category
MarkovMC model used in the experiments. The derivation follows the EM approach
for the Markov chain choice model in \cite{csimcsek2018expectation}, adapted to
our two-category model with choice-dependent downstream arrival distributions.

Recall that in the MarkovMC specification, the within-category choice model in
category $A$ is an MC model with parameters $(\psi^A,\rho^A)$, while the
within-category choice model in category $B$ has a shared transition matrix
$\rho^B$. After observing the choice $i\in \NN_A^+$ in category $A$, the initial
arrival distribution into category $B$ is $\lambda^i=(\lambda_{i,\ell})_{\ell
\in \NN_B^+}$. Thus the full parameter vector is
\(
    \Theta
    =
    \left(
    \psi^A,\rho^A,\lambda,\rho^B
    \right),
\)
where
\(
    \sum_{i\in \NN_A^+}\psi_i^A=1,
    \sum_{\ell\in \NN_B^+}\lambda_{i,\ell}=1,
 \forall i\in \NN_A^+,
\)
and the rows of $\rho^A$ and $\rho^B$ are stochastic. We set the no-purchase
states $0_A$ and $0_B$ to be absorbing. Any additional structural restrictions,
such as zero self-transitions for real products, can be imposed by fixing the
corresponding entries and normalizing over the remaining feasible entries.

Let
\(
    \mathcal H
    =
    \{(S_A^t,S_B^t,a^t,b^t)\}_{t\in[T]}
\)
denote the observed purchase history, where $S_A^t\subseteq \NN_A$ and
$S_B^t\subseteq \NN_B$ are the offered assortments, and
$a^t\in S_A^{t,+}$ and $b^t\in S_B^{t,+}$ are the observed choices. The
observed-data log-likelihood is
\begin{equation}
    \log L(\Theta;\mathcal H)
    =
    \sum_{t=1}^T
    \log p_A(a^t,S_A^t;\psi^A,\rho^A)
    +
    \sum_{t=1}^T
    \log p_B(b^t,S_B^t;\lambda^{a^t},\rho^B),
    \label{eq:markovmc-observed-ll}
\end{equation}
where $p_A(a,S_A;\psi^A,\rho^A)$ is the probability that the category-$A$ MC
model chooses $a$ from $S_A^+$, and
$p_B(b,S_B;\lambda^a,\rho^B)$ is the probability that the category-$B$ MC model,
with initial distribution $\lambda^a$ and transition matrix $\rho^B$, chooses
$b$ from $S_B^+$.

\paragraph{A single-category MC calculation.}
We first state the quantities needed for the E-step for a generic category
$C\in\{A,B\}$. Let $\NN_C^+=\NN_C\cup\{0_C\}$, let $S_C^+=S_C\cup\{0_C\}$, and let
$k\in S_C^+$ be the observed final choice. Consider an MC model with initial
arrival distribution $\mu\in\Delta(\NN_C^+)$ and transition matrix $\rho^C$.

For each state $i\in \NN_C^+$, let
\(
    M_i^C(S_C;\mu,\rho^C)
\)
denote the expected number of visits to state $i$ up to and including
the first visit to $S_C^+$. These quantities solve the linear system
\begin{equation}
    M_i^C(S_C;\mu,\rho^C)
    =
    \mu_i
    +
    \sum_{u\in \NN_C^+\setminus S_C^+}
    M_u^C(S_C;\mu,\rho^C)\rho_{u,i}^C,
    \qquad i\in \NN_C^+.
    \label{eq:mc-visit-system}
\end{equation}
Since states in $S_C^+$ are absorbing, for $k\in S_C^+$ the choice probability
of $k$ is
\(
    p_C(k,S_C;\mu,\rho^C)
    =
    M_k^C(S_C;\mu,\rho^C).
\)

Next, define
\(
    \Psi_{k,i}^C(S_C;\rho^C)
\)
as the probability that the chain is eventually absorbed at $k$, conditional on
starting from state $i$. For $i\in S_C^+$, this probability is
\(
    \Psi_{k,i}^C(S_C;\rho^C)=\mathbf 1\{i=k\}.
\)
For $i\notin S_C^+$, it solves
\begin{equation}
    \Psi_{k,i}^C(S_C;\rho^C)
    =
    \rho_{i,k}^C
    +
    \sum_{j\in \NN_C^+\setminus S_C^+}
    \rho_{i,j}^C
    \Psi_{k,j}^C(S_C;\rho^C).
    \label{eq:mc-hitting-system}
\end{equation}
Equivalently, one can solve the linear system over the unavailable states
$\NN_C^+\setminus S_C^+$ and use the boundary values
$\Psi_{k,i}^C(S_C;\rho^C)=\mathbf 1\{i=k\}$ for $i\in S_C^+$.

We now define the missing data. Let $F_i^C$ be the indicator that the initial
arrival state in category $C$ is $i$, and let $X_{i,j}^C(S_C)$ be the number of
times the Markov chain transitions from state $i$ to state $j$ before absorption
when assortment $S_C$ is offered. Conditional on observing final choice
$k\in S_C^+$, the posterior probability of initial state $i$ is
\begin{equation}
    \mathbb E[F_i^C\mid k,S_C,\mu,\rho^C]
    =
    \frac{
    \mu_i\Psi_{k,i}^C(S_C;\rho^C)
    }{
    p_C(k,S_C;\mu,\rho^C)
    }.
    \label{eq:mc-posterior-initial}
\end{equation}
Similarly, the expected number of transitions from $i$ to $j$ is
\begin{equation}
    \mathbb E[X_{i,j}^C(S_C)\mid k,S_C,\mu,\rho^C]
    =
    \begin{cases}
    \displaystyle
    \frac{
    M_i^C(S_C;\mu,\rho^C)\rho_{i,j}^C
    \Psi_{k,j}^C(S_C;\rho^C)
    }{
    p_C(k,S_C;\mu,\rho^C)
    },
    & i\notin S_C^+, \\[2ex]
    0,
    & i\in S_C^+ .
    \end{cases}
    \label{eq:mc-posterior-transitions}
\end{equation}
The second case follows because once the chain reaches an available state, it
absorbs and no further transition is made.

\paragraph{Complete-data log-likelihood.}
For each transaction $t$, let $F_{A,i}^t$ denote the initial-state indicator in
category $A$, and let $X_{A,i,j}^t$ denote the number of category-$A$ transitions
from $i$ to $j$. Similarly, let $F_{B,\ell}^t$ denote the initial-state indicator
in category $B$, and let $X_{B,\ell,m}^t$ denote the number of category-$B$
transitions from $\ell$ to $m$.

If these latent variables were observed, the complete-data log-likelihood would
be
\begin{equation}
\begin{aligned}
    \log L_C(\Theta;\mathcal H,F,X)
    =
    &
    \sum_{t=1}^T
    \sum_{i\in \NN_A^+}
    F_{A,i}^t\log \psi_i^A
    +
    \sum_{t=1}^T
    \sum_{i,j\in \NN_A^+}
    X_{A,i,j}^t\log \rho_{i,j}^A
    \\
    &
    +
    \sum_{t=1}^T
    \sum_{\ell\in \NN_B^+}
    F_{B,\ell}^t\log \lambda_{a^t,\ell}
    +
    \sum_{t=1}^T
    \sum_{\ell,m\in \NN_B^+}
    X_{B,\ell,m}^t\log \rho_{\ell,m}^B .
\end{aligned}
\label{eq:markovmc-complete-ll}
\end{equation}
The first two terms correspond to the path in category $A$. The third term
corresponds to the choice-dependent initial arrival distribution into category
$B$, and the fourth term corresponds to the substitution path in category $B$.

\paragraph{E-step.}
At iteration $s$, suppose the current parameter estimate is
\(
    \Theta^{(s)}
    =
    \left(
    \psi^{A,(s)},\rho^{A,(s)},\lambda^{(s)},\rho^{B,(s)}
    \right).
\)
For each transaction $t$, we first apply
\eqref{eq:mc-visit-system}--\eqref{eq:mc-posterior-transitions} to category
$A$ with
\(
    C=A,
    S_C=S_A^t,
    k=a^t,
    \mu=\psi^{A,(s)},
    \rho^C=\rho^{A,(s)}.
\)
This gives
\begin{equation}
    \widehat F_{A,i}^{t,(s)}
    =
    \frac{
    \psi_i^{A,(s)}
    \Psi_{a^t,i}^{A,t,(s)}
    }{
    p_A^{t,(s)}
    },
    \qquad i\in \NN_A^+,
    \label{eq:estep-A-initial}
\end{equation}
where
\(
    p_A^{t,(s)}
    =
    p_A(a^t,S_A^t;\psi^{A,(s)},\rho^{A,(s)}),
\)
and we use $\Psi_{a^t,i}^{A,t,(s)}$ as shorthand for
$\Psi_{a^t,i}^A(S_A^t;\rho^{A,(s)})$. The expected transition counts are
\begin{equation}
    \widehat X_{A,i,j}^{t,(s)}
    =
    \begin{cases}
    \displaystyle
    \frac{
    M_i^{A,t,(s)}\rho_{i,j}^{A,(s)}
    \Psi_{a^t,j}^{A,t,(s)}
    }{
    p_A^{t,(s)}
    },
    & i\notin S_A^{t,+}, \\[2ex]
    0,
    & i\in S_A^{t,+},
    \end{cases}
    \label{eq:estep-A-transition}
\end{equation}
where
\[
    M_i^{A,t,(s)}
    =
    M_i^A(S_A^t;\psi^{A,(s)},\rho^{A,(s)}).
\]

Next, we apply the same calculation to category $B$. Conditional on the observed
category-$A$ choice $a^t$, the initial arrival distribution into category $B$ is
the row $\lambda_{a^t,\cdot}^{(s)}$. Thus we use
\(
    C=B,
    S_C=S_B^t,
    k=b^t,
    \mu=\lambda_{a^t,\cdot}^{(s)},
    \rho^C=\rho^{B,(s)}.
\)
This gives
\begin{equation}
    \widehat F_{B,\ell}^{t,(s)}
    =
    \frac{
    \lambda_{a^t,\ell}^{(s)}
    \Psi_{b^t,\ell}^{B,t,(s)}
    }{
    p_B^{t,(s)}
    },
    \qquad \ell\in \NN_B^+,
    \label{eq:estep-B-initial}
\end{equation}
where
\(
    p_B^{t,(s)}
    =
    p_B(b^t,S_B^t;\lambda_{a^t,\cdot}^{(s)},\rho^{B,(s)}),
\)
and
\(
    \Psi_{b^t,\ell}^{B,t,(s)}
    =
    \Psi_{b^t,\ell}^B(S_B^t;\rho^{B,(s)}).
\)
The expected category-$B$ transition counts are
\begin{equation}
    \widehat X_{B,\ell,m}^{t,(s)}
    =
    \begin{cases}
    \displaystyle
    \frac{
    M_\ell^{B,t,(s)}
    \rho_{\ell,m}^{B,(s)}
    \Psi_{b^t,m}^{B,t,(s)}
    }{
    p_B^{t,(s)}
    },
    & \ell\notin S_B^{t,+}, \\[2ex]
    0,
    & \ell\in S_B^{t,+},
    \end{cases}
    \label{eq:estep-B-transition}
\end{equation}
where
\(
    M_\ell^{B,t,(s)}
    =
    M_\ell^B(S_B^t;\lambda_{a^t,\cdot}^{(s)},\rho^{B,(s)}).
\)

\paragraph{M-step.}
Substituting the conditional expectations from the E-step into
\eqref{eq:markovmc-complete-ll}, the expected complete-data log-likelihood
decomposes as
\[
    Q(\Theta\mid \Theta^{(s)})
    =
    Q_{\psi^A}(\psi^A)
    +
    Q_{\rho^A}(\rho^A)
    +
    Q_{\lambda}(\lambda)
    +
    Q_{\rho^B}(\rho^B),
\]
where
\[
    Q_{\psi^A}(\psi^A)
    =
    \sum_{t=1}^T
    \sum_{i\in \NN_A^+}
    \widehat F_{A,i}^{t,(s)}\log \psi_i^A,\quad
    Q_{\rho^A}(\rho^A)
    =
    \sum_{t=1}^T
    \sum_{i,j\in \NN_A^+}
    \widehat X_{A,i,j}^{t,(s)}\log \rho_{i,j}^A,
\]
\[
    Q_{\lambda}(\lambda)
    =
    \sum_{t=1}^T
    \sum_{\ell\in \NN_B^+}
    \widehat F_{B,\ell}^{t,(s)}\log \lambda_{a^t,\ell}, \quad
    Q_{\rho^B}(\rho^B)
    =
    \sum_{t=1}^T
    \sum_{\ell,m\in \NN_B^+}
    \widehat X_{B,\ell,m}^{t,(s)}\log \rho_{\ell,m}^B.
\]
Without positive lower bounds, the M-step separates into normalized expected-count updates.

The category-$A$ initial arrival probabilities are updated by
\begin{equation}
    \psi_i^{A,(s+1)}
    =
    \frac{
    \sum_{t=1}^T \widehat F_{A,i}^{t,(s)}
    }{
    \sum_{t=1}^T\sum_{h\in \NN_A^+}\widehat F_{A,h}^{t,(s)}
    },
    \qquad i\in \NN_A^+.
    \label{eq:mstep-psiA}
\end{equation}
Since $\sum_{h\in \NN_A^+}\widehat F_{A,h}^{t,(s)}=1$ for every $t$, the
denominator is equal to $T$.

For every $i\in \NN_A^+$, the transition probabilities in category $A$ are updated
row by row:
\begin{equation}
    \rho_{i,j}^{A,(s+1)}
    =
    \frac{
    \sum_{t=1}^T \widehat X_{A,i,j}^{t,(s)}
    }{
    \sum_{t=1}^T\sum_{h\in \NN_A^+}\widehat X_{A,i,h}^{t,(s)}
    },
    \qquad j\in \NN_A^+.
    \label{eq:mstep-rhoA}
\end{equation}
If the denominator in \eqref{eq:mstep-rhoA} is zero, then row $i$ is not
identified by the current expected paths; in that case, we leave the row
unchanged or use a small smoothing constant.

For every observed category-$A$ outcome $i\in \NN_A^+$, the downstream arrival row
$\lambda_i$ is updated by
\begin{equation}
    \lambda_{i,\ell}^{(s+1)}
    =
    \frac{
    \sum_{t:a^t=i}\widehat F_{B,\ell}^{t,(s)}
    }{
    \sum_{t:a^t=i}\sum_{m\in \NN_B^+}\widehat F_{B,m}^{t,(s)}
    },
    \qquad \ell\in \NN_B^+.
    \label{eq:mstep-lambda}
\end{equation}
If no transaction has $a^t=i$, then the row $\lambda_i$ is not identified from
the data and can be left unchanged or smoothed.

Finally, the shared category-$B$ transition matrix is updated by aggregating
expected transition counts across all transactions:
\begin{equation}
    \rho_{\ell,m}^{B,(s+1)}
    =
    \frac{
    \sum_{t=1}^T \widehat X_{B,\ell,m}^{t,(s)}
    }{
    \sum_{t=1}^T\sum_{q\in \NN_B^+}\widehat X_{B,\ell,q}^{t,(s)}
    },
    \qquad m\in \NN_B^+.
    \label{eq:mstep-rhoB}
\end{equation}
Again, rows with zero expected transition counts are left unchanged or updated
with a small smoothing constant. If some entries of $\rho^A$ or $\rho^B$ are
fixed by modeling restrictions, then the sums in
\eqref{eq:mstep-rhoA} and \eqref{eq:mstep-rhoB} are taken only over the free
entries in the corresponding row.

\paragraph{Algorithm summary.}
The EM algorithm proceeds as follows.

\begin{algorithm}[H]
\caption{EM Algorithm for the Two-Category MarkovMC Model}
\label{alg:markovmc-em}
\KwIn{Observed data $\mathcal H=\{(S_A^t,S_B^t,a^t,b^t)\}_{t\in[T]}$; initial
parameters $\Theta^{(0)}=(\psi^{A,(0)},\rho^{A,(0)},\lambda^{(0)},\rho^{B,(0)})$.}
Set $s\leftarrow 0$\;
\Repeat{convergence}{
    \For{$t=1,\ldots,T$}{
        Compute $M^{A,t,(s)}$ from \eqref{eq:mc-visit-system} with
        $(S_C,\mu,\rho^C)=(S_A^t,\psi^{A,(s)},\rho^{A,(s)})$\;
        Compute $\Psi^{A,t,(s)}$ from \eqref{eq:mc-hitting-system} with
        $(S_C,k,\rho^C)=(S_A^t,a^t,\rho^{A,(s)})$\;
        Compute $\widehat F_A^{t,(s)}$ and $\widehat X_A^{t,(s)}$ from
        \eqref{eq:estep-A-initial} and \eqref{eq:estep-A-transition}\;

        Compute $M^{B,t,(s)}$ from \eqref{eq:mc-visit-system} with
        $(S_C,\mu,\rho^C)=(S_B^t,\lambda_{a^t,\cdot}^{(s)},\rho^{B,(s)})$\;
        Compute $\Psi^{B,t,(s)}$ from \eqref{eq:mc-hitting-system} with
        $(S_C,k,\rho^C)=(S_B^t,b^t,\rho^{B,(s)})$\;
        Compute $\widehat F_B^{t,(s)}$ and $\widehat X_B^{t,(s)}$ from
        \eqref{eq:estep-B-initial} and \eqref{eq:estep-B-transition}\;
    }
    Update $\psi^{A,(s+1)}$ by \eqref{eq:mstep-psiA}\;
    Update $\rho^{A,(s+1)}$ by \eqref{eq:mstep-rhoA}\;
    Update $\lambda^{(s+1)}$ by \eqref{eq:mstep-lambda}\;
    Update $\rho^{B,(s+1)}$ by \eqref{eq:mstep-rhoB}\;
    If positive lower bounds are imposed, replace these normalized-count updates with the exact solutions of the corresponding constrained row-wise M-step problems\;
    Set $s\leftarrow s+1$\;
}
\Return{$\Theta^{(s)}$}
\end{algorithm}

\paragraph{Computational complexity.}
For a category $C$ and transaction $t$, the E-step requires solving the visit
system \eqref{eq:mc-visit-system} and the hitting system
\eqref{eq:mc-hitting-system} over the unavailable states
$\NN_C^+\setminus S_C^{t,+}$. Therefore, a direct implementation requires solving
linear systems of dimension at most $|\NN_C^+|$. The naive per-transaction cost is
$O(|\NN_A^+|^3+|\NN_B^+|^3)$, although matrix factorizations can be reused for
repeated assortments within each EM iteration. Category-$B$ visit vectors also
depend on the observed upstream choice $a^t$, while hitting probabilities depend
on the observed final choice $b^t$. The M-step consists of row-wise multinomial
problems, with normalized-count solutions when no bound binds, and is polynomial
in $T$, $|\NN_A^+|$, and $|\NN_B^+|$.

\paragraph{Convergence.}
The EM algorithm above is a standard EM procedure with the unobserved Markov
paths as missing data. If each M-step exactly maximizes
$Q(\Theta\mid \Theta^{(s)})$ over the feasible parameter space, then the
observed-data log-likelihood in \eqref{eq:markovmc-observed-ll} is
non-decreasing across iterations. To obtain a compact interior parameter space,
one can impose a lower bound $\epsilon>0$ on all free probabilities:
\(
    \psi_i^A\ge \epsilon,
    \lambda_{i,\ell}\ge \epsilon,
    \rho_{i,j}^A\ge \epsilon,
  \rho_{\ell,m}^B\ge \epsilon
\)
for all free entries, together with the row-sum constraints, choosing $\epsilon$
so that the feasible set is nonempty. With these bounds, each M-step solves the
corresponding constrained multinomial maximum-likelihood problem; the
normalized-count formulas apply only when they satisfy the bounds. These
row-wise problems remain polynomial-time solvable. Under the standard EM
regularity conditions, the likelihood sequence converges to a finite limit, and
every accumulation point of the parameter sequence is a constrained stationary
point.

\begin{theorem}
\label{thm:markovmc-em}
Let
$\{\Theta^{(s)}=(\psi^{A,(s)},\rho^{A,(s)},\lambda^{(s)},\rho^{B,(s)})\}_{s\ge 0}$
be the sequence generated by Algorithm~\ref{alg:markovmc-em}. If each M-step is
solved exactly over the feasible parameter space, then
\(
    \log L(\Theta^{(s+1)};\mathcal H)
    \ge
    \log L(\Theta^{(s)};\mathcal H)
\) for all $s\ge 0$.
Moreover, if the parameter space is restricted by a positive lower bound
$\epsilon>0$ on all free probabilities and the standard EM regularity conditions
hold, then the likelihood sequence converges to a finite limit, and every
accumulation point of the parameter sequence is a constrained stationary point.
\end{theorem}

\proof{Proof.}
The complete-data log-likelihood in
\eqref{eq:markovmc-complete-ll} is obtained by treating the initial arrival
states and transition counts in both categories as missing data. In the E-step,
\eqref{eq:estep-A-initial}--\eqref{eq:estep-B-transition} compute the
conditional expectations of these missing-data sufficient statistics under the
current parameter vector. In the M-step, maximizing the expected complete-data
log-likelihood separates into row-wise multinomial maximum-likelihood problems.
Without positive lower bounds, these problems yield the normalized expected-count
updates \eqref{eq:mstep-psiA}--\eqref{eq:mstep-rhoB}. When lower bounds are
imposed, we solve the constrained row-wise problems exactly. The monotonicity of
the observed-data log-likelihood then follows from the standard EM inequality.
Under the stated compactness and regularity conditions, the likelihood sequence
converges to a finite limit, and every accumulation point of the parameter
sequence is a constrained stationary point.
\Halmos
\endproof

\subsection{General DAG Estimation}\label{app:est-dag}
The EM framework described for the two-category model can be generalized to estimate the parameters of a model defined on any DAG $\mathcal{G} = (\mathcal{V}, \mathcal{E})$. To show the core estimation mechanism, we focus our description on DAGs where each node has at most one inbound edge. The framework extends to arbitrary DAGs by treating the anchor parent of each multi-parent category as an additional latent variable (see Appendix~\ref{app:multi-parent}); the pooling weights can be kept fixed or updated in closed form from posterior responsibilities. The remaining details are omitted for conciseness.

The full parameter set is $\Theta = ( \{\mathbf{v}^U\}_{U \in \mathcal{V}}, \{\boldsymbol{\lambda}^{U,U'}\}_{(U,U') \in \mathcal{E}})$. Here, $\mathbf{v}^U$ represents the parameters for the internal choice model $\phi_U$ of category $U$ (e.g., an MNL model), and $\boldsymbol{\lambda}^{U,U'}$ is the matrix of transition probabilities for the edge $(U,U')$. {The observed data are $\mathcal H=\{(S_U^t,y_U^t)\}_{U\in\mathcal V,\,t\in[T]}$, where $S_U^t\subseteq\NN_U$ is the offered assortment and $y_U^t\in S_U^t\cup\{0_U\}$ is the observed choice in category $U$ for transaction $t$.}

The core of the EM algorithm is to introduce latent variables that describe the unobserved parts of a customer's journey. For each transaction $t$ that traverses an edge $(U, U')$, and for each product $i \in \NN_U^+$ chosen in category $U$, we define a latent variable $X_{m}^{t,(U,U')}$ to be $1$ if the customer's subsequent attention is directed to product/option $m \in \NN_{U'}^+$, and $0$ otherwise. Let $\mathbf{X}$ be the collection of all such latent variables for all relevant transactions and edges.

The key property of this formulation is that the expected complete-data log-likelihood, $Q(\Theta|\Theta^{(\ell)}) = \mathbb{E}_{\mathbf{X}|\mathcal{H},\Theta^{(\ell)}}[\log L_C(\Theta; \mathcal{H}, \mathbf{X})]$, decomposes into a sum of independent terms:
\[
Q(\Theta|\Theta^{(\ell)}) = \sum_{U \in \mathcal{V}} Q_U(\mathbf{v}^U|\Theta^{(\ell)}) + \sum_{(U,U') \in \mathcal{E}} Q_{U,U'}(\boldsymbol{\lambda}^{U,U'}|\Theta^{(\ell)}).
\]
Each component function $Q_U$ depends only on the internal parameters of a single category $U$, and each $Q_{U,U'}$ depends only on the transition probabilities for a single edge $(U,U')$.

\noindent\textbf{E-step:} At iteration $\ell$, the E-step involves computing the conditional expectation of each relevant latent variable: $\hat{X}_{m}^{t,(U,U'),(\ell)} = \mathbb{E}[X_{m}^{t,(U,U')} | \mathcal{H}, \Theta^{(\ell)}]$. This conditional expectation is calculated using the current parameter estimates $\Theta^{(\ell)}$ and the observed choices in the transaction path, analogous to the two-category case.

\noindent\textbf{M-step:} The M-step uses the decomposition of $Q$ to update all parameters by solving a set of smaller, independent maximization problems.
\begin{itemize}
    \item For each category $U \in \mathcal{V}$, the internal model parameters are updated by solving:
    \[
    \mathbf{v}^{U,(\ell+1)} = \arg\max_{\mathbf{v}^U} Q_U(\mathbf{v}^U|\Theta^{(\ell)}).
    \]
    This subproblem is equivalent to a weighted maximum likelihood estimation problem for the choice model $\phi_U$. For a non-root category, each conditional-choice likelihood contribution is weighted by the posterior probability of its latent initial attraction, as computed in the E-step. Root categories contribute their ordinary observed-choice likelihoods.

    \item For each edge $(U, U') \in \mathcal{E}$, the update for the transition probabilities has a closed-form solution. The objective $Q_{U,U'}(\boldsymbol{\lambda}^{U,U'}|\Theta^{(\ell)})$ is a sum of weighted logarithms:
    \[ \sum_{t=1}^T \sum_{m\in \NN_{U'}^+} \hat{X}_{m}^{t,(U,U'),(\ell)} \log \lambda_{y_U^t,m}^{U,U'}. \]
    This function is jointly concave in $\bfl^{U,U'}$, and thus can be maximized efficiently subject to the constraint $\sum_{m \in \NN_{U'}^+} \lambda_{i,m}^{U,U'} = 1$ for each $i \in \NN_U^+$.
\end{itemize}
\section{Bidirectional Extension: Estimation and Optimization}
\label{app:bi-opt}

This appendix provides the estimation and optimization details for the bidirectional extension introduced in Section~\ref{sec:bidirectional}. The model is a mixture of two directional models. With probability \(\gamma\), a customer follows the direction \(A\to B\); with probability \(1-\gamma\), the customer follows the direction \(B\to A\).
\subsection{EM Algorithm for the Bidirectional MNL Special Case}

We present the EM algorithm for the MNL special case, which is the version used for direct comparison with the MNL-based complementarity literature. Let
\(
\Theta=(\gamma,\Lambda^{AB},\Lambda^{BA},\bfv^A,\bfv^B),
\)
where \(\Lambda^{AB}=(\lambda^{AB}_{i,m})_{i\in\NN_A^+,m\in\NN_B^+}\) and
\(\Lambda^{BA}=(\lambda^{BA}_{j,k})_{j\in\NN_B^+,k\in\NN_A^+}\). The within-category MNL weights \(\bfv^A\) and \(\bfv^B\) are shared across the two directions.

For transaction \(t\), define
\(
\eta_t^A
=
\frac{v_{a^t}^A}{1+V^A(S_A^t)},
\eta_t^B
=
\frac{v_{b^t}^B}{1+V^B(S_B^t)},
\)
with the convention that \(v_{0_A}^A=v_{0_B}^B=1\). Conditional on the direction \(A\to B\), the likelihood contribution is
\(
L_t^{AB}
=
\eta_t^A
\left(
\lambda^{AB}_{a^t,b^t}
+
\eta_t^B
\sum_{m\in \NN_B\setminus S_B^t}
\lambda^{AB}_{a^t,m}
\right).
\)
Conditional on the direction \(B\to A\), the likelihood contribution is
\(
L_t^{BA}
=
\eta_t^B
\Bigl(
\lambda^{BA}_{b^t,a^t}
+
\eta_t^A
\sum_{k\in \NN_A\setminus S_A^t}
\lambda^{BA}_{b^t,k}
\Bigr).
\)
Thus, the observed-data log-likelihood is
\begin{equation}
\label{eq:bi-obs-ll}
\log L(\Theta;\mathcal H)
=
\sum_{t=1}^T
\log\left(
\gamma L_t^{AB}
+
(1-\gamma)L_t^{BA}
\right).
\end{equation}

\paragraph{Latent variables.}
We introduce a latent direction variable \(Y^t\in\{0,1\}\), where \(Y^t=1\) indicates that transaction \(t\) follows direction \(A\to B\), and \(Y^t=0\) indicates direction \(B\to A\). We also introduce direction-specific latent initial-attraction indicators. For direction \(A\to B\), let \(X_{t,m}^{AB}=1\) if the customer is initially attracted to \(m\in\NN_B^+\) after choosing \(a^t\). For direction \(B\to A\), let \(X_{t,k}^{BA}=1\) if the customer is initially attracted to \(k\in\NN_A^+\) after choosing \(b^t\).

\paragraph{E-step.}
At iteration \(\ell\), define
\(
\Delta_t^{AB,(\ell)}
=
\lambda^{AB,(\ell)}_{a^t,b^t}
+
\eta_t^{B,(\ell)}
\sum_{m\in \NN_B\setminus S_B^t}
\lambda^{AB,(\ell)}_{a^t,m},
\)
and
\(
\Delta_t^{BA,(\ell)}
=
\lambda^{BA,(\ell)}_{b^t,a^t}
+
\eta_t^{A,(\ell)}
\sum_{k\in \NN_A\setminus S_A^t}
\lambda^{BA,(\ell)}_{b^t,k}.
\)
The posterior probability that transaction \(t\) follows direction \(A\to B\) is
\begin{equation}
\label{eq:bi-e-y}
\hat y_t^{(\ell)}
=
\frac{
\gamma^{(\ell)}L_t^{AB}(\Theta^{(\ell)})
}{
\gamma^{(\ell)}L_t^{AB}(\Theta^{(\ell)})
+
(1-\gamma^{(\ell)})L_t^{BA}(\Theta^{(\ell)})
}.
\end{equation}
For \(m\in\NN_B^+\), the posterior expected initial-attraction indicator under direction \(A\to B\) is
\[
\hat x_{t,m}^{AB,(\ell)}
=
\begin{cases}
\displaystyle
\hat y_t^{(\ell)}
\frac{\lambda^{AB,(\ell)}_{a^t,b^t}}
{\Delta_t^{AB,(\ell)}},
& m=b^t,\\[8pt]
\displaystyle
\hat y_t^{(\ell)}
\frac{\lambda^{AB,(\ell)}_{a^t,m}\eta_t^{B,(\ell)}}
{\Delta_t^{AB,(\ell)}},
& m\in\NN_B\setminus S_B^t,\\[8pt]
0, & \text{otherwise.}
\end{cases}
\]
Similarly, for \(k\in\NN_A^+\), the posterior expected initial-attraction indicator under direction \(B\to A\) is
\[
\hat x_{t,k}^{BA,(\ell)}
=
\begin{cases}
\displaystyle
(1-\hat y_t^{(\ell)})
\frac{\lambda^{BA,(\ell)}_{b^t,a^t}}
{\Delta_t^{BA,(\ell)}},
& k=a^t,\\[8pt]
\displaystyle
(1-\hat y_t^{(\ell)})
\frac{\lambda^{BA,(\ell)}_{b^t,k}\eta_t^{A,(\ell)}}
{\Delta_t^{BA,(\ell)}},
& k\in\NN_A\setminus S_A^t,\\[8pt]
0, & \text{otherwise.}
\end{cases}
\]

\paragraph{M-step.}
The mixture weight has the closed-form update
\(
\gamma^{(\ell+1)}
=
\frac{1}{T}\sum_{t=1}^T \hat y_t^{(\ell)}.
\)
The transition probabilities are updated by normalized posterior counts:
\[
\lambda_{i,m}^{AB,(\ell+1)}
=
\frac{
\sum_{t:a^t=i}\hat x_{t,m}^{AB,(\ell)}
}{
\sum_{t:a^t=i}\hat y_t^{(\ell)}
},
\qquad
i\in\NN_A^+,\ m\in\NN_B^+,
\]
and
\(
\lambda_{j,k}^{BA,(\ell+1)}
=
\frac{
\sum_{t:b^t=j}\hat x_{t,k}^{BA,(\ell)}
}{
\sum_{t:b^t=j}(1-\hat y_t^{(\ell)})
},
j\in\NN_B^+,\ k\in\NN_A^+.
\)
Rows with zero posterior mass can be left unchanged or updated with a small smoothing constant; this does not affect the likelihood on the observed data.

It remains to update the MNL weights. Let \(v_i^A=\exp(\alpha_i)\) and \(v_j^B=\exp(\beta_j)\), with \(\alpha_{0_A}=\beta_{0_B}=0\). The relevant transaction weights are
\(
w_t^{A,(\ell)}
=
\hat y_t^{(\ell)}
+
\sum_{k\in\NN_A\setminus S_A^t}
\hat x_{t,k}^{BA,(\ell)},
\)
and
\(
w_t^{B,(\ell)}
=
1-\hat y_t^{(\ell)}
+
\sum_{m\in\NN_B\setminus S_B^t}
\hat x_{t,m}^{AB,(\ell)}.
\)
Then \(\boldsymbol{\alpha}\) is updated by solving the weighted MNL problem
\(
\max_{\boldsymbol{\alpha}}
\sum_{t=1}^T
w_t^{A,(\ell)}
\left[
\alpha_{a^t}
-
\log\left(1+\sum_{i\in S_A^t}\exp(\alpha_i)\right)
\right],
\)
and \(\boldsymbol{\beta}\) is updated by solving
\(
\max_{\boldsymbol{\beta}}
\sum_{t=1}^T
w_t^{B,(\ell)}
\left[
\beta_{b^t}
-
\log\left(1+\sum_{j\in S_B^t}\exp(\beta_j)\right)
\right].
\)
Both problems are concave in the log-weights. If each M-step attains its maximum and is solved exactly, the observed-data log-likelihood in~\eqref{eq:bi-obs-ll} is non-decreasing and converges to a finite limit. Under the standard EM regularity conditions, every accumulation point of the parameter sequence is a constrained stationary point of the observed-data log-likelihood.

\subsection{A \texorpdfstring{\(\frac{1}{2}\)}{1/2}-Approximation for Bidirectional Assortment Optimization}
\label{app:bi-opt-approx}

We now consider assortment optimization under the bidirectional mixture model. For a pair of assortments \((S_A,S_B)\), let \(R^{AB}(S_A,S_B)\) denote the expected revenue under the \(A\to B\) directional model, and let \(R^{BA}(S_A,S_B)\) denote the expected revenue under the \(B\to A\) directional model. The bidirectional objective is
\begin{equation}
\label{eq:bi-mix-revenue}
R^{\mathrm{bi}}(S_A,S_B)
=
\gamma R^{AB}(S_A,S_B)
+
(1-\gamma)R^{BA}(S_A,S_B).
\end{equation}

When the within-category models are MC models, each directional problem can be solved in polynomial time by Theorem~\ref{thm:opt}. Specifically, the \(A\to B\) problem is solved exactly by Algorithm~\ref{alg:mc_assortment}, and the \(B\to A\) problem is solved by applying the same algorithm after swapping the roles of the two categories.

Let
\(
(S_A^{AB},S_B^{AB})
\in
\arg\max_{S_A,S_B} R^{AB}(S_A,S_B),
\)
and
\(
(S_A^{BA},S_B^{BA})
\in
\arg\max_{S_A,S_B} R^{BA}(S_A,S_B).
\)
The following simple procedure computes both directional optima and returns the one with the larger bidirectional objective:
\begin{equation}
\label{eq:bi-alg}
(S_A^{\ALG},S_B^{\ALG})
\in
\underset{(S_A,S_B)\in
\{(S_A^{AB},S_B^{AB}),(S_A^{BA},S_B^{BA})\}}
{\arg\max}
R^{\mathrm{bi}}(S_A,S_B).
\end{equation}

\begin{theorem}
\label{thm:bi-half-approx}
If both within-category choice models are MC models, the procedure in~\eqref{eq:bi-alg} is a polynomial-time \(\frac{1}{2}\)-approximation for the bidirectional assortment optimization problem~\eqref{eq:bi-mix-revenue}.
\end{theorem}

\proof{Proof.}
Let \((S_A^*,S_B^*)\) be an optimal solution to~\eqref{eq:bi-mix-revenue}. By optimality of \((S_A^{AB},S_B^{AB})\) for the \(A\to B\) directional problem,
\(
R^{AB}(S_A^{AB},S_B^{AB})
\ge
R^{AB}(S_A^*,S_B^*).
\)
Similarly,
\(
R^{BA}(S_A^{BA},S_B^{BA})
\ge
R^{BA}(S_A^*,S_B^*).
\)
Define
\(
X=\gamma R^{AB}(S_A^*,S_B^*),
Y=(1-\gamma)R^{BA}(S_A^*,S_B^*).
\)
Then the optimal bidirectional revenue is \(X+Y\). Evaluating the bidirectional objective at the \(A\to B\) directional optimum gives
\(
R^{\mathrm{bi}}(S_A^{AB},S_B^{AB})
\ge
\gamma R^{AB}(S_A^{AB},S_B^{AB})
\ge
X.
\)
Similarly,
\[
R^{\mathrm{bi}}(S_A^{BA},S_B^{BA})
\ge
(1-\gamma)R^{BA}(S_A^{BA},S_B^{BA})
\ge
Y.
\]
Since the algorithm returns the better of these two solutions,
\[
R^{\mathrm{bi}}(S_A^{\ALG},S_B^{\ALG})
\ge
\max\{X,Y\}
\ge
\frac{1}{2}(X+Y)
=
\frac{1}{2}R^{\mathrm{bi}}(S_A^*,S_B^*).
\]
The two directional problems are solved in polynomial time under MC choice models by Theorem~\ref{thm:opt}, so the overall procedure is polynomial-time.
\Halmos
\endproof
\section{Additional Experimental Results}\label{app:syn}
In this section, we present all omitted tables and additional experimental results referenced in Sections~\ref{sec:synthetic_experiments} and~\ref{sec:real_data}.
\subsection{Tables for Synthetic Data Results}
 The results are presented in Tables \ref{tab:train_ll_synthetic}, \ref{tab:test_ll_synthetic}, \ref{tab:hit_rate_synthetic} and \ref{tab:mean_rank_synthetic}.
\begingroup
\begin{table*}[htbp]
    \centering
    \smallXnoCE
    \setlength{\tabcolsep}{4pt}
    \renewcommand{\arraystretch}{1.15}
    \caption{In-sample (train) log-likelihood on synthetic data for varying $\theta$. Values in parentheses indicate the percentage improvement relative to \textsc{IndMNL}.}
    \label{tab:train_ll_synthetic}
    \begin{tabular}{@{}lccccc@{}}
        \toprule
        $\theta$ & \textsc{IndMNL} & \textsc{MultiMNL} & \textsc{MarkovMNL} & \textsc{IndMC} & \textsc{MarkovMC} \\
        \midrule
        0.0 & -30123.08 & -30048.38 ($+0.25\,\%$) & -29718.71 ($+1.34\,\%$) & -27774.41 ($+7.80\,\%$) & -27668.55 ($+8.15\,\%$) \\
        1.0 & -30286.99 & -30037.19 ($+0.82\,\%$) & -29730.07 ($+1.84\,\%$) & -28030.84 ($+7.45\,\%$) & -27690.18 ($+8.57\,\%$) \\
        2.0 & -30450.87 & -29711.46 ($+2.43\,\%$) & -29460.01 ($+3.25\,\%$) & -28361.71 ($+6.86\,\%$) & -27474.58 ($+9.77\,\%$) \\
        3.0 & -30663.80 & -29342.52 ($+4.31\,\%$) & -29168.83 ($+4.88\,\%$) & -28636.90 ($+6.61\,\%$) & -27168.38 ($+11.40\,\%$) \\
        4.0 & -30777.86 & -28837.56 ($+6.30\,\%$) & -28772.65 ($+6.52\,\%$) & -28824.65 ($+6.35\,\%$) & -26759.94 ($+13.05\,\%$) \\
        5.0 & -30816.80 & -28211.80 ($+8.45\,\%$) & -28294.12 ($+8.19\,\%$) & -28871.14 ($+6.31\,\%$) & -26241.14 ($+14.85\,\%$) \\
        6.0 & -30909.25 & -27862.61 ($+9.86\,\%$) & -28081.26 ($+9.15\,\%$) & -28985.83 ($+6.22\,\%$) & -25988.62 ($+15.92\,\%$) \\
        7.0 & -30905.56 & -27534.10 ($+10.91\,\%$) & -27888.99 ($+9.76\,\%$) & -29000.63 ($+6.16\,\%$) & -25749.11 ($+16.68\,\%$) \\
        8.0 & -30934.95 & -27285.07 ($+11.80\,\%$) & -27754.78 ($+10.28\,\%$) & -29041.88 ($+6.12\,\%$) & -25592.92 ($+17.27\,\%$) \\
        9.0 & -30989.88 & -27069.99 ($+12.65\,\%$) & -27656.16 ($+10.76\,\%$) & -29067.26 ($+6.20\,\%$) & -25432.93 ($+17.93\,\%$) \\
        \bottomrule
    \end{tabular}
\end{table*}

\begin{table*}[htbp]
    \centering
    \smallXnoCE
    \setlength{\tabcolsep}{4pt}
    \renewcommand{\arraystretch}{1.15}
       \caption{Out-of-sample (test) log-likelihood on synthetic data for varying $\theta$. Values in parentheses indicate the percentage improvement relative to \textsc{IndMNL}.}
    \label{tab:test_ll_synthetic}
    \begin{tabular}{@{}lccccc@{}}
        \toprule
        $\theta$ & \textsc{IndMNL} & \textsc{MultiMNL} & \textsc{MarkovMNL} & \textsc{IndMC} & \textsc{MarkovMC} \\
        \midrule
        0.0 & -12889.44 & -12949.52 ($-0.47\,\%$) & -12770.66 ($+0.92\,\%$) & -11910.44 ($+7.60\,\%$) & -11914.80 ($+7.56\,\%$) \\
        1.0 & -12987.93 & -12972.57 ($+0.12\,\%$) & -12800.53 ($+1.44\,\%$) & -12050.07 ($+7.22\,\%$) & -11947.64 ($+8.01\,\%$) \\
        2.0 & -13073.84 & -12856.65 ($+1.66\,\%$) & -12693.01 ($+2.91\,\%$) & -12199.11 ($+6.69\,\%$) & -11859.38 ($+9.29\,\%$) \\
        3.0 & -13122.42 & -12640.19 ($+3.67\,\%$) & -12519.91 ($+4.59\,\%$) & -12293.28 ($+6.32\,\%$) & -11702.28 ($+10.82\,\%$) \\
        4.0 & -13164.21 & -12430.93 ($+5.57\,\%$) & -12343.88 ($+6.23\,\%$) & -12362.10 ($+6.09\,\%$) & -11517.15 ($+12.51\,\%$) \\
        5.0 & -13235.91 & -12235.67 ($+7.56\,\%$) & -12204.51 ($+7.79\,\%$) & -12433.19 ($+6.06\,\%$) & -11351.42 ($+14.24\,\%$) \\
        6.0 & -13272.26 & -12069.24 ($+9.06\,\%$) & -12087.72 ($+8.92\,\%$) & -12475.76 ($+6.00\,\%$) & -11220.18 ($+15.46\,\%$) \\
        7.0 & -13271.28 & -11948.93 ($+9.96\,\%$) & -12016.00 ($+9.46\,\%$) & -12468.95 ($+6.05\,\%$) & -11126.53 ($+16.16\,\%$) \\
        8.0 & -13250.36 & -11800.46 ($+10.94\,\%$) & -11928.98 ($+9.97\,\%$) & -12472.38 ($+5.87\,\%$) & -11038.60 ($+16.69\,\%$) \\
        9.0 & -13319.05 & -11741.11 ($+11.85\,\%$) & -11923.56 ($+10.48\,\%$) & -12521.37 ($+5.99\,\%$) & -11009.91 ($+17.34\,\%$) \\
        \bottomrule
    \end{tabular}
\end{table*}

\begin{table*}[htbp]
    \centering
    \smallXnoCE
    \setlength{\tabcolsep}{4pt}
    \renewcommand{\arraystretch}{1.15}
    \caption{Top-3 hit rate (test set) on synthetic data. Values in parentheses represent the difference in percentage points (p.p.) relative to \textsc{IndMNL}.}
    \label{tab:hit_rate_synthetic}
    \begin{tabular}{@{}lccccc@{}}
        \toprule
        $\theta$ & \textsc{IndMNL} & \textsc{MultiMNL} & \textsc{MarkovMNL} & \textsc{IndMC} & \textsc{MarkovMC} \\
        \midrule
        0.0 & 0.7770 & 0.7740 (-0.30 p.p.) & 0.7867 (+0.97 p.p.) & 0.8194 (+4.24 p.p.) & 0.8178 (+4.07 p.p.) \\
        1.0 & 0.7680 & 0.7748 (+0.68 p.p.) & 0.7864 (+1.84 p.p.) & 0.8077 (+3.97 p.p.) & 0.8153 (+4.73 p.p.) \\
        2.0 & 0.7624 & 0.7869 (+2.44 p.p.) & 0.7989 (+3.65 p.p.) & 0.7969 (+3.44 p.p.) & 0.8244 (+6.19 p.p.) \\
        3.0 & 0.7518 & 0.7962 (+4.44 p.p.) & 0.8037 (+5.18 p.p.) & 0.7861 (+3.43 p.p.) & 0.8327 (+8.09 p.p.) \\
        4.0 & 0.7459 & 0.8076 (+6.16 p.p.) & 0.8141 (+6.82 p.p.) & 0.7787 (+3.27 p.p.) & 0.8409 (+9.49 p.p.) \\
        5.0 & 0.7426 & 0.8157 (+7.31 p.p.) & 0.8187 (+7.61 p.p.) & 0.7731 (+3.04 p.p.) & 0.8463 (+10.37 p.p.) \\
        6.0 & 0.7371 & 0.8212 (+8.41 p.p.) & 0.8221 (+8.50 p.p.) & 0.7686 (+3.14 p.p.) & 0.8493 (+11.22 p.p.) \\
        7.0 & 0.7394 & 0.8244 (+8.49 p.p.) & 0.8258 (+8.64 p.p.) & 0.7693 (+2.98 p.p.) & 0.8510 (+11.15 p.p.) \\
        8.0 & 0.7396 & 0.8299 (+9.02 p.p.) & 0.8284 (+8.88 p.p.) & 0.7705 (+3.09 p.p.) & 0.8540 (+11.44 p.p.) \\
        9.0 & 0.7350 & 0.8319 (+9.69 p.p.) & 0.8316 (+9.66 p.p.) & 0.7678 (+3.27 p.p.) & 0.8555 (+12.05 p.p.) \\
        \bottomrule
    \end{tabular}
\end{table*}

\begin{table*}[htbp]
    \centering
    \smallXnoCE
    \setlength{\tabcolsep}{4pt}
    \renewcommand{\arraystretch}{1.15}
    \caption{Mean rank of the chosen item (test set) on synthetic data; lower is better. Values in parentheses give the percentage change relative to \textsc{IndMNL} (negative $=$ improvement).}
    \label{tab:mean_rank_synthetic}
    \begin{tabular}{@{}lccccc@{}}
        \toprule
        $\theta$ & \textsc{IndMNL} & \textsc{MultiMNL} & \textsc{MarkovMNL} & \textsc{IndMC} & \textsc{MarkovMC} \\
        \midrule
        0.0 & 2.414 & 2.429 (+0.62\%)  & 2.357 (-2.36\%)  & 2.222 (-7.95\%) & 2.229 (-7.66\%) \\
        1.0 & 2.453 & 2.421 (-1.30\%)  & 2.361 (-3.75\%)  & 2.278 (-7.13\%) & 2.243 (-8.56\%) \\
        2.0 & 2.486 & 2.363 (-4.95\%)  & 2.309 (-7.12\%)  & 2.331 (-6.23\%) & 2.201 (-11.46\%) \\
        3.0 & 2.524 & 2.307 (-8.60\%)  & 2.271 (-10.02\%) & 2.378 (-5.78\%) & 2.159 (-14.46\%) \\
        4.0 & 2.550 & 2.249 (-11.80\%) & 2.223 (-12.82\%) & 2.409 (-5.53\%) & 2.112 (-17.18\%) \\
        5.0 & 2.573 & 2.202 (-14.42\%) & 2.193 (-14.77\%) & 2.444 (-5.01\%) & 2.078 (-19.24\%) \\
        6.0 & 2.597 & 2.169 (-16.48\%) & 2.170 (-16.44\%) & 2.463 (-5.16\%) & 2.052 (-20.99\%) \\
        7.0 & 2.591 & 2.148 (-17.10\%) & 2.149 (-17.06\%) & 2.458 (-5.13\%) & 2.036 (-21.42\%) \\
        8.0 & 2.585 & 2.113 (-18.26\%) & 2.128 (-17.68\%) & 2.458 (-4.91\%) & 2.019 (-21.90\%) \\
        9.0 & 2.606 & 2.101 (-19.38\%) & 2.119 (-18.69\%) & 2.478 (-4.91\%) & 2.008 (-22.95\%) \\
        \bottomrule
    \end{tabular}
\end{table*}

\subsection{Real-Data Case Study: Extension to Three Product Categories}\label{app:three_group_extension}
To demonstrate our framework's scalability, we extend the `Meats and Buns' case study by introducing a third group of relevant condiments (category $C$), including cheese, pickles, and mayonnaise. We model this three-group dependency using a tree structure ($B\leftarrow A\rightarrow C$), which reflects the consumer insight that the primary choice of meat ($A$) is the main driver for the selection of both the bun ($B$) and the condiment ($C$). We acknowledge that other dependency structures, such as a sequence $A \to B \to C$, may exist; however, we choose the tree model because the purchase of condiments is typically more dependent on the meat than on the bun. {For comparison, we adapt the three benchmarks. For \textsc{IndMNL} and \textsc{IndMC}, we estimate three independent MNL and MC models, respectively, one for each category. For \textsc{MultiMNL}, we extend the two-category model to fit our tree structure. We introduce a set of preference weights $\{w_{i,k}\}_{i\in \NN_A^+,k\in \NN_C}$ for category $C$. We set $w_{0_A}=w_{i,0_B}=w_{i,0_C}=1$ for all $i\in\NN_A^+$. The probability of choosing $i\in S_A^+$, $j\in S_B^+$ and $k\in S_C^+$ is:
\(
P_{\textsc{MultiMNL}}(i,j,k, S_A, S_B,S_C) = \frac{w_i}{\sum_{m\in S_A}w_m+1}\frac{w_{i,j}}{\sum_{\ell \in S_B} w_{i,\ell}+1}\frac{w_{i,k}}{\sum_{q \in S_C} w_{i,q}+1}.
\)}
It is important to note that while this extension provides a valuable benchmark, the associated assortment optimization problem is NP-hard even for the original two-category case. {For each transaction, we sum the conditional log-likelihoods and the ranks of the realized choices in $B$ and $C$, conditional on the realized choice in $A$. We report the total conditional log-likelihood and the average sum of ranks across transactions. For the top-$3$ hit rate, we determine separately whether the realized choice in each of $B$ and $C$ is among the top three alternatives under its conditional distribution given the realized choice in $A$, and report the average of the two hit rates.}

The performance of our proposed \textsc{MarkovMC} and \textsc{MarkovMNL} models and the benchmarks is summarized in Table~\ref{tab:model_comparison_three_groups}. The results clearly demonstrate the advantage of modeling cross-category dependencies. Our \textsc{MarkovMC} model provides a substantial improvement over \textsc{IndMNL}, increasing the out-of-sample log-likelihood by $2.35\%$ and the top-$3$ hit rate by $1.3$ percentage points. Both \textsc{MultiMNL} and \textsc{MarkovMNL} also outperform the independent benchmarks across all reported metrics.

\begingroup
\begin{table}[htbp]
\centering
\caption{Model fit and prediction accuracy. For conditional LL, parenthetical values are relative improvements over \textsc{IndMNL}. For the mean sum of ranks, they are relative percentage changes from \textsc{IndMNL}. For HR, they are absolute differences in percentage points (p.p.).}
\label{tab:model_comparison_three_groups}
\resizebox{0.95\textwidth}{!}{%
\begin{tabular}{@{}lllll@{}}
\toprule
Model & In-Sample Conditional LL & Out-of-Sample Conditional LL & Top-3 HR & Mean Sum of Ranks \\
\midrule
\textsc{IndMNL}    & -116465.50           & -49956.41           & 65.37\%              & 15.28 \\
\textsc{IndMC}     & -115720.28 (+0.64\%) & -49670.14 (+0.57\%) & 65.84\% (+0.47 p.p.) & 15.17 (-0.72\%) \\
\textsc{MultiMNL}    & -113342.02 (+2.68\%) & -49035.43 (+1.84\%) & 66.15\% (+0.78 p.p.) & 14.96 (-2.09\%) \\
\textsc{MarkovMNL} & -113172.19 (+2.83\%) & -49183.04 (+1.55\%) & 66.19\% (+0.82 p.p.) & 14.95 (-2.16\%) \\
\textsc{MarkovMC}  & -112522.56 (+3.39\%) & -48783.30 (+2.35\%) & 66.67\% (+1.30 p.p.) & 14.85 (-2.81\%) \\
\bottomrule
\end{tabular}%
}
\end{table}
\endgroup

\section{Implementation Details for Cardinality-Constrained Optimization}
\label{sec:heuristic}
This appendix gives the implementation details for Algorithm~\ref{alg:alt_opt} and the \textsc{AltOpt} experiment summarized in Section~\ref{sec:cardinality_experiments}. The benchmark formulations, solver settings, and further implementation details are provided in Appendix~\ref{app:cardinality_computation}.

For the numerical study, we generate $50$ random two-category \textsc{MarkovMC} instances for each $N \in \{10,20,40,80,160,320,640\}$ with $n_A=n_B=N$. We benchmark \textsc{AltOpt} against MILP solutions certified optimal to a relative tolerance of $10^{-4}$ whenever the solver terminates within two hours ($N\le 160$). For larger instances, we report ratios against the LP relaxation upper bound. The table in Section~\ref{sec:cardinality_experiments} reports the resulting averages and worst-case ratios.

\endgroup

\subsection{Random Instance Generation for the Cardinality-Constrained Experiments}
\label{app:random_instances}

We describe the random instance-generation procedure used in Section~\ref{sec:cardinality_experiments}. For each problem size $N\in\{10,20,40,80,160,320,640\}$, we set $n_A=n_B=N$ and generate independent two-category \textsc{MarkovMC} instances directly in the model parameter space. To lighten notation in this subsection, we write $0$ for both no-purchase options $0_A$ and $0_B$ when the category is clear from context, so that $\NN_A^+=\NN_A\cup\{0\}$ and $\NN_B^+=\NN_B\cup\{0\}$. The no-purchase state is absorbing in both within-category Markov chains.

For category~$A$, the initial arrival vector is generated by drawing independent random variables $Y_i^A\sim U[0,1]$ for all $i\in\NN_A^+$ and normalizing:
\(
\psi_i^A
=
\frac{Y_i^A}{\sum_{\ell\in\NN_A^+}Y_\ell^A}, i\in\NN_A^+\).
The within-category transition matrix $\rho^A$ is generated row by row. For each product $i\in\NN_A$, we draw independent random variables $X_{ij}^A\sim U[0,1]$ for all $j\in\NN_A^+\setminus\{i\}$, set $X_{ii}^A=0$, and normalize:
\(
\rho_{ij}^A
=
\frac{X_{ij}^A}
{\sum_{\ell\in\NN_A^+\setminus\{i\}}X_{i\ell}^A},
j\in\NN_A^+\setminus\{i\}\).
Thus, real products have no self-loops and may transition directly to the no-purchase state. The no-purchase state is absorbing:
\(
\rho_{00}^A=1\), \(
\rho_{0j}^A=0, j\in\NN_A\).
The within-category transition matrix $\rho^B$ for category~$B$ is generated analogously.

We next generate the conditional arrival distribution into category~$B$. First, we draw a shared baseline arrival vector $q\in\Delta(\NN_B^+)$ by sampling independent random variables $Y_j^B\sim U[0,1]$ for all $j\in\NN_B^+$ and normalizing:
\(
q_j
=
\frac{Y_j^B}{\sum_{\ell\in\NN_B^+}Y_\ell^B}\),
\( j\in\NN_B^+\).
This vector represents the category-$B$ arrival distribution in the absence of cross-category dependence. At complementarity level $\beta=0$, every category-$A$ outcome induces the same category-$B$ arrival vector $q$, so the two categories are independent.

To introduce cross-category complementarity, we generate a complementarity map between real products in the two categories. As in the synthetic experiment in Section~\ref{sec:synthetic_experiments}, for each pair $(i,j)\in\NN_A\times\NN_B$, we draw $
c_{ij}\sim \mathrm{Bernoulli}(0.3)$,
where $c_{ij}=1$ indicates that product $j$ in category~$B$ is complementary to product $i$ in category~$A$. For each complementary pair, we draw an independent positive shock \(
\epsilon_{ij}\sim \mathrm{Exp}(1)\),
and set the effective shock to zero when $c_{ij}=0$. Given a complementarity-strength parameter $\beta\ge 0$, the conditional arrival probability of category-$B$ product $j\in\NN_B$ after category-$A$ product $i\in\NN_A$ is chosen is
\(
\lambda_{ij}^B(\beta)
=
\frac{q_j+\beta c_{ij}\epsilon_{ij}}
{q_0+\sum_{\ell\in\NN_B}
\left(q_\ell+\beta c_{i\ell}\epsilon_{i\ell}\right)}\).
The corresponding no-purchase probability is
\(
\lambda_{i0}^B(\beta)
=
\frac{q_0}
{q_0+\sum_{\ell\in\NN_B}
\left(q_\ell+\beta c_{i\ell}\epsilon_{i\ell}\right)}\).
If the customer makes no purchase in category~$A$, the category-$B$ arrival distribution remains the baseline distribution: \(
\lambda_{0j}^B(\beta)=q_j\), \(
 j\in\NN_B^+\).
Therefore, $\beta=0$ corresponds exactly to the independent two-category setting, while larger values of $\beta$ shift the downstream arrival distribution toward products in category~$B$ that are complementary to the realized purchase in category~$A$.

Product prices are generated independently from lognormal distributions. Specifically, for each $i\in\NN_A$ and $j\in\NN_B$, we draw
\(
r_i^A\sim \mathrm{Lognormal}(1.5,1.0)\), \(
r_j^B\sim \mathrm{Lognormal}(1.0,1.0)\),
independently across products and categories.

After generating the unconstrained two-category \textsc{MarkovMC} instance, we compute the unconstrained joint optimal assortment $(S_A^*,S_B^*)$. The cardinality limits are then sampled based on the sizes of these unconstrained optimal assortments, following the experimental protocol of \citet{desir2020constrained}. Specifically, we draw
\(
K_A
\sim
\mathrm{Unif}\left(
\left\{1,\ldots,\max\left\{1,\left\lceil |S_A^*|/2\right\rceil\right\}\right\}
\right)\),
and
\(
K_B
\sim
\mathrm{Unif}\left(
\left\{1,\ldots,\max\left\{1,\left\lceil |S_B^*|/2\right\rceil\right\}\right\}
\right)\).
The resulting cardinality-constrained instances are used to compare \textsc{AltOpt} with the numerical MILP benchmark when Gurobi certifies optimality to the stated relative tolerance within the time limit, and with the LP relaxation upper bound otherwise.

\subsection{MILP and LP Benchmarks and Computational Details}
\label{app:cardinality_computation}

\paragraph{Exact MILP formulation.}
Recall that $\NN_A^+=\NN_A\cup\{0_A\}$ and $\NN_B^+=\NN_B\cup\{0_B\}$, where the outside options are absorbing and \(r_{0_A}^A=r_{0_B}^B=0\), and let \(K_A\) and \(K_B\) denote the cardinality limits.

For \(i\in\NN_A^+\) and \(j\in\NN_B^+\), let
\(x_i^A\) and \(x_j^B\) indicate whether the corresponding options are
offered. The variable \(y_i^A\) is the probability that the category-\(A\)
chain is absorbed at \(i\), while \(z_i^A\) is its expected transient flow
through state \(i\). The variable \(w_{ij}\) is the joint probability of
final outcomes \((i,j)\), and \(u_{ij}\) is the probability-weighted
transient flow through state \(j\) in category \(B\), conditional on the
category-\(A\) outcome being \(i\).

The exact formulation is the following flow-based MILP:
\begin{align}
\max\quad
  &\sum_{i\in\NN_A^+}\sum_{j\in\NN_B^+}
    (r_i^A+r_j^B)w_{ij}
    \label{eq:cc-milp-objective}\\
\text{s.t.}\quad
  &\psi_i^A+
    \sum_{\ell\in\NN_A^+}\rho_{\ell i}^A z_\ell^A
    =y_i^A+z_i^A,
    &&i\in\NN_A^+,
    \label{eq:cc-flow-A}\\
  &y_i^A\leq x_i^A,
    &&i\in\NN_A^+,
    \label{eq:cc-link-A1}\\
  &z_i^A\leq M_i^A(1-x_i^A),
    &&i\in\NN_A^+,
    \label{eq:cc-link-A2}\\
  &\lambda_{ij}^B y_i^A+
    \sum_{\ell\in\NN_B^+}\rho_{\ell j}^B u_{i\ell}
    =w_{ij}+u_{ij},
    &&i\in\NN_A^+,\ j\in\NN_B^+,
    \label{eq:cc-flow-B}\\
  &w_{ij}\leq x_j^B,
    &&i\in\NN_A^+,\ j\in\NN_B^+,
    \label{eq:cc-link-B1}\\
  &w_{ij}\leq y_i^A,
    &&i\in\NN_A^+,\ j\in\NN_B^+,
    \label{eq:cc-link-B2}\\
  &u_{ij}\leq M_{ij}^B(1-x_j^B),
    &&i\in\NN_A^+,\ j\in\NN_B^+,
    \label{eq:cc-link-B3}\\
  &\sum_{i\in\NN_A}x_i^A\leq K_A,\qquad
    \sum_{j\in\NN_B}x_j^B\leq K_B,
    \label{eq:cc-cardinality}\\
  &x_{0_A}^A=x_{0_B}^B=1,
    \label{eq:cc-outside}\\
  &x_i^A,x_j^B\in\{0,1\},
    &&i\in\NN_A^+,\ j\in\NN_B^+,
    \label{eq:cc-integrality}\\
  &0\leq y_i^A\leq1,\quad z_i^A\geq0,
    &&i\in\NN_A^+,\\
  &0\leq w_{ij}\leq1,\quad u_{ij}\geq0,
    &&i\in\NN_A^+,\ j\in\NN_B^+.
\end{align}

The big-\(M\) constants are computed using the fundamental matrices of the
two transient Markov chains. Define
\(
Q^A=(\rho_{\ell k}^A)_{\ell,k\in\NN_A}\), \(
Q^B=(\rho_{\ell k}^B)_{\ell,k\in\NN_B},
\)
and let \(\boldsymbol{\psi}_{\NN_A}^A\) and \(\boldsymbol{\lambda}_{\NN_B}^{i}\) denote the restrictions of the corresponding arrival vectors to the transient states. The implementation uses
\[
M_i^A=
\begin{cases}
0, & i=0_A,\\[1mm]
\displaystyle
\max\left\{
\left[
(\boldsymbol{\psi}_{\NN_A}^A)^\top
(I-Q^A)^{-1}
\right]_i,0
\right\}+10^{-9},
& i\in\NN_A,
\end{cases}
\]
and, for every \(i\in\NN_A^+\),
\[
M_{ij}^B=
\begin{cases}
0, & j=0_B,\\[1mm]
\displaystyle
\max\left\{
\left[
(\boldsymbol{\lambda}_{\NN_B}^{i})^\top
(I-Q^B)^{-1}
\right]_j,0
\right\}+10^{-9},
& j\in\NN_B.
\end{cases}
\]

\paragraph{LP-relaxation benchmark.}
The LP benchmark retains the objective and all flow, linkage, and
cardinality constraints in
\eqref{eq:cc-milp-objective}--\eqref{eq:cc-outside}, but replaces
\eqref{eq:cc-integrality} with
\(
0\leq x_i^A\leq1\) and \(
0\leq x_j^B\leq1
\)
for
\(
i\in\NN_A^+,\ j\in\NN_B^+
\).
The fractional \(x\)-values are not interpreted as an assortment; only the
LP objective is used as an upper bound. Let \(Z^\star\) denote the exact
optimal value of the MILP. If \(R^{\textsc{AltOpt}}\) is the revenue of the
feasible \textsc{AltOpt} solution, then
\(
R^{\textsc{AltOpt}}
\leq Z^\star
\leq Z_{\mathrm{LP}}
\).
Thus,
\(
R^{\textsc{AltOpt}}/Z_{\mathrm{LP}}
\leq
R^{\textsc{AltOpt}}/Z^\star
\),
so the ratio against the LP bound is a conservative lower bound on the
unknown exact optimality ratio.

\paragraph{Performance measures.}
For instance \(h\), let \(\widehat Z_h^{\mathrm{MILP}}\) denote the incumbent
objective value returned by Gurobi. Define
\[
Z_h^{\mathrm{bench}}=
\begin{cases}
\widehat Z_h^{\mathrm{MILP}}, & \text{if Gurobi reports optimal termination at the stated tolerance},\\
Z_h^{\mathrm{LP}}, & \text{otherwise}.
\end{cases}
\]
We define the performance ratio and percentage gap as
\[
q_h=
\frac{R_h^{\textsc{AltOpt}}}{Z_h^{\mathrm{bench}}},
\qquad
g_h=100(1-q_h)
=100\frac{Z_h^{\mathrm{bench}}-R_h^{\textsc{AltOpt}}}
          {Z_h^{\mathrm{bench}}}.
\]
For each problem size, the reported average and worst-case ratios are
\(
\overline q=\frac{1}{50}\sum_{h=1}^{50}q_h\),
\(
q_{\min}=\min_{1\leq h\leq50}q_h.
\)

\paragraph{Complementarity parameter.}
For the instances in Table~\ref{tab:cardinality-altopt}, the
complementarity parameter is fixed at \(\beta=2.0\).
For each instance, the shared baseline arrival vector, complementarity map,
and exponential shocks are generated once, as described in
Appendix~\ref{app:random_instances}, and then used to construct
\(\boldsymbol{\lambda}(\beta)\).

\paragraph{\textsc{AltOpt} implementation.}
The implemented heuristic is single-start. It initializes
\(S_A^{(0)}=\NN_A\),
computes the first category-\(B\) update conditional on this full
category-\(A\) assortment, and then computes the category-\(A\) update.
Each block update uses Algorithm~3 of \citet{desir2020constrained} with
parameters \(\epsilon=0.1\) and \(\alpha=0.5\) (in the notation of that paper); this is a polynomial-time approximation algorithm for the APX-hard cardinality-constrained MC assortment problem, so block updates are approximate rather than exact best responses.
After each complete \(B\)--\(A\) update, the exact joint revenue of the
resulting assortment pair is evaluated, and the best pair encountered is
retained. The algorithm terminates when
\(S_A^{(t)}=S_A^{(t-1)}\)
or after \(T=50\) complete iterations, whichever occurs first. Revenue
improvements smaller than \(\eta=10^{-9}\) are not used to update the incumbent.

Because the block subproblems are solved approximately, the sequence of
iterates is not guaranteed to improve the joint objective monotonically.
Accordingly, we do not claim that the returned solution is a coordinate-wise
local optimum; the algorithm returns the best solution encountered before
the fixed-point or iteration-limit condition is reached.

\paragraph{Solver and computing environment.}
The MILP and LP formulations were implemented in Python using Gurobi. The computational environment was Python~3.11.8 and Gurobi~11.0.0 on an
Apple M2 processor with eight cores and 16~GB of memory.
The MILP relative optimality tolerance was set to
\(\texttt{MIPGap}=10^{-4}\),
and the per-instance time limit was \(7{,}200\) seconds. Instances for which
Gurobi reported optimal termination at this tolerance were used as MILP
benchmarks.

\section{Beyond Exact Arrival-Invariance}
\label{sec:approx-invariance}

Theorem~\ref{thm:opt} uses a downstream assortment that is optimal for every arrival distribution. We now extend the adjusted-price argument to assortments that approximate optimal downstream revenue uniformly over distributions of upstream choices. We keep the sequential structure but allow general conditional downstream choice models $\phi_{B|i}$. Product revenues are nonnegative, and the outside options have zero revenue. For a distribution $\mathbf p\in\Delta(\NN_A^+)$ over upstream choices, let $
R_i(S_B)=\sum_{j\in S_B}r_j^B\phi_{B|i}(j,S_B)$, which is the downstream revenue conditional on choice $i$.

\begin{theorem}[Approximate Invariance]
\label{thm:approx-decoupling}
Let $\alpha,\beta\in(0,1]$. Suppose there exists an assortment $\tilde{S}_B$, independent of $\mathbf p$, such that, for every $\mathbf p\in\Delta(\NN_A^+)$,
\begin{equation}\label{eq:approx-arrival-invariance}
\sum_{i\in\NN_A^+}p_iR_i(\tilde{S}_B)\ge\alpha\max_{S_B\subseteq\NN_B} \sum_{i\in\NN_A^+}p_iR_i(S_B).
\end{equation}
If the assortment optimization problem for category $A$ admits a $\beta$-approximation, then Algorithm~\ref{alg:mc_assortment}, with $S_B^{\ALG}$ set to $\tilde{S}_B$ and the category-$A$ step solved by this $\beta$-approximation, gives an $\alpha\beta$-approximation.
\end{theorem}

\proof{Proof.}
Let $(S_A^*,S_B^*)$ be an optimal solution, and let $\mathbf p^*$ be the distribution of upstream outcomes induced by $S_A^*$. By \eqref{eq:approx-arrival-invariance} and nonnegative upstream revenues, $\tilde{S}_B$ satisfies
\begin{equation}\label{eq:candidate-downstream-bound}
R(S_A^*,\tilde{S}_B)
=\sum_{i\in\NN_A^+}p_i^*r_i^A+\sum_{i\in\NN_A^+}p^*_iR_i(\tilde{S}_B)
\ge\alpha R(S_A^*,S_B^*).
\end{equation}
Let $S_B^{\ALG}=\tilde{S}_B$, and let $S_A^{\ALG}$ denote the output of Algorithm \ref{alg:mc_assortment}. Since the assortment for $A$ is a $\beta$-approximate solution given $S_B^{\ALG}$, we have
$R(S_A^{\ALG},S_B^{\ALG})
\ge\beta R(S_A^*,\tilde{S}_B)
\ge\alpha\beta R(S_A^*,S_B^*)$.
\Halmos
\endproof

The original MC specification satisfies the condition with a single invariant optimal assortment and $\alpha=1$. The following corollary shows the extension works when $\phi_B$ follows the independent-demand/MNL mixture of \citet{cao2023revenue}, which is not a special case of the MC model.

\begin{corollary}[Mixture of Independent Demand and MNL]
\label{cor:cao-mixture-half}
Suppose the choice model for category $B$ follows the
independent-demand/MNL mixture of \citet{cao2023revenue}:
\(
\phi_{B\mid i}(j,S_B)
=
q\theta_{i,j}
+
(1-q)\frac{v_j}{1+\sum_{k\in S_B}v_k},
i\in\NN_A^+,j\in S_B\),
where \(q\in[0,1]\), the weights \(v_j\geq 0\) for $j\in \NN_B$, and
\(\boldsymbol{\theta}_i\in\Delta(\NN_B^+)\) may depend on \(i\).
Suppose further that the
category-\(A\) assortment optimization problem
admits a polynomial-time \(\beta\)-approximation, where
\(\beta\in(0,1]\). Then the joint assortment optimization problem
admits a polynomial-time \(\beta/2\)-approximation.
\end{corollary}

\proof{Proof of Corollary~\ref{cor:cao-mixture-half}.}
For any distribution
\(\mathbf{p}\in\Delta(\NN_A^+)\) over upstream choices, define
\(
\bar{\theta}_j
=
\sum_{i\in\NN_A^+}p_i\theta_{i,j}.
\)
The resulting downstream revenue is
\(
G_{\mathbf{p}}(S_B)
=
q\sum_{j\in S_B}r_j^B\bar{\theta}_j
+
(1-q)
\frac{\sum_{j\in S_B}r_j^Bv_j}
     {1+\sum_{j\in S_B}v_j}\).
Denote the two nonnegative terms by
\(I_{\mathbf{p}}(S_B)\) and \(M(S_B)\), respectively. Let
\(
S^{\mathrm{ID}}
=
\{j\in\NN_B:r_j^B>0\},
\)
which maximizes \(I_{\mathbf{p}}\) for every \(\mathbf{p}\), and let
\(S^{\mathrm{MNL}}\) maximize \(M\). Both assortments are independent
of \(\mathbf{p}\) and can be computed in polynomial time.

Fix \(\mathbf{p}\), and let \(S_B^*\) maximize
\(G_{\mathbf{p}}\). Since each candidate maximizes one of the two
components,
\begin{equation}
\label{eq:candidate-mixture-half}
\begin{aligned}
\max\bigl\{
G_{\mathbf{p}}(S^{\mathrm{ID}}),
G_{\mathbf{p}}(S^{\mathrm{MNL}})
\bigr\}
&\geq
\max\bigl\{
I_{\mathbf{p}}(S_B^*),
M(S_B^*)
\bigr\} \geq
\frac{1}{2}
\bigl(
I_{\mathbf{p}}(S_B^*)+M(S_B^*)
\bigr) =
\frac{1}{2}G_{\mathbf{p}}(S_B^*).
\end{aligned}
\end{equation}

We now apply this bound to the joint problem. For each
\(C\in\{S^{\mathrm{ID}},S^{\mathrm{MNL}}\}\), fix \(S_B=C\), apply the
\(\beta\)-approximation algorithm to the resulting category-\(A\)
problem with adjusted revenues \(r_i^A+R_i(C)\), and denote its output
by \(\widehat{S}_A^C\). Return the candidate pair with the larger joint
revenue. Let \((S_A^*,S_B^*)\) be an optimal joint assortment, and let
\(\mathbf{p}^*\) be the distribution of upstream choices induced by
\(S_A^*\). The returned solution satisfies
\begin{align*}
R(S_A^{\mathrm{ALG}},S_B^{\mathrm{ALG}})
&\geq
\beta
\max_{C\in\{S^{\mathrm{ID}},S^{\mathrm{MNL}}\}}
R(S_A^*,C) =
\beta\biggl(
\sum_{i\in\NN_A^+}p_i^*r_i^A
+
\max_{C\in\{S^{\mathrm{ID}},S^{\mathrm{MNL}}\}}
G_{\mathbf{p}^*}(C)
\biggr) \\
&\geq
\beta\left(
\sum_{i\in\NN_A^+}p_i^*r_i^A
+
\frac{1}{2}G_{\mathbf{p}^*}(S_B^*)
\right) \geq
\frac{\beta}{2}R(S_A^*,S_B^*),
\end{align*}
where the final inequality follows from the nonnegativity of upstream
revenues. Thus, the procedure is a polynomial-time
\(\beta/2\)-approximation.
\Halmos
\endproof

\end{APPENDICES}
\end{document}